\documentclass[11pt]{article}
\usepackage{amsmath}
\usepackage{amssymb}
\usepackage{amscd}
\usepackage{amsthm}

\def\proof{\vspace{2ex}\noindent{\bf Proof.} }

\def\endproof{\relax\ifmmode\expandafter\endproofmath\else
  \unskip\nobreak\hfil\penalty50\hskip.75em\hbox{}\nobreak\hfil\bull
  {\parfillskip=0pt \finalhyphendemerits=0 \bigbreak}\fi}
\def\endproofmath$${\eqno\bull$$\bigbreak}
\def\bull{\vbox{\hrule\hbox{\vrule\kern3pt\vbox{\kern6pt}\kern3pt\vrule}\hrule}}
\newtheorem{theorem}{Theorem}[subsection]

\newtheorem{proposition}[theorem]{Proposition}
\newtheorem{lemma}[theorem]{Lemma}

\newtheorem{corollary}[theorem]{Corollary}

\newtheorem{D}[theorem]{Definition}
\newenvironment{defn}{\begin{D} \rm }{\end{D}}

\newtheorem{R}[theorem]{Remark}
\newenvironment{remark}{\begin{R}\rm }{\end{R}}

\def\Zee{\mathbb{Z}}

\def\Cee{\mathbb{C}}
\def\Aff{\mathbb{A}}
\def\Pee{\mathbb{P}}
\def\WP{\mathbb{WP}}

\def\scrO{\mathcal{O}}
\def\ov{\overline}
\def\spcheck{^{\vee}}
\def\frak{\mathfrak}
\def\Pic{\operatorname{Pic}}
\def\Ker{\operatorname{Ker}}
\def\Coker{\operatorname{Coker}}
\def\Sym{\operatorname{Sym}}
\def\Hom{\operatorname{Hom}}
\def\End{\operatorname{End}}
\def\Aut{\operatorname{Aut}}
\def\Ext{\operatorname{Ext}}
\def\Id{\operatorname{Id}}

\def\Spec{\operatorname{Spec}}
\def\ad{\operatorname{ad}}
\def\Ad{\operatorname{Ad}}
\def\Lie{\operatorname{Lie}}

\title{Holomorphic Principal Bundles Over Elliptic Curves III:\\
Singular Curves and Fibrations}
\author{Robert Friedman\thanks{The first author was partially
    supported by NSF grant DMS-99-70437.}
\  and John W. Morgan\thanks{The second author was partially supported
   by NSF grant DMS-97-04507.}}

\begin{document}

\maketitle

\section*{Introduction}

Let $G$ be a reductive complex linear algebraic group, which in this paper for
simplicity we shall always assume  to be simple and simply connected, and let
$\Lie (G) =\frak g$. This paper is the third in a series   concerned with the
study of holomorphic principal
$G$-bundles over elliptic fibrations. The first paper in the series studied the
moduli space $\mathcal{M}$ of $G$-bundles on a smooth elliptic curve $E$ by means
of flat connections, or equivalently by finding a holomorphic reduction of the
structure group to a Cartan subgroup $H$ of $G$. This approach, while in many
ways the most natural, is not adapted to finding universal $G$-bundles over
$E\times
\mathcal{M}$ or to dealing with the case of singular elliptic curves or
fibrations. In the second paper, we introduced a new method for studying
$G$-bundles, by finding a reduction of structure to  a certain maximal parabolic
subgroup $P$. In this construction, one sees naturally that the moduli space is a
weighted projective space, and that the weights are given as follows: let
$\alpha_1, \dots, \alpha_r$ be a set of simple roots for $(G,H)$,   let
$\widetilde
\alpha$ be the highest root, and set $\alpha_0 = -\widetilde \alpha$. Then
there is a unique linear relation $\sum _{i=0}^rg_i\alpha_i\spcheck =0$ with
$g_0=1$, and the weights of the weighted projective space are  the
positive integers
$g_i$, $0\leq i\leq r$. We also used this construction to exhibit universal
bundles  over 
$E\times\mathcal{M}$, at least
for an appropriate conformal form of $G$ and away from the singular points of
$\mathcal{M}$. In this paper, we shall continue our study by analyzing how this
construction behaves for singular elliptic curves and how it can be made in
families $\pi\colon Z\to B$ of Weierstrass cubics. One goal is to construct a
relative moduli space, at least if either $G$ is not of type $E_8$ or $\pi$ has
no cuspidal fibers, which is a weighted projective space bundle over
$B$. If $G$ is not of type $E_8$, then this bundle is the bundle of weighted
projective spaces associated to 
$$\scrO_B\oplus \mathcal{L}^{-d_1} \oplus \cdots \oplus \mathcal{L}^{-d_r},$$
where the $d_i$ are the Casimir weights of $G$, i.e.\ the numbers $m_i+1$,
where the $m_i$ are the exponents of $G$, and the line bundle $\mathcal{L}=
R^0\pi_*\omega_{Z/B}$ is the direct image of the relative dualizing sheaf.
Moreover,
$\Cee^*$ acts on
$\mathcal{L}^{-d_i}$ with weight $g_i$ for some ordering of the  simple roots,
and $\Cee^*$ acts on $\scrO_B$ with weight one. A very closely related result was
proved by Wirthm\"uller \cite{Wirt}, by somewhat \emph{ad hoc} methods. One goal
of this paper is to find a conceptual explanation for the appearance of the
integers $d_i$. A second goal is to understand the construction for cuspidal
or nodal curves, and to see that it gives  a
section of the adjoint quotient morphism into the set of regular elements of
$\frak g$ or $G$ as well as a compactification of the adjoint quotient.  Thus,
this part of our results should be viewed as both another approach to the
work of Kostant
\cite{Kostant} and Steinberg
\cite{Steinberg} on the existence of such sections, as well as a 
generalization, in the sense that we construct natural compactifications
of the adjoint quotient and show that the Steinberg section is in some sense
a deformation of the Kostant section, and that both may be deformed to the
parabolic construction on a smooth elliptic curve.

The basic idea behind the construction for singular curves and families is to
extend the parabolic construction to singular curves. Let us recall this
construction for a smooth elliptic curve $E$, with origin $p_0$.  Let $R$ be
the root system of the pair
$(G,H)$ and let
$\Delta$ be a set of simple roots. For each  $\alpha\in \Delta$, there is a
corresponding maximal parabolic subgroup $P$, with unipotent radical
$U$ and Levi factor $L$. We define $\alpha$ to be \textsl{special} if 
\begin{enumerate}
\item[(i)] The Dynkin diagram associated to $\Delta -\{\alpha\}$ is a union of
diagrams of $A$-type;
\item[(ii)] The simple root $\alpha$ meets each component of the Dynkin diagram
associated to $\Delta -\{\alpha\}$ at an end of the component;
\item[(iii)] The root $\alpha$ is a long root.
\end{enumerate}
If $G$ is of type $A_n$, then every root is special. In all other cases, there is
a unique special root. It is easy to check that the derived subgroup of $L$ is
isomorphic to a product of groups of $A$-type, and that $L$ is the subgroup
of a product of groups of the form $GL_{n_i}(\Cee)$ consisting of matrices of
equal determinant. Thus an $L$-bundle over $E$ is the same as a collection of
vector bundles with equal determinant, and hence there is a
unique  $L$-bundle $\eta_0$ which corresponds to a collection of stable
vector bundles of determinant $\scrO_E(-p_0)$. In this case, lifts of
$\eta_0$ to a
$P$-bundle are classified by the nonabelian cohomology set $H^1(E;
U(\eta_0))$, where $U(\eta_0)$ is the sheaf of unipotent groups associated to
the principal $L$-bundle $\eta_0$ and the action of $L$ on $U$. While this 
cohomology set is hard to analyze directly, one can look at the linearized
version $H^1(E; \frak u(\eta_0))$, where $\frak u$ is the Lie algebra of $U$. It
is easy to see that the $L$-module $\frak u$ is a direct sum of irreducible
representations $\frak u^k$ of
$L$ on which the center acts via positive weights, and thus $\frak u(\eta_0)$ is
a direct sum of semistable vector bundles $\frak u^k(\eta_0)$ over $E$ of 
negative degrees, which can be explicitly calculated. From this, one can
prove that the set 
$H^1(E; U(\eta_0))$ carries the structure of an affine space, on which 
$\Cee^*$ acts with the appropriate weights. The quotient $(H^1(E;
U(\eta_0))-\{0\})/\Cee^*$ is a weighted projective space, corresponding to
isomorphism classes of $P$-bundles whose reduction mod $U$ is isomorphic to
$\eta_0$. Of course, given a $P$-bundle, one can form the associated
$G$-bundle. One of the main results in
\cite{FMII} is that every S-equivalence class in $\mathcal{M}$ arises
uniquely in this way. In particular, up to S-equivalence, every semistable
holomorphic $G$-bundle over $E$ has a natural reduction of structure  to a
$P$-bundle with certain very special properties.

How does this picture change when we consider, instead of smooth elliptic curves
$E$, singular curves $C$ of arithmetic genus one? One fundamental difference is
that, while one can still define the notion of a stable or semistable vector
bundle over a singular curve $C$, there is no reasonable candidate at present
for a definition of a semistable principal $G$-bundle over $C$. Closely
related to this is the fact that the tensor product of two semistable vector
bundles on a singular curve can be unstable. At best, one can consider the
set of all
$G$-bundles on $C$ whose pullback to the normalization $\widetilde C$ is
semistable (and hence is trivial, in our case,  since $\widetilde C\cong
\Pee^1$). While a coarse moduli space of such bundles exists, it is never
compact. On the other hand, by the explicit description of  the group
$L$, there is a natural definition of the
$L$-bundle $\eta_0$ over $C$ and it is uniquely specified by requiring that it
correspond to $t$ stable vector bundles of  given determinant
$\scrO_C(-p_0)$, where $p_0$ is an origin for the group $C_{\rm reg}$ of
smooth points of $C$. We can then form the set 
$H^1(C; U(\eta_0))$ and its linearized version $H^1(C; \frak u(\eta_0))$.
However, the vector bundles $\frak u^k(\eta_0)$ are typically no longer
semistable. Nonetheless, after a long and somewhat painful calculation, these
vector bundles can be directly tabulated and their cohomology can be shown to
behave just like that of a semistable bundle of the same degree, except in
one case where
$C$ is cuspidal and $G$ is of type $E_8$. Unfortunately, we have no conceptual
reason why the calculation fails in this case, nor why it gives the right
answer in all other cases. In a future paper \cite{FMIV}, we will discuss the
corresponding picture for the moduli space of del Pezzo surfaces of degree
one, or equivalently rational elliptic surfaces, and give  a geometric
explanation for the anomalous nature of $E_8$.

Once we have constructed the bundle $\eta_0$, or more precisely a conformal
variant of $\eta_0$ in the case of a fibration $\pi\colon Z\to B$, then we can
make the parabolic construction for the family $Z$ as long as the cohomology
calculations described above give the expected answer, and hence whenever $G$ is
not of type $E_8$ or $\pi$ has no cuspidal fibers. To further analyze the result,
we consider the  universal family $\mathcal{E}
\subseteq \Pee^2 \times \Aff^2$ defined by
$$\mathcal{E} = \{([x,y,z],g_2, g_3): y^2z = x^3 + g_2xz^2 + g_3z^3\}.$$
Of crucial importance is the fact that the base $\Aff^2$ has a $\Cee^*$-action
which lifts to $\mathcal{E}$, with a unique fixed point at $0$ corresponding to a
cuspidal curve. The main idea now is to use the $\Cee^*$-action to localize the
picture around the cuspidal curve, and to relate the weights there to the
Chevalley generators of the graded ring $(\Sym^*\frak h^*)^W$ of Weyl-invariant
polynomial functions on a Cartan subalgebra $\frak h =\Lie(H)$.

The organization of this paper is as follows. In Sections 1, 2, and 3 we
collect many preliminary results, most of which seem to be well-known but for
which we could not find adequate references. In Section 1, it is convenient
to consider the moduli stack $\mathbf{E}$ of elliptic fibrations with a
section. We use very little of the general theory of stacks. Instead the main
point is to describe the group of line bundles over $\mathbf{E}$ and their
sections. An elementary argument shows that every vector bundle over
$\mathbf{E}$ is a direct sum of line bundles. Finally, we show that a
conformal form of the stable $L$-bundle $\eta_0$ defined above exists over
$\mathbf{E}$.   Section 2 contains elementary facts concerning vector bundles
over nodal and cuspidal curves. These facts are used to make the necessary
cohomology calculations and to analyze the stability of vector bundles
associated to certain $G$-bundles. Related results about vector bundles over
singular curves may be found in \cite{BDG} and the references therein. In
Section 3, we describe principal
$G$-bundles over nodal and cuspidal curves in terms of $G$-bundles on the
normalization. There is a also a very brief discussion of a partial moduli
space for
$G$-bundles over singular curves.

We begin Section 4 with the necessary cohomology calculations for the
parabolic construction. Once these are established, it follows from the theory
of
\cite{FMII} that, except in the case of $E_8$ and cuspidal fibers, the
constructions fit together in families to give a bundle of weighted
projective spaces. Turning to the nodal or cuspidal case, we show that the
set of bundles inside the weighted projective space which become trivial on
the normalization is an affine space, which we then  identify with the
adjoint quotient of $G$ or of
$\frak g$. Finally, by linearizing the
$\Cee^*$-action for the family $\mathcal{E}$, we show that the relative moduli
space is the quotient of the vector bundle associated to $\scrO_B\oplus
\mathcal{L}^{-d_1} \oplus \cdots \oplus \mathcal{L}^{-d_r}$ by a diagonal
$\Cee^*$-action. 

In Section 5, we prove a very general result on families of weighted
projective bundles which shows that our compactification is unique in an
appropriate sense. We use this result to compare our construction to that of
Wirthm\"uller, and will also use it in a future paper \cite{FMIV} to relate
our construction to ones that can be made  for del Pezzo surfaces.
In Section 6, we study the way that the parabolic construction behaves under
inclusions of subgroups $G'\subseteq G$. Similar methods give more insight into
the case of $E_8$. In particular, we are able to construct a bundle of
$9$-dimensional weighted projective spaces and an action of a one-dimensional
unipotent group   on this bundle, such that the quotient, if it existed,
would be the correct moduli space. However, the quotient does not exist in any
reasonable sense if there are cuspidal fibers. Still, given a cuspidal
curve, there is an eight-dimensional affine space inside of the corresponding
nine-dimensional weighted projective space which is a slice for the action of
the unipotent group on the open subset of bundles which become trivial on the
normalization and which plays the role of the Kostant section. The situation
is very closely related to the correspondence between del Pezzo surfaces of
degree one and rational elliptic surfaces, as we will show in a future paper.
Finally, in Section 7, in the case of a singular curve $C$ of arithmetic
genus one, we analyze when the vector bundle
$V$ associated to a $G$-bundle
$\xi$ coming from the parabolic construction and an irreducible representation
$\rho\colon G \to GL_N(\Cee)$ is unstable.  If $\xi$ pulls back to the
trivial bundle on the normalization of $C$, then $V$ is always semistable.
Otherwise, there are only a very few representations $\rho$ for which $V$
might possibly be semistable. We classify all such representations and
identify those bundles $\xi$ arising  from the parabolic construction for
which $V$ is unstable. For $G$ of type $E_6$, $E_7$, or
$F_4$, there exist bundles arising  from the parabolic construction such that,
for every irreducible representation
$\rho$, the associated vector bundle $V$ is unstable. For $G$ of type $E_8$,
the situation is even more striking: for every $G$-bundle $\xi$ such that the
pullback of $\xi$ to the normalization is not trivial and every irreducible 
representation
$\rho$, the associated vector bundle $V$ is unstable.

There are many remaining open questions. One of the deepest is the problem of
finding an intrinsic definition of semistability for $G$-bundles on a singular
curve, and of a generalized form of S-equivalence, which would be broad enough to
include those bundles coming from the parabolic construction. We emphasize
that, even though the parabolic construction gives a very \emph{ad hoc}
construction of a compact moduli space, the bundles which arise are quite
natural from many points of view. For example, either they pull back to the
trivial bundle on the normalization, or they are in some sense generic
among bundles whose pullback to the normalization is nontrivial. One feature
of curves of arithmetic genus one, which is already apparent in the study of
vector bundles, is that for $G$ semisimple the moduli space of
$G$-bundles can be compactified by using only bundles. For example, in the
case of
$SL_n$-bundles, we do not have to consider   torsion free, non-locally free
sheaves in order to compactify the moduli space. However, this result fails
in higher genus. Another open problem is to find an explicit link between
Wirthm\"uller's toric compactification of $\Cee^*\otimes \Lambda=H$ and that
given by the parabolic construction.

Finally, we would like to thank Ed Witten for many stimulating discussions
concerning the subject of this paper and its predecessors, and Titus
Teodorescu for several conversations about some of the material in Section 2.

\section{Bundles over the moduli stack of Weierstrass fibrations}

Recall that a \textsl{Weierstrass cubic} $C$ is a reduced irreducible curve
$C$ of arithmetic genus one, together with a point $p_0$ on the smooth locus
of $C$. Every such curve is either a smooth
elliptic curve or is isomorphic to an irreducible plane cubic with either a
node or a cusp. We shall refer to such curves (always understood to have
arithmetic genus one) as
\textsl{nodal} or \textsl{cuspidal} respectively. Each such curve is embedded
as a plane cubic curve by the linear system $|3p_0|$ and thus has a
homogeneous equation $y^2z = x^3 + g_2xz^2 + g_3z^3$, where $g_2$ and $g_3$
are specified up to the action of $\Cee^*$ defined by $\mu\cdot (g_2, g_3) =
(\mu^4g_2, \mu^6g_3)$. We will refer to such an equation as a
\textsl{Weierstrass model}.

Our goal in
this section is to study the moduli stack of Weierstrass cubics and certain
bundles over this stack. We begin by defining the stack and showing that it
is a quotient stack. Next, we classify line bundles and vector bundles
over the stack. Finally, we turn to the existence of certain vector bundles
over the total space of the stack, which are the building blocks for
constructing unstable $G$-bundles (or rather conformal $G$-bundles) over this
total space.

\subsection{Generalities on the moduli stack}

\begin{defn}\label{def1} A \textsl{Weierstrass fibration} consists of a
triple $(Z, \pi, \sigma)$, where $Z$ is a scheme, $\pi\colon Z \to B$ is  a
flat proper morphism such that every fiber of
$\pi$ is a reduced irreducible curve of arithmetic genus one, and
$\sigma\colon B \to Z$ is a section
 whose image (which we shall also write as $\sigma$) is
contained in the smooth locus of $\pi$. Every fiber is a Weierstrass cubic. As 
is well-known, there is a one-to-one correspondence between  isomorphism
classes of Weierstrass fibrations over $B$ and line bundles
$\mathcal{L}$ over $B$, together with sections $G_2\in H^0(B;
\mathcal{L}^{\otimes 4})$ and $G_3\in H^0(B; \mathcal{L}^{\otimes 6})$,
modulo the action of $H^0(B; \scrO_B^*)$ on $H^0(B; \mathcal{L}^{\otimes
4})\times  H^0(B; \mathcal{L}^{\otimes 6})$ defined by $f
\cdot (G_2, G_3) = (f^4G_2, f^6G_3)$. Here, $\mathcal{L} ^{-1} =
R^1\pi_*\scrO_Z$, which is isomorphic to the normal bundle of $\sigma$ in
$Z$ under the isomorphism $\pi\colon \sigma \to B$, and by relative duality
$\mathcal{L}
\cong R^0\pi_*\omega_{Z/B}$. Moreover,
$Z$ is the subscheme of $\Pee(\mathcal{L}^2\oplus
\mathcal{L}^3\oplus \scrO_B)$ defined by the
section $y^2z - (x^3 + G_2xz^2 + G_3z^3)$ of $\Sym^3(\mathcal{L}^2\oplus
\mathcal{L}^3\oplus \scrO_B)\spcheck \otimes \mathcal{L}^{6}$. We shall also
refer to this subscheme of $\Pee(\mathcal{L}^2\oplus
\mathcal{L}^3\oplus \scrO_B)$ as a \textsl{Weierstrass model}. 

We
define the
\textsl{moduli stack} $\mathbf{E}$ as the category whose objects are triples 
$(Z,
\pi, \sigma)$ as above, and such that the set of morphisms from
$(Z',\pi',\sigma')$ to
$(Z,\pi, \sigma)$ is the set of Cartesian diagrams
$$\begin{CD} Z' @>{g}>> Z\\ @V{\pi'}VV @VV{\pi}V\\ B' @>{f}>> B,
\end{CD}$$ such that $g^{-1}(\sigma) = \sigma'$, viewing $\sigma$ as a
subvariety of $Z$, or equivalently $g\circ \sigma ' =\sigma \circ f$. Thus
the set of morphisms from $(Z',\pi',\sigma')$ to $(Z,\pi, \sigma)$ covering a
given $f\colon B' \to B$ is a principal homogeneous space over the group of
automorphisms of $Z'$ covering the identity on $B'$ and preserving $\sigma'$.
\end{defn}

The stack $\mathbf{E}$ is not representable. However, there is an almost
universal family defined as follows. Let $\Aff^2 = \Spec \Cee[g_2, g_3]$ be
the affine plane with coordinates $g_2, g_3$. Define $\mathcal{E}
\subseteq \Pee^2 \times \Aff^2$ via
$$\mathcal{E} = \{([x,y,z],g_2, g_3): y^2z = x^3 + g_2xz^2 + g_3z^3\}.$$ Then
the induced morphism $\pi\colon\mathcal{E} \to \Aff^2$ and the section
$\sigma(g_2,g_3) =  ([0,1,0], g_2, g_3)$ exhibit
$\mathcal{E}$ as a Weierstrass fibration over $\Aff^2$. There is a cuspidal
fiber over $(0,0)$ and a nodal fiber over
$(-3a^2, -2a^3)$ for
$a\neq 0$. All other fibers are smooth. 

\begin{defn}\label{action} Let $\lambda\colon \Cee^*\times \mathcal{E} \to
\mathcal{E}$ be the  $\Cee^*$-action on
$\mathcal{E}$ defined by
$$\lambda(\mu, ([x,y,z],g_2, g_3)) = \mu\cdot ([x,y,z],g_2, g_3) = ([\mu^2x,
\mu^3y, z], \mu^4g_2,
\mu^6g_3).$$ This action covers the $\Cee^*$-action $\lambda_0$ on $\Aff^2$
given by
$\lambda_0(\mu , (g_2, g_3))=\mu\cdot(g_2,g_3) = (\mu^4g_2, \mu^6g_3)$. The
section
$\sigma$ is invariant under the
$\Cee^*$-action, in the sense that $\lambda^*\sigma =
\pi_2^*\sigma$, and hence $\mu^*\sigma =\sigma$ for all $\mu\in \Cee^*$.
\end{defn}

 The differential
$dx/y$ generates the dualizing sheaf of every fiber. Thus
$R^0\pi_*\omega_{\mathcal{E}/ \Aff^2} =\scrO_{\Aff^2}\cdot [dx/y]$ is the
trivial line bundle over $\Aff^2$, and $\mu\in \Cee^*$ acts on this line
bundle by sending the everywhere generating section $dx/y$ to $\mu^{-1}dx/y$.
It follows from relative duality that the action of $\mu$ on
$R^1\pi_*\scrO_{\mathcal{E}}=(R^0\pi_*\omega_{\mathcal{E}/ \Aff^2})\spcheck$
sends an everywhere generating section $s$ to
$\mu\cdot s$. Using the section $dx/y$ to trivialize the bundle
$R^0\pi_*\omega_{\mathcal{E}/ \Aff^2}$, an equivalent formulation is that the
$\Cee^*$-linearization on $R^0\pi_*\omega_{\mathcal{E}/ \Aff^2}$ is
equivalent to the linearization on the trivial bundle $\Aff^2 \times
\Cee$ defined by
$\mu\cdot (g_2, g_3, z) = (\mu^4g_2, \mu^6g_3, \mu\cdot z)$, and the
$\Cee^*$-invariant sections of $R^0\pi_*\omega_{\mathcal{E}/ \Aff^2}$
correspond to the functions $f\in \Cee[g_2, g_3]$ such that $f(\mu^4g_2,
\mu^6g_3) = \mu\cdot f(g_2, g_3)$ (so that there are no invariant sections)
and similarly for Zariski open subsets of
$\Aff^2$.

\begin{theorem}\label{Weiermodel} There is a one-to-one correspondence between
\begin{enumerate}
\item[\rm (i)] Isomorphism classes of Weierstrass fibrations over $B$;
\item[\rm (ii)] Principal $\Cee^*$-bundles $B^* \to B$, together with
$\Cee^*$-equivariant morphisms $F\colon B^* \to \Aff^2$, modulo the natural
action of the automorphism group $H^0(B; \scrO_B^*)$ of the principal bundle
$B^*$.
\end{enumerate} If $B^*$  and $F$ correspond to $\pi\colon Z \to B$, then the
line bundle corresponding to
$B^*$ is
$ R^1\pi_*\scrO_Z=\mathcal{L}^{-1}$.
This correspondence is natural under pullback. Finally, if $F\colon B^* \to
\Aff^2$ corresponds to $Z$ and $Z^*$ is the pullback of $Z$ to $B^*$, then
$Z^*$ is $\Cee^*$-equivariantly isomorphic to $F^*\mathcal{E}$ and $Z$ is
isomorphic to $F^*\mathcal{E}/\Cee^*$, which exists as a scheme since $B^*\to
B$ is locally trivial in the Zariski topology.
\end{theorem}
\begin{proof} As we have seen above, there is a one-to-one correspondence
between  isomorphism classes of Weierstrass fibrations over $B$ and line
bundles
$\mathcal{L}$ over $B$, together with sections $G_2\in H^0(B;
\mathcal{L}^{\otimes 4})$ and $G_3\in H^0(B; \mathcal{L}^{\otimes 6})$,
modulo the action of $H^0(B; \scrO_B^*)$ on $H^0(B; \mathcal{L}^{\otimes
4})\times  H^0(B; \mathcal{L}^{\otimes 6})$ defined by $\lambda
\cdot (G_2, G_3) = (\lambda^4G_2, \lambda^6G_3)$. Such a pair of sections is
equivalent to a $\Cee^*$-equivariant morphism $F\colon B^*\to \Aff^2$, where
$B^*$ is the total space of the line bundle $\mathcal{L}^{-1}$. In one
direction, the morphism
$F=(f_1, f_2)$ is defined via
$$f_1(s) = G_2(s^{\otimes 4}) \in \Cee;\ \ f_2(s) = G_3(s^{\otimes 6}) \in
\Cee.$$  Conversely, if say $f \colon B^* \to \Cee$ satisfies $f(\mu\cdot s) =
\mu^af(s)$, then $f$ corresponds to a section of the line bundle $B^*\times
_{\Cee^*}\Cee$, where $\Cee^*$ acts on $\Cee$ via the character $\mu \mapsto
\mu^{-a}$, and thus $f$ is a section of the line bundle $B^*\times
_{\Cee^*}\Cee = \mathcal{L}^a$. Hence a $\Cee^*$-equivariant morphism
$F\colon B^*\to
\Aff^2$ defines sections $G_2\in H^0(B; \mathcal{L}^{\otimes 4})$ and $G_3\in
H^0(B; \mathcal{L}^{\otimes 6})$, and these two constructions are easily
checked to be inverses of each other, and to be compatible with the natural
action of
$H^0(B; \scrO_B^*)$. The remaining statements of the theorem are clear.
\end{proof}

Recall that the quotient stack $[\Aff^2/\Cee^*]$ is the functor whose objects
are principal $\Cee^*$-bundles $p\colon B^*\to B$, together with
$\Cee^*$-equivariant morphisms $F\colon B^*\to \Aff^2$, and whose morphisms
are Cartesian diagrams compatible with the morphisms to $\Aff^2$. We shall
denote an object of $[\Aff^2/\Cee^*]$ as a triple $(B^*, p, F)$. One way to
rephrase the preceding theorem, which makes more precise the sense in which
the family $\pi\colon\mathcal{E} \to \Aff^2$ is almost universal, is as
follows:

\begin{theorem} The functor $\mathbf{w}\colon  [\Aff^2/\Cee^*] \to
\mathbf{E}$ which assigns to a triple $(B^*, p, F)$ the Weierstrass model 
$F^*\mathcal{E}/\Cee^*$ is an equivalence of categories, i.e.\ the stacks
$\mathbf{E}$ and $[\Aff^2/\Cee^*]$ are isomorphic.
\end{theorem}
\begin{proof} Every   $(Z,\pi,  \sigma)$ is isomorphic to a  Weierstrass
model, i.e.\ to an object of the form $\mathbf{w}(o)$ for some object $o$ of
$[\Aff^2/\Cee^*]$. Thus, it is enough to check that the group of automorphisms
of a triple
$(B^*, p, F)$ covering the identity on $B$ is isomorphic to the automorphism
group of the corresponding Weierstrass model.  In both cases, this group is
isomorphic to the group of $\lambda \in H^0(B; \scrO_B^*)$ such that
$\lambda^4G_2=G_2$ and $\lambda^6G_3=G_3$, and hence the groups are the same.
\end{proof}

\subsection{Vector bundles over the moduli stack}

\begin{defn} A \textsl{scheme $\mathbf{X} \to
\mathbf{E}$ over the moduli stack}  consists of the following: for each
triple $(Z,B,\sigma)$, we are given a scheme
$X_Z\to B$, and, for each morphism $(g,f) \in \Hom ((Z',\pi', \sigma'),
(Z,\pi,
\sigma))$, i.e.\ a commutative diagram as in Definition~\ref{def1}, we are
given a 
$B'$-isomorphism
$\varphi_{(g,f)}\colon X_{Z'} \to f^*X_Z$, satisfying:
\begin{enumerate}
\item $\varphi _{(\Id,\Id)} = \Id $;
\item $\varphi_{(g_2\circ g_1, f_2 \circ f_1)} = \varphi_{(g_2,f_2)}\circ
f_2^*\varphi_{(g_1,f_1)}$ if the domain of $(g_2,f_2)$ is the range of
$(g_1,f_1)$.
\end{enumerate} We define a \textsl{line bundle over the moduli stack} or a
\textsl{coherent sheaf  over the moduli stack} similarly. The group of all
isomorphism classes of line bundles over $\mathbf{E}$ will be denoted $\Pic
\mathbf{E}$.

We can also define a  \textsl{scheme over the total space of the moduli stack}
$\mathbf{Y} \to \mathrm{Tot}(\mathbf{E})$ as follows: for each
$(Z,B,\sigma)$, we are given a scheme
$Y_Z\to Z$, and, for each morphism $(g,f) \in \Hom ((Z',\pi', \sigma'),
(Z,\pi,
\sigma))$,  we are given  a 
$Z'$-isomorphism
$\psi_{(g,f)}\colon Y_{Z'} \to g^*Y_Z$, satisfying:
\begin{enumerate}
\item $\psi _{(\Id,\Id)} = \Id $;
\item $\psi_{(g_2\circ g_1, f_2 \circ f_1)} = \psi_{(g_2,f_2)}\circ
g_2^*\psi_{(g_1,f_1)}$ if the domain of $(g_2,f_2)$ is the range of
$(g_1,f_1)$.
\end{enumerate} Line bundles, coherent sheaves and principal bundles over the
total space of the moduli stack are defined similarly.
\end{defn}

Since $\mathbf{E}$ is isomorphic to $[\Aff^2/\Cee^*]$, we can identify
schemes over $\mathbf{E}$ with schemes over $[\Aff^2/\Cee^*]$. To see this
concretely, let
$\mathbf{X}
\to\mathbf{E}$ be a scheme over the moduli stack. Then in particular, we can
define the scheme
$X_{\mathcal{E}}
\to
\Aff^2$, and by functoriality it has a natural action of $\Cee^*$ lifting the
action on
$\Aff^2$. Given a triple $(Z,B,\sigma)$, with corresponding morphism
$F\colon B^*\to
\Aff^2$, the scheme $F^*X_{\mathcal{E}}$ has a $\Cee^*$-action covering that
on $B^*$, coming from that on $X_{\mathcal{E}}$, and this agrees with the
action viewing $F^*X_{\mathcal{E}}$ as the pullback of $X_Z$ to $B^*$.
Conversely, given a scheme $X_{\mathcal{E}} \to \Aff^2$ together with an
action of $\Cee^*$ lifting the action on $\Aff^2$, there is an induced scheme
$F^*X_{\mathcal{E}}$ over $B^*$. Since $B^*\to B$ is a Zariski locally
trivial fibration, the quotient $F^*X_{\mathcal{E}}/\Cee^*$ exists as a
scheme $X_Z$, and it is easy to see that the $X_Z$ fit together to give a
scheme over $\mathbf{E}$.  Similarly, schemes $\mathbf{Y} \to
\mathrm{Tot}(\mathbf{E})$ may be identified with schemes
$Y_{\mathcal{E}} \to  \mathcal{E}$ together with a lifting of the
$\Cee^*$-action on $\mathcal{E}$. Thus we see:

\begin{theorem} There is a one-to-one correspondence between schemes
$\mathbf{X} \to\mathbf{E}$ and schemes $X_{\mathcal{E}} \to \Aff^2$ together
with an action of $\Cee^*$ on $X_{\mathcal{E}}$ lifting the action on
$\Aff^2$. A similar statement holds for coherent sheaves over $\mathbf{E}$.
Likewise, there is a one-to-one correspondence between  schemes $\mathbf{Y}
\to
\mathrm{Tot}(\mathbf{E})$ and schemes
$Y_{\mathcal{E}} \to  \mathcal{E}$ together with a lifting of the
$\Cee^*$-action on $\mathcal{E}$.
\qed
\end{theorem}

For example, let us consider line bundles $\mathbf{L}$ over the moduli stack.
These correspond to line bundles over $\Aff^2$ together with a
$\Cee^*$-linearization. Every line bundle over $\Aff^2$ is trivial.  Since
$\scrO_{\Aff^2}^* =\Cee^*$, it is easy to see that every
$\Cee^*$-linearization on the trivial bundle is given by a character
$\chi\colon \Cee^*\to \Cee^*$. We denote the corresponding linearized line
bundle by $(\scrO_{\Aff^2}, \chi)$. In terms of the total space
$\Aff^2\times \Cee$, the action is via $\mu\cdot (g_2, g_3, z) = (\mu^4g_2,
\mu^6g_3, \chi(\mu)z)$. Suppose that
$F\colon B^* \to \Aff^2$ is a
$\Cee^*$-equivariant morphism corresponding to $Z\to B$. Pulling back the
linearized bundle $(\scrO_{\Aff^2}, \chi)$ gives a linearized bundle
$(\scrO_{B^*}, \chi)$. Clearly, the line bundle over $B$
corresponding to the $\Cee^*$-linearized bundle $(\scrO_{B^*}, \chi)$ is
exactly $B^*\times_{\Cee^*}\Cee$, where $\Cee^*$ acts on $\Cee$ via
$\chi^{-1}$. The characters of
$\Cee^*$ are exactly the $\chi_a$ defined by $\chi_a(\mu) =\mu^a$, and for
each $a\in \Zee$ we let $\mathbf{L}^a$ be the line bundle over
$\mathbf{E}$ corresponding to the linearization on
$\scrO_{\Aff^2}$ given by
$\chi_a$.  Since
$B^*$ is the
$\Cee^*$-bundle corresponding to $\mathcal{L}^{-1}$, where $\mathcal{L}$ is
the line bundle over $B$ defined in Theorem~\ref{Weiermodel}, it follows that
the line bundle over $B$ induced by $\mathbf{L}^a$ is $\mathcal{L}^a$. Thus:

\begin{theorem} The map $a\in \Zee \mapsto \mathbf{L}^a$ defines an
isomorphism from $\Zee$ to $\Pic \mathbf{E}$. Given
$a\in
\Zee$,  the corresponding bundle over $B$ is $\mathcal{L}^a$.
\qed
\end{theorem}

For example, using the remarks at the beginning of this section we recover
the well-known fact that
$\mathcal{L} = R^0\pi_*\omega_{Z/B}$.

Similarly, one can determine the algebra of the sections of $\bigoplus _{a\in
\Zee}\mathbf{L}^a$ over $\mathbf{E}$ (in an obvious sense):

\begin{proposition}\label{sectionsoverstack} 
$\mathbf{H}^0(\mathbf{E};\bigoplus _{a\in
\Zee}\mathbf{L}^a) = \Cee[G_2, G_3]$. Thus, if $\boldsymbol{\Phi} \colon
\bigoplus _i\mathbf{L}^{a_i}\to \bigoplus_j\mathbf{L}^{b_j}$ is a morphism of
vector bundles over $\mathbf{E}$ and if, for all $i,j$, either $b_i-a_j < 4$
or $b_i-a_j$ is odd,  then
$\boldsymbol{\Phi}=0$.\qed
\end{proposition}

We turn now to the classification of vector bundles over $\mathbf{E}$. 

\begin{theorem}\label{splits} Let $\mathbf{V}$ be a vector bundle over
$\mathbf{E}$. Then
$\mathbf{V} \cong \bigoplus _{i=1}^n\mathbf{L}^{a_i}$ for integers $a_i$,
uniquely determined up to order.
\end{theorem}
\begin{proof} We begin with the following lemma (which is a partial 
generalization of  Horrocks' splitting criterion for vector bundles over
$\Pee^n$ \cite[Theorem 2.3.2]{OSS}):

\begin{lemma}\label{graded} Let $R=\bigoplus _{k\geq 0}R_k$ be a graded ring
such that
$R_0=\Cee$, and let $\mathfrak{m} = \bigoplus _{k> 0}R_k$. Suppose that $M$
is a finitely generated graded $R$-module such that $M_{\mathfrak{m}}$ is a
free
$R_{\mathfrak{m}}$-module. Then $M$ is free as a graded $R$-module, i.e.\
there exists a homogeneous basis $e_1,
\dots, e_n$ for
$M$.
\end{lemma}
\begin{proof} The proof is by induction on the rank of $M_{\mathfrak{m}}$. If
$M_{\mathfrak{m}}=0$, then $M=\mathfrak{m}\cdot M$. We claim in this case that
$M=0$, i.e.\ that the graded Nakayama lemma holds. For, if $M\neq 0$, then 
since
$M$ is finitely generated there is a nonzero homogeneous element of smallest
degree in $M$, which cannot lie in $\mathfrak{m} \cdot M$. This
contradiction shows that $M=0$. 

If the rank of $M_{\mathfrak{m}}$ is positive, choose $f_1, \dots, f_n \in M$
mapping to a basis $\{\ov f_i\}$ of $M_{\mathfrak{m}}$. Let $e_i$ be the
initial part of
$f_i$ , i.e.\ the nonzero homogeneous component of $f_i$ of minimal degree.
Possibly after changing the basis
$\{\ov f_i\}$, we may clearly assume that, for some $k\geq 1$, the elements 
$e_1,
\dots, e_k$ have the smallest degree among the possible $e_i$ and that they
are linearly independent over
$\Cee$. Write $e_1 = a_1f_1 + \sum _{i\geq 2}a_if_i$. Then $a_1=
1+m$ with $m\in \frak m$. Thus the images of $e_1, f_2, \dots, f_n$ are also a
basis of $M_{\mathfrak{m}}$.  An easy argument shows that, if $re_1=0$, then
$r=0$ (since $e_1$ maps to a basis element of the free module
$M_{\mathfrak{m}}$  and the map $R\to R_{\mathfrak{m}}$ is injective by a
homogeneity argument).

 Let $Q = M/Re_1$. Clearly $Q$ is a finitely generated graded $R$-module, and
$Q_{\mathfrak{m}} = M_{\mathfrak{m}}/R_{\mathfrak{m}}e_1$. By the above,
$Q_{\mathfrak{m}}$ is a free $R_{\mathfrak{m}}$-module of rank $n-1$, and the
sequence
$$0 \to R \to M \to Q \to 0$$ is exact, where the first map is defined by
$r\mapsto re_1$.  By induction, there exists a homogeneous basis $\hat e_2,
\dots, \hat e_n$ for $Q$. Lift the
$\hat e_i$ to homogeneous elements $e_i\in M$. It follows that $e_1, \dots,
e_n$ is a homogeneous basis of $M$, as desired.
\end{proof}

Returning to the proof of the theorem, let $R = \Cee[g_2, g_3]$.  A
$\Cee^*$-linearized vector bundle
$V$ over $\Spec R$ corresponds to a finitely generated projective graded
$R$-module
$M$, and thus $M_{\mathfrak{m}}$ is free. The lemma then implies that there
is a homogeneous basis for $M$, and hence $V$ is a direct sum of
$\Cee^*$-linearized line bundles. The uniqueness statement is clear.
\end{proof}

The following lemma, whose proof is clear, shows how to determine the
integers $a_i$ in Theorem~\ref{splits}:

\begin{lemma}\label{linwts} Suppose that $\mathbf{V}$ corresponds to the 
$\Cee^*$-linearized vector bundle $V$ over
$\Aff^2$. Since $0\in \Aff^2$ is a fixed point for the $\Cee^*$-action,
there is an induced linear action of $\Cee^*$ on the fiber $V_0$ of $V$
over $0$. Let $a_1, \dots, a_n$ be the weights of this action. Then
$\mathbf{V} \cong \bigoplus _{i=1}^n\mathbf{L}^{a_i}$. 
\qed
\end{lemma}

\subsection{The stable bundle over the total space}

We begin by recalling the following general facts (see e.g \cite{FMW}). For
each Weierstrass cubic  $C$, with origin $p_0$, there is a unique stable
vector bundle $W_n$ of rank $n$ and such that $\det W_n =
\scrO_C(p_0)$. Given the fibration $Z\to B$, there is a unique vector bundle
$\mathcal{W}_n$ over $Z$ whose restriction to every fiber is isomorphic to
$W_n$ and such that there exists a  surjection
$p\colon \mathcal{W}_n\to \scrO_Z(\sigma)$. The bundle
$\mathcal{W}_n$ is a successive extension of $\scrO_Z(\sigma)$ by
$\mathcal{L}$,
$\mathcal{L}^2$, \dots, $\mathcal{L}^{n-1}$. Every other vector bundle
$\mathcal{W}$ over $Z$ whose restriction to every fiber is isomorphic to
$W_n$ is of the form
$\mathcal{W}_n\otimes \pi^*M$ for some line bundle $M$ over $B$, with $\det
(\mathcal{W}_n\otimes \pi^*M) = \scrO_Z(\sigma) \otimes
\pi^*\mathcal{L}^{n(n-1)/2}\otimes M^n$. In particular, 
$$\det \mathcal{W}_n\otimes \pi^*\mathcal{L}^a = \scrO_Z(\sigma) \otimes
\pi^*\mathcal{L}^{n(n-1)/2 + an}.$$ Since $W_n$ is 
simple, $\pi_*Hom(\mathcal{W}_n, \mathcal{W}_n) =\scrO_B$, and the
sheaf of relative automorphisms of $\mathcal{W}_n$ is
$\scrO_B^*$. To identify
$\mathcal{W}_n$ up to a unique isomorphism, it is best to consider pairs
$(\mathcal{W}, p)$, where $\mathcal{W}$ is a bundle whose restriction to
every slice is isomorphic to $W_n$ and $p\colon \mathcal{W}\to
\scrO_Z(\sigma)$ is a surjection.  It is easy to see that the pair
$(\mathcal{W}_n, p)$ is in fact a vector bundle  over the total space of
the stack
$\mathbf{E}$ (in a slightly generalized
sense, because of the choice of the surjection $p$). More generally,  the
vector bundles $\mathcal{W}_n\otimes
\pi^*\mathcal{L}^a$ for   $a\in
\Zee$ are also naturally bundles over the total space of the moduli stack.  

We now wish to do the same thing equivariantly for the family $\mathcal{E}\to
\Aff^2$. The section $\sigma$ is $\Cee^*$-equivariant, and thus there is a
$\Cee^*$-linearization on $\scrO_{\mathcal{E}}(-\sigma)$ compatible with the
inclusion $\scrO_{\mathcal{E}}(-\sigma) \subseteq \scrO_{\mathcal{E}}$ (where
$\scrO_{\mathcal{E}}$ is given the $\Cee^*$-linearization coming from
pullback of functions). From the isomorphism
$R^1\pi_*\scrO_{\mathcal{E}}(-\sigma) \cong R^1\pi_*\scrO_{\mathcal{E}}$, we
see that $\Cee^*$ acts on $R^1\pi_*\scrO_{\mathcal{E}}(-\sigma)$ via the
character $\chi_{-1}$. Thus, letting $(\scrO_{\Aff^2}, \chi_1)$ denote the
trivial bundle with the linearization given by $\chi_1$, and similarly for
the pulled back linearization $(\scrO_{\mathcal{E}}, \chi_1)$, we see that
there is an invariant extension class in
$\Ext^1(\scrO_{\mathcal{E}}(\sigma), (\scrO_{\mathcal{E}}, \chi_1))\cong
H^0(R^1\pi_*\scrO_{\mathcal{E}}(-\sigma) )$, and hence that there is a
$\Cee^*$-linearization on the extension
$\mathcal{W}_2$ of $\scrO_{\mathcal{E}}(\sigma)$ by $\scrO_{\mathcal{E}}$
which is compatible with the linearizations on the factors. Likewise, we can
obtain  a $\Cee^*$-linearization inductively on $\mathcal{W}_n$, such that
there is an exact sequence of linearized vector bundles
$$0 \to (\scrO_{\mathcal{E}}, \chi_{n-1}) \to \mathcal{W}_n \to
\mathcal{W}_{n-1} \to 0$$ (compare \cite[Prop.\ 4.4]{FMW}).  Thus, using the
equivalence of categories between $\mathbf{E}$ and $[\Aff^2/\Cee^*]$, there
is a (non-canonical) vector bundle $\mathbf{W}_n$ over the total space of 
$\mathbf{E}$, defined as follows: given a triple
$(B^*,p,F)\in [\Aff^2/\Cee^*]$ corresponding to $(Z, \pi,\sigma)$, define
$\mathcal{W}_n\to Z$ as the vector bundle
$F^*\mathcal{W}_n/\Cee^*$ over $F^*\mathcal{E}/\Cee^* = Z$. By construction,
there is a surjection $\mathbf{W}_n \to \mathbf{W}_{n-1}$ for every $n$.

\subsection{A universal unstable conformal bundle over the moduli
stack}\label{conformbund}

Let $G$ be a simple and simply connected complex linear group. In the
introduction, we have defined a \textsl{special root} $\alpha$ of $G$, cf.\
also  
 \cite[\S3.1]{FMII}.  Let
$P$ be a maximal parabolic subgroup of $G$ associated to a special root
and let
$L$ be the Levi factor. These groups have been studied in \cite[\S1.2,
\S3.2]{FMII}. In particular 
$$L \cong \left\{(A_1, \dots, A_t) \in \prod
  _{i=1}^tGL_{n_i}(\Cee): \det A_1= \cdots = \det A_t\right\},$$ where $t\leq
3$.  Thus, the bundles $W_{n_i}\spcheck$, $1\leq i\leq t$, define an
$L$-bundle, which was denoted by $\eta_0$ in \cite{FMII}. Of course, we can
define the corresponding $L$-bundle for a singular Weierstrass cubic as well,
and  we shall continue to denote it by
$\eta_0$.

Let
$\pi\colon Z
\to B$ be an elliptic fibration. We would like to construct, in a functorial
way, an
$L$-bundle $\eta \to Z$ which restricts on every smooth fiber to the bundle
$\eta_0$ constructed in \cite{FMII}. While this is possible for many groups,
it is not always possible, and instead we shall have to work with a conformal
extension of the form $\hat G = G\times _{\Zee/n\Zee}\Cee^*$, where
$\Zee/n\Zee$ is embedded in $G$ as a subgroup of the center.

An $L$-bundle over $Z$ consists of $t$ vector bundles $V_1, \dots, V_t$ over
$Z$, where the rank of $V_i$ is $n_i$, and such that $\det V_1 = \cdots =
\det V_t$. We must have $V_i = (\mathcal{W}_n\otimes
\pi^*\mathcal{L}^{a_i})\spcheck$. Thus the condition that the $V_i$ define an
$L$-bundle is that it is possible to choose integers $a_i$ so that the
integers
$n_i(n_i-1)/2 + a_in_i$ are all equal. This is trivially satisfied if $t=1$. 
If $t=2$ and at least one of $n_1, n_2$ is odd, then we can again find
integers $a_1, a_2$ such that 
$$n_1(n_1-1)/2 + a_1n_1=n_2(n_2-1)/2 + a_2n_2.$$ To see this, if both $n_1$
and $n_2$ are odd, then this equation becomes
$A_1n_1 = A_2n_2$, where $A_i =(n_i-1)/2 + a_i$. Choosing $a_i$ so that $A_1 =
n_2$ and $A_2=n_1$ gives a solution. Likewise, if say $n_1$ is odd and $n_2$
is even, then the equation becomes $A_1n_1 = A_2n_2/2$, where $A_1 =(n_1-1)/2
+ a_1$ and $A_2 = n_2-1+2a_2$. Choosing $a_i$ so that $A_1 = n_2/2$ and 
$A_2=n_1$, which is possible since $n_1$ is odd and $n_2$ even, gives a
solution. This constructs the desired $L$-bundle for  $A_n$, $n$ even and
half of the special vertices for $A_n$ and $n$ odd, $B_n$ and $n$ even,
$C_n$, $F_4$, $G_2$. In case $t=3$, if two of $n_1, n_2, n_3$ are equal and
at least one is odd, the same argument goes through. This handles $D_n$,
where $(n_1, n_2, n_3) = (2,2,n-2)$, for $n$ odd, as well as $E_6$, where
$(n_1, n_2, n_3) = (2,3,3)$. The case of $E_8$ can be handled directly. Here
$(n_1, n_2, n_3) = (2,3,5)$, so we seek $a_1, a_2, a_3$ such that 
$$1+2a_1 = 3+3a_2 = 10 + 5a_3.$$ This equation can be solved by taking
$a_1=7, a_2 = 4, a_3 = 1$.

Some of the remaining cases can also be solved directly. For $B_n$, a
solution exists exactly when $n$ is not congruent to $1$ mod $4$. For $D_n$, a
solution exists exactly when $n$ is not congruent to $2$ mod $4$. In the case
of $SL_n(\Cee)$, where we consider the maximal parabolic with Levi factor
$SL_{2k}(\Cee) \times SL_{2\ell}(\Cee)$ we can solve the problem if and only
if  $k$ and $\ell$ are both divisible by the same power of $2$. But in
general we need to take a conformal extension of
$G$. In the following, we omit the case of $SL_n$ since that has been
discussed extensively elsewhere \cite{FMW}.

\begin{theorem}\label{conf1} Let $\pi\colon Z \to B$ be a Weierstrass
fibration. Suppose that
$G$ is not of type $A_n$. For $G$ of type  $B_n$,
$n\equiv 1 \mod 4$, $D_n$, $n\equiv 2 \mod 4$, or $E_7$, let
$\hat G = G\times _{\Zee/2\Zee}\Cee^*$, where $\Zee/2\Zee$ is included in
$\Cee^*$, is the full center for $G$ of type $B_n$ or $E_7$, and is the
subgroup whose quotient is $SO(2n)$ in case $G=Spin(2n)$, and let $\hat L =
L\times _{\Zee/2\Zee}\Cee^* \subseteq \hat G$. In all other cases, let $\hat
G = G$ and $\hat L =L$. Then there exists an $\hat L$-bundle $\hat
\eta_0$ over $Z$ which
lifts to the $L$-bundle $\eta_0$ on every fiber.
\end{theorem}
\begin{proof} In the case of $B_n$, $n$ odd, the group $\hat L\cong
GL_2(\Cee)\times GL_{n-1}(\Cee)$. In the case of $D_n$, $n$ even, the group
$$\hat L\cong
\{ (B_1, B_2, B_3) \in GL_2(\Cee)\times GL_2(\Cee)\times GL_{n-2}(\Cee):
\det B_1 = \det B_2\}. $$ In each of these cases the result is clear. In the
case of $E_7$, 
$$\hat L\cong
\{ (B_1, B_2, B_3) \in GL_2(\Cee)\times GL_3(\Cee)\times GL_4(\Cee):
\det B_2\det B_3 = (\det B_1)^2\}. $$ Choosing bundles
$V_1=(\mathcal{W}_2\otimes \pi^*\mathcal{L}^{a_1})\spcheck$,
$V_2=(\mathcal{W}_3\otimes \pi^*\mathcal{L}^{a_2})\spcheck$,
$V_3=(\mathcal{W}_4\otimes
\pi^*\mathcal{L}^{a_3})\spcheck$, we see that we want to choose $a_1, a_2,
a_3$ so that
$$3+3a_2 + 6+4a_3 = 2+4a_1.$$ This can be solved with $a_1=0, a_2 = a_3 =-1$.
\end{proof}

We can again make an equivariant version of the above construction over
$\mathcal{E}\to\Aff^2$. In this version, we replace the bundle
$\mathcal{W}_n\otimes \pi^*\mathcal{L}^a$ by the linearized bundle
$\mathcal{W}_n\otimes (\scrO_{\mathcal{E}}, \chi_a)$, where
$\mathcal{W}_n$ is the linearized bundle over $\mathcal{E}$ constructed in
\S 1.3, with the natural linearization on the tensor product. The
constructions of Theorem~\ref{conf1} then give, in this context:

\begin{theorem}\label{conf2} Let $\hat G$ and $\hat L$ be as in
Theorem~\ref{conf1}. Then there exists a $\Cee^*$-linearized $\hat L$-bundle
$\hat
\eta_0$ over $\mathcal{E}$  which
lifts to the $L$-bundle $\eta_0$ on every fiber. \qed
\end{theorem}

Once we have the  linearized $\hat L$-bundle $\hat \eta_0$ over
$\mathcal{E}$, by using the equivalence of categories between
$\mathbf{E}$ and $[\Aff^2/\Cee^*]$, there is a (non-canonical) principal
$\hat L$-bundle $\boldsymbol{\hat\eta_0}$ over the total space of 
$\mathbf{E}$.

The unipotent radical $U$ of $L$ is still the unipotent radical of $\hat L$.
Thus, we can form the sheaf of unipotent groups $U(\hat \eta_0)$ over $Z$ or
$\mathcal{E}$, depending on the context. If $\hat \eta_0$ reduces to an
$L$-bundle $\eta_0$, then clearly $U(\hat \eta_0)=U( \eta_0)$. We shall
discuss the representability of the corresponding functor in \S\ref{functor}.

\section{Vector bundles over singular curves}

In this section, we collect standard facts about vector bundles on cuspidal
and nodal curves of arithmetic genus one.  

\subsection{Vector bundles over cuspidal curves}

We begin with the cuspidal case. Let $C$ be cuspidal, with normalization
$\iota \colon \widetilde C\to C$, and suppose that $0\in C$ is the singular
point. Let
$R_0 =\scrO_{C,0}$ be the (analytic) local ring at $0$ and let
$\widetilde R_0$ be the normalization of $R_0$. Fixing once and for all a
coordinate
$t$ at the preimage of $0$ in $\widetilde C$, we can identify $\widetilde
R_0$ with
$\Cee\{t\}$ and
$R_0$ with $\Cee\{t^2, t^3\}$. Let $\widetilde {\frak m}_0 =(t)$ be the
maximal ideal of
$0$ in
$\widetilde R_0$, and let
$\frak m_0 =(t^2,t^3)R_0$ be the maximal ideal of $R_0$.

\begin{lemma}\label{cusp1} Let $U$ be a finite-dimensional complex vector
space and let $\widetilde M$ be the free
$\widetilde R_0$-module $U\otimes _\Cee \widetilde R_0$, which we can view as
the $\widetilde R_0$-module of germs of morphisms $m\colon \widetilde C\to U$
at $0$. Then there is a one-to-one correspondence between
projective {\rm(}i.e.\ free{\rm)} $R_0$-submodules
$M$ of
$\widetilde M$ such that the natural map
$M\otimes _{R_0}\widetilde R_0 \to \widetilde M$ is an isomorphism, and
endomorphisms $A$ of $U$. If $M$ corresponds to the endomorphism $A$, then 
$$M = \left\{m\in \widetilde M: \frac{dm}{dt}(0) = Am(0)\right\}.$$ 
\end{lemma}
\begin{proof} Suppose that $M$ is a projective $R_0$-submodule   of
$\widetilde M$ such that the natural map
$\widetilde M\otimes _{R_0}\widetilde R_0 \to \widetilde M$ is an
isomorphism. Since $R_0/\frak m_0 \cong \widetilde R_0/\widetilde {\frak
m}_0$, there is a canonical isomorphism
$$M/\frak m_0 M \cong M\otimes _{R_0}(R_0/\frak m_0)\cong M\otimes
_{R_0}(\widetilde R_0/\widetilde{\frak m}_0) \cong \widetilde M/\widetilde
{\frak m}_0\widetilde M
\cong U.$$ Of course, the image of $m\in M$ under the composition $M \to
M/\frak m_0 M
\cong U$ is just $m(0)$, viewing $M$ as a submodule of $\widetilde M$.  Define
the endomorphism
$A$ as follows: given
$u\in U$, choose
$m\in M$ such that $m(0)=u$, and set
$$Au = \left.\frac{dm}{dt}\right|_{t=0}.$$
It is easy to see that $Au$ only depends on $u$, not the choice of $m$, and
defines a linear map from $U$ to itself. Clearly 
$$M = \left\{m\in \widetilde M: \frac{dm}{dt}(0) = Am(0)\right\}.$$
Thus, $A$ determines $M$ as well, and every possible $A$ arises in this way.
\end{proof} 

For example, the submodule $U\otimes _\Cee R_0$ corresponds to 
$0\in \End U$.

In the above situation, suppose that $U_1$ and $U_2$ are two
finite-dimensional vector spaces and that $T\colon U_1\otimes _\Cee
\widetilde R_0 \to U_2\otimes _\Cee
\widetilde R_0$ is a homomorphism of free $\widetilde R_0$-modules. Let $M_1$
be a projective $R_0$-submodule of $U_1\otimes _\Cee
\widetilde R_0$ and let $M_2$ be a projective submodule of $U_2\otimes _\Cee
\widetilde R_0$, such that the natural homomorphism $M_i\otimes _{R_0}
\widetilde R_0 \to U_i\otimes _\Cee
\widetilde R_0$ is an isomorphism for $i=1,2$, and such that $T(M_1)
\subseteq M_2$. Then clearly, for all $m\in M_1$,
$$\left.\frac{d(Tm)}{dt}\right|_{t=0}= \left.\frac{dT}{dt}\right|_{t=0}\cdot
m(0) + T(0)\cdot\left.\frac{dm}{dt}\right|_{t=0}.$$
In particular, if $U_1=U_2 = U$ and $M_1=M_2$, suppose that $T=\Id +t^2B$ for
some $B$. Then automatically $T$ has coefficients in $R_0$, so that
$T(M_1)=M_1$. By the above remarks, $d(Tm)/dt|_{t=0}=dm/dt|_{t=0}$. Hence
the endomorphism
$A$ only depends on the isomorphism $\widetilde M_0 \otimes _{R_0}\widetilde
R_0/(\widetilde{\frak m}_0)^2 \cong U\otimes _\Cee \Cee\{t\}/(t^2)$, which we
refer to as a \textsl{first order trivialization} of $\widetilde M_0$.
However, we shall not use this in what follows.

Passing from the local ring of a cuspidal curve to the global setting of the
curve itself, we immediately obtain:

\begin{corollary}\label{sections} Let $\widetilde V$ be a vector bundle on
$\widetilde C$, and fix a  trivialization at $0$, i.e.\ an isomorphism from 
the stalk of
$\widetilde V$ at
$0$ to $U\otimes _\Cee\Cee\{t\}$, where $U$ is a
finite-dimensional vector space. Then there is a one-to-one correspondence
between locally free
$\scrO_C$-subsheaves $V$ of
$\iota_*\widetilde V$, such that the natural map $\iota^*V \to \widetilde V$
is an isomorphism, and elements $A\in \End U$. Under this correspondence,
sections $s$ of $V$ correspond to sections $\widetilde s$ of
$\widetilde V$ such that 
$$\frac{d\widetilde s}{dt}(0) = A\widetilde s(0),$$ where the derivative is
with respect to the  trivialization of
$\widetilde V$ at $0$, and similarly for local sections.
\qed
\end{corollary}

\begin{corollary}\label{linalg} Let $\widetilde V_1$ and $\widetilde V_2$ be  
vector bundles on
$\widetilde C$,  and fix isomorphisms from the stalk of $\widetilde V_i$ at
$0$ to $U_i\otimes _\Cee\Cee\{t\}$. Let $V_1$ and
$V_2$ be locally free subsheaves of $\iota_*\widetilde V_1$ and
$\iota_*\widetilde V_2$ respectively, such that the natural map $\iota^*V_i
\to \widetilde V_i$ is an isomorphism, and let
$A_{V_i}\in
\End U_i$ be the corresponding elements. Then, using the induced
trivializations of the stalks of the corresponding bundles,
\begin{enumerate}
\item[\rm (i)] The subbundle $V_1\spcheck$ of $\iota_*\widetilde V_1\spcheck$
corresponds to the element $-A_{V_1}^*\in \End U_1^*$;
\item[\rm (ii)] The subbundle $V_1 \otimes V_2$ of $\iota_*(\widetilde
V_1\otimes \widetilde V_2)$ corresponds to the element $ A_{V_1}\otimes \Id +
\Id \otimes A_{V_2}\in
\End (U_1\otimes U_2)$;
\item[\rm (iii)]  Homomorphisms
$f\colon V_1 \to V_2$ correspond to homomorphisms
$\widetilde f\colon \widetilde V_1 \to \widetilde V_2$ such that 
$$\frac {d\widetilde f}{dt}(0) = A_{V_2}\circ \widetilde f - \widetilde f
\circ A_{V_1}.$$
\item[\rm (iv)] The subbundle $\bigwedge ^kV_1$ of
$\iota_*(\bigwedge^k\widetilde V_1)$ corresponds to the element
$A_k\in
\End(\bigwedge^k U_1)$ which is the coefficient of $t$ in $\wedge^k(\Id +
tA)$.
\end{enumerate}
\end{corollary}
\begin{proof} These are standard constructions via elementary linear algebra.
\end{proof}

We first apply the corollary to bundles which become trivial on the
normalization. In general, if $\widetilde V$ is a vector bundle over
$\widetilde C$ with a given trivialization $\widetilde V \cong U\otimes_\Cee
\scrO_{\widetilde C}$, where $U$ is a vector space, we always use the
fixed global trivialization to define a  trivialization of
$\widetilde V$ at
$0$. By the remarks after Lemma~\ref{cusp1}, changing the choice of
global trivialization has the effect of conjugating the matrix $A$ by an
invertible matrix $T$.

\begin{corollary}\label{conditions} 
\begin{enumerate}
\item[\rm (i)] Let $\widetilde V = U\otimes_\Cee \scrO_{\widetilde C}$ be  the
trivial rank $n$ vector bundle over $\widetilde C$, where $U$ is a
finite-dimensional vector space. Let 
$V$ be the locally free subsheaf of $\iota_*\widetilde V$ corresponding to
$A\in \End U$. Then the composed map
$$H^0(C;V) \to H^0(\widetilde C; \iota^*V) \cong H^0(\widetilde C; \widetilde
V) \cong U$$ is injective and its image is equal to $\Ker A \subseteq U$.
\item[\rm (ii)] Suppose that, for $i=1,2$,  $\widetilde V_i = U_i\otimes_\Cee
\scrO_{\widetilde C}$  are trivial vector bundles over $\widetilde C$, where
$U_i$ are finite-dimensional vector spaces, and that $V_i\subseteq
\iota_*\widetilde V_i$ correspond to $A_i\in \End U_i$. Then 
$$\Hom (V_1, V_2)  \cong\{f\in \Hom (U_1, U_2): A_2\circ f = f\circ A_1\}.$$
\item[\rm (iii)] There is a one-to-one correspondence between isomorphism
classes of rank $n$ vector bundles
$V$ on $C$ such that $\iota^*V$ is a trivial  vector bundle on
$\widetilde C$ and conjugacy classes of $n\times n$ matrices.
\end{enumerate}
\end{corollary}
\begin{proof} The corollary follows easily by using the fact that every
section of the trivial bundle over $\widetilde C$ is constant, and so has
derivative zero.
\end{proof}

Next let us describe the bundle $W_n\spcheck$. In \cite{FMW}, $W_n\spcheck$ is
defined inductively as follows: $W_1\spcheck=\scrO_C(-p_0)$, and for $n>1$,
there is a unique nonsplit extension
$$0 \to W_{n-1}\spcheck \to W_n\spcheck\to \scrO_C \to 0.$$  It follows that
$\iota^*W_n\spcheck \cong \scrO_{\widetilde C}(-p_0) \oplus
\scrO_{\widetilde C}^{n-1}$. Choose a local trivialization $e_1, \dots, e_n$
of 
$\scrO_{\widetilde C}(-p_0) \oplus
\scrO_{\widetilde C}^{n-1}$ at $0$ which is compatible with the direct sum
decomposition, is standard on the trivial factors, and is the pullback to
$\scrO_{\widetilde C}(-p_0) = \iota^* \scrO_C(-p_0)$ of a local trivialization
of   $\scrO_C(-p_0)$ at $0$. For example,   we can take for $e_1$
the standard local section $1$ which has a pole at $p_0$, viewed as a
section of  $\scrO_C(-p_0)$. We now use this choice to describe
$W_n\spcheck$. Of course, the dual description gives $W_n$.

\begin{lemma}\label{cuspcomps} With this choice of local trivialization, the
matrix
$A$ corresponding to $W_n\spcheck$ is the upper triangular matrix with ones
along the superdiagonal and zeroes elsewhere.
\end{lemma}
\begin{proof} The proof is by induction. For $n=1$ the result is clear.
Suppose that we have established the result for $W_n\spcheck$. Let $X_n$ be
the bundle defined by the $n\times n$ matrix described in the statement of the
lemma. Then there is an exact sequence
$$0\to X_n \to X_{n+1} \to \scrO_C \to 0,$$ and it will suffice to prove that
this extension is not split. A splitting would be a homomorphism $f\colon
\scrO_C\to X_{n+1}$, or equivalently a section of $X_{n+1}$, projecting to the
identity on the last factor. Using the standard trivialization of $\scrO_C$
with basis
$e$ and the given local trivialization $e_1, \dots, e_{n+1}$ of 
$\scrO_{\widetilde C}(-p_0) \oplus
\scrO_{\widetilde C}^{n}$ at $0$, we see that 
$$\widetilde f(e)= \sum _{i=1}^n f_i\cdot e_i + e_{n+1}.$$ Since $\Hom
(\scrO_{\widetilde C}, \scrO_{\widetilde C}(-p_0)) =0$, $f_1 = 0$. Similarly,
the remaining coefficients $f_i$ are constant. Thus
$\displaystyle\frac{d\widetilde f}{dt}(0) = 0$. By (i) of
Corollary~\ref{conditions} applied to the section $f$, we must have
$A_1\circ
\widetilde f(0) = 0$, where
$A_1$ defines
$X_{n+1}$. But 
$$A_1\circ \widetilde f(0) = \sum _{i=2}^{n}f_ie_{i-1} + e_n \neq 0,$$
a contradiction. Hence the sequence is not split.
\end{proof}

\subsection{The nodal case}

Here we let $C$ be a nodal Weierstrass cubic with normalization
$\iota\colon \widetilde C\to C$. Suppose that the singular point of $C$ is $0$
and that $\iota^*(0) =\{x,y\}$. Let
$R_0 =\scrO_{C,0}$ be the (analytic) local ring at $0$ and let $\widetilde
R_0$ be the normalization. Fixing once and for all  coordinates $t_1, t_2$ at
$x$ and $y$ respectively in $\widetilde C$, we can identify $R_0$ with
$\Cee\{t_1, t_2\}/(t_1t_2)$ and
$\widetilde  R_0$ with $\Cee\{t_1\}\oplus \Cee\{t_2\}$. Let $\frak m_0
=(t_1, t_2)R_0$ be the maximal ideal of $R_0$, let $\frak m_x= t_1\widetilde
R_0$ be the maximal ideal in
$\widetilde R_0$ corresponding to the point $x$, and similarly for $\frak
m_y$.

\begin{lemma} Let $\widetilde M$ be a free $\widetilde R_0$-module. Then
there is a one-to-one correspondence between projective $R_0$-submodules $M$
of $\widetilde M$ such that the natural map
$M\otimes _{R_0}\widetilde R_0 \to \widetilde M$ is an isomorphism, and
isomorphisms $A\colon \widetilde M/ \frak m_x\widetilde M \to \widetilde M/
\frak m_y\widetilde M$. Moreover, under this correspondence, the $R_0$-module
$M$ corresponding to $A$ is given by
$$M = \left\{m\in \widetilde M: Am(x) = m(y)\right\}.$$
\end{lemma}
\begin{proof} In this case, $R_0/\frak m_0 \cong  \widetilde R_0/\frak m_x
\cong \widetilde R_0/\frak m_y$. As before, given a projective
$R_0$-submodule 
$M$   of $\widetilde M$ such that the natural map
$\widetilde M\otimes _{R_0}\widetilde R_0 \to \widetilde M$ is an isomorphism,
the induced homomorphism $M/\frak m_0M \to \widetilde M/\frak m_x\widetilde
M$ is an isomorphism. Given $u\in \widetilde M/\frak m_x\widetilde
M$, choose $m\in M$ which maps to $u$ under the composition $M \to M/\frak
m_0M \to \widetilde M/\frak m_x\widetilde M$, and let $Au$ be the image of
$m$ in $\widetilde M/ \frak m_y\widetilde M$ under the analogous composition 
$M \to M/\frak m_0M \to \widetilde M/\frak m_y\widetilde M$. Clearly, this
image only depends on $u$, not on the choice of $m$, and defines a linear
isomorphism $A$ from  $\widetilde M/ \frak m_x\widetilde M$ to $\widetilde M/
\frak m_y\widetilde M$, and all such isomorphisms arise in this way.
\end{proof}

In case we have an isomorphism $\widetilde M \cong U\otimes _\Cee \widetilde
R_0$, where $U$ is a finite-dimensional vector space, the fibers $\widetilde
M/ \frak m_x\widetilde M$ and $\widetilde M/
\frak m_y\widetilde M$ are both canonically identified with $U$, and thus $A$
corresponds to a linear automorphism of $U$.

\begin{proposition}\label{nodegluing} Let $\widetilde V$ be a vector bundle on
$\widetilde C$. Then there is a one-to-one correspondence between locally free
$\scrO_C$-subsheaves $V$ of
$\iota_*\widetilde V$, such that the natural map $\iota^*V \to \widetilde V$
is an isomorphism, and isomorphisms $A\colon \widetilde V(x)\to \widetilde
V(y)$. Under this correspondence, sections $s$ of $V$ correspond to
sections
$\widetilde s$ of
$\widetilde V$ such that 
$$\widetilde s(y) = A\widetilde s(x),$$
and similarly for local sections.\qed
\end{proposition}

\begin{corollary} Let $\widetilde V_1$ and $\widetilde V_2$ be   vector
bundles on
$\widetilde C$. Let $V_1$ and
$V_2$ be locally free subsheaves of $\iota_*\widetilde V_1$ and
$\iota_*\widetilde V_2$ respectively, satisfying the hypotheses of the above
corollary, and let $A_{V_i}\colon V_i(x) \to V_i(y)$ be the corresponding
elements. Then
\begin{enumerate}
\item[\rm (i)] The subbundle $V_1\spcheck$ of $\iota_*\widetilde V_1\spcheck$
corresponds to the element $(A_{V_1} ^{-1})^*\colon \widetilde V_1(x)^* \to
\widetilde V_1(y)^*$;
\item[\rm (ii)] The subbundle $V_1 \otimes V_2$ of $\iota_*(\widetilde
V_1\otimes \widetilde V_2)$ corresponds to the element $ A_{V_1}\otimes  
A_{V_2}$;
\item[\rm (iii)]  Homomorphisms
$f\colon V_1 \to V_2$ correspond to homomorphisms
$\widetilde f\colon \widetilde V_1 \to \widetilde V_2$ such that 
$$  A_{V_2}\circ \widetilde f(x) = \widetilde f(y) \circ A_{V_1}. \qed$$
\end{enumerate}
\end{corollary}

We shall usually be interested in the case where there is a fixed
trivialization $\widetilde V = U\otimes \scrO_{\widetilde C}$, where $U$ is a
finite-dimensional vector space. Using this trivialization to identify the
fibers  $\widetilde V(x)$ and $\widetilde V(y)$ with $U$ means that we can
identify $A$ with a fixed automorphism of $U$. Changing the trivialization
has the effect of conjugating $A$ by an element of $\Aut U$.

\begin{corollary} 
\begin{enumerate}
\item[\rm (i)] Let $\widetilde V = U\otimes \scrO_{\widetilde C}$ be  the
trivial rank $n$ vector bundle over $\widetilde C$, where $U$ is a
finite-dimensional vector space. Let 
$V$ be the locally free subsheaf of $\iota_*\widetilde V$ corresponding to
$A\in \Aut U$. Then the composed map
$$H^0(C;V) \to H^0(\widetilde C; \iota^*V) \cong H^0(\widetilde C; \widetilde
V) \cong U$$ is injective and its image is equal to $\Ker (A-\Id) \subseteq
U$.
\item[\rm (ii)] Suppose that, for $i=1,2$,  $\widetilde V_i = U_i\otimes
\scrO_{\widetilde C}$  are trivial vector bundles over $\widetilde C$, where
$U_i$ are vector spaces, and that $V_i\subseteq \iota_*\widetilde V_i$
correspond to $A_i\in \Aut U_i$. Then 
$$\Hom (V_1, V_2)  \cong\{f\in \Hom (U_1, U_2): A_2\circ f = f\circ A_1\}.$$
\item[\rm (iii)] There is a one-to-one correspondence between isomorphism
classes of  rank $n$ vector bundles
$V$ on $C$ such that $\iota^*V$ is a trivial vector bundle on
$\widetilde C$ and conjugacy classes of invertible $n\times n$ matrices.\qed 
\end{enumerate}
\end{corollary}

We have the description of the bundle $W_n\spcheck$ corresponding  to Lemma
~\ref{cuspcomps}: given the isomorphism
$\iota^*W_n\spcheck \cong \scrO_{\widetilde C}(-p_0) \oplus
\scrO_{\widetilde C}^{n-1}$, choose a local trivialization $e_1, \dots,
e_n$ of 
$\scrO_{\widetilde C}(-p_0) \oplus
\scrO_{\widetilde C}^{n-1}$ at $x$ and $y$ as in the cuspidal case. Then:

\begin{lemma}\label{nodecomps} With this choice of local trivialization, the
matrix
$A$ corresponding to $W_n\spcheck$ is the matrix satisfying $Ae_1 = e_1$
and $Ae_i = e_i+ e_{i-1}$ for $i>1$.\qed
\end{lemma}

\subsection{Semistable vector bundles supported at the singular point}

Let $V$ be a torsion free sheaf on $C$ which is semistable of degree zero (but not
necessarily with trivial determinant). Via the Fourier-Mukai correspondence
\cite{Mukai, Teo, Spectral}, there is an equivalence of categories between the
category of all such
$V$ and the category of torsion $\scrO_C$-modules $M$. We define the
\textsl{support} of $V$ to be the support of $M$. 

Assume that $C$ is singular, with singular point $0$. Let $\mathcal{F}$ be the
unique rank one torsion free, non-locally free sheaf on $C$ of degree zero, so
that the support of $\mathcal{F}$ is $\{0\}$ and $\mathcal{F}$ corresponds to
the
$\scrO_C$-module
$\scrO_C/\frak m_0$. We shall only consider torsion sheaves supported at $0$, 
and let $R_0= \scrO_{C,0}$ be the analytic local ring of $C$ at $0$. Hence
$R_0=\Cee\{x,y\}/(xy)$ if $C$ is nodal and $R_0= \Cee\{x,y\}/(x^2-y^3) =
\Cee\{t^2,t^3\}$ if $C$ is cuspidal. 
Finally, recall that a semistable torsion free sheaf $V$ of degree zero on $C$
and with support equal to $\{0\}$ is \textsl{strongly indecomposable} if and
only if $\dim \Hom (V, \mathcal{F}) =1$ if and only if $\dim \Hom (\mathcal{F},V)
=1$. Every strongly indecomposable vector bundle $V$ with $\det V\cong
\scrO_C$ arises from the parabolic construction \cite{FMW} and has a
nontrivial pullback to the normalization of
$C$. Thus $\iota^*V \cong \scrO_{\Pee^1}(1) \oplus \scrO_{\Pee^1}^{n-2} \oplus
\scrO_{\Pee^1}(-1)$. More generally, the vector bundles $V$ with $\det V
\cong \scrO_C$ arising from the parabolic construction for $SL_n(\Cee)$ are
exactly the vector bundles of the form $V_0\oplus V_1$, where $V_0$ is a
strongly indecomposable vector bundle supported at the singular point, and
$V_2$ has no support at the singular point and is regular.

\begin{lemma} Let $V$ be a semistable torsion free sheaf  of degree zero on $C$
and with support equal to $\{0\}$, and let $M$ be the corresponding $R_0$-module.
Then:
\begin{enumerate}
\item[\rm (i)] The sheaf $V$ is locally free if and only if $M$ has homological
dimension one as an $R_0$-module.
\item[\rm (ii)] The sheaf $V$ is strongly indecomposable if and only if
$M\cong R_0/I$ for some ideal $I$ in $R_0$ such that $R_0/I$ has finite length.
Hence $V$ is
locally free and strongly indecomposable if and only if
$M\cong R_0/(f)$, where $f$ is not a zero divisor or a unit in $R_0$.
\end{enumerate}
\end{lemma}
\begin{proof} The proof of (i) is in \cite{Teo}. To see (ii), the equivalence of 
categories implies that $V$ is strongly indecomposable if and only if $\dim
\Hom_{R_0}(M, R_0/\frak m_0) = 1$. First suppose that $M=R_0/I$. Then clearly
$\dim
\Hom_{R_0}(R_0/I,R_0/\frak m_0) = 1$. Conversely, suppose that $\dim
\Hom_{R_0}(M, R_0/\frak m_0) = 1$, and choose
$\varphi \neq 0$ in $\Hom_{R_0}(M, R_0/\frak m_0)$, and  $m\in M$ such that
$\varphi(m) \neq 0$. The map $r\in R_0 \mapsto rm$ defines an exact sequence
$$0 \to R_0/I \to M \to \ov M \to 0,$$
where $I$ is the annihilator of $m$ and $\ov M =M/R_0m$. Hence there is an
exact sequence
$$0 \to \Hom_{R_0}(\ov M, R_0/\frak m_0) \to \Hom_{R_0}(M, R_0/\frak m_0) \to
\Hom_{R_0}(R_0/I, R_0/\frak m_0).$$
By construction, the image of $\varphi$ in $\Hom_{R_0}(R_0/I,
R_0/\frak m_0)$ is nonzero. By assumption,  $\dim
\Hom_{R_0}(M, R_0/\frak m_0) = 1$, and thus $\Hom_{R_0}(\ov M, R_0/\frak m_0) =0$,
which implies that
$\ov M =0$ and
$M=R_0/I$. Finally, if $M$ has homological dimension one, then $I$ is a rank one
torsion free $R_0$-module of homological dimension zero, and it then follows
directly from the classification of all such (cf.\ \cite[(0.2)]{FMW}) that $I$ is
principal.
\end{proof}

Direct calculation then gives the following   \cite{Teo}:

\begin{corollary} \begin{enumerate}
\item[\rm (i)] If $C$ is nodal, then every semistable  vector bundle $V$ of 
rank $n$ on $C$, supported at the singular point and strongly
indecomposable, is of the form
$V_{a,b;\lambda}$ for
$a,b$ are positive integers such that $a+b=n$ and $\lambda \in \Cee^*$,
corresponding to $M=R_0/(x^a+\lambda y^b)$. For an appropriate choice of the
coordinates $x,y$ and an identification of $\Pic^0C$ with $\Cee^*$, we may assume
that $\det V_{a,b;\lambda} =\lambda$, and shall abbreviate $V_{a,b;1}=V_{a,b}$.
\item[\rm (ii)] If $C$ is cuspidal, then every semistable vector bundle $V$ 
of  rank $n$ on $C$, supported at the singular point and strongly
indecomposable, is of the form
$V_{n;\lambda}$ for some $n\geq 2$ and  $\lambda \in \Cee$,
corresponding to $M=\Cee[[t^2,t^3]]/(t^n+ \lambda t^{n+1})$. Here $\lambda =0$
corresponds to trivial determinant, and we shall abbreviate
$V_{n;\lambda}=V_n$.\qed
\end{enumerate}
\end{corollary}

It is easy to check that the dual of $V_{a,b;\lambda}$ is $V_{b,a;
\lambda^{-1}}$, and the dual of $V_{n;\lambda} $ is $V_{n;-\lambda}$. 

The next corollary also follows easily by direct calculation.

\begin{corollary} With notation as in the preceding corollary, suppose that $C$ is
nodal and that $a,b\geq 2$. For each pair $c,d$ of positive integers such that
$a-c\geq 1$ and $b-d\geq 1$ and for $\lambda', \lambda'' \in \Cee^*$ with
$\lambda =\lambda'\lambda''$, there is an exact sequence
$$0 \to V_{c,d;\lambda'} \to V_{a,b;\lambda} \to V_{a-c, b-d;\lambda''} \to
0.$$ Likewise, suppose that $C$ is cuspidal, that $n\geq 4$, and that $2\leq
n' \leq n-2$. Then, for $\lambda', \lambda''\in \Cee$ and $\lambda
=\lambda'+\lambda''$, there is an exact sequence
$$0\to V_{n';\lambda'} \to V_{n;\lambda} \to V_{n-n';\lambda''} \to 0.\qed$$
\end{corollary}

The following is the main application of the above description:

\begin{corollary}\label{wedgeunstable} In case $C$ is nodal, suppose that $a,b\geq
2$ and that
$n=a+b$. For all
$k$ such that $2\leq k\leq n-2$, the vector bundle
$\bigwedge^kV_{a,b;\lambda}$ is unstable. Likewise, if $C$ is cuspidal and
$n\geq 4$, then for all
$k$ such that $2\leq k\leq n-2$, the vector bundle $\bigwedge^kV_{n;\lambda}$
is unstable.
\end{corollary}
\begin{proof} We shall just write down the proof in the nodal case and for
$\lambda=1$; the proof in the cuspidal case is essentially identical.
Beginning with the exact sequence
$$0 \to V_{1,1} \to V_{a,b} \to V_{a-1, b-1} \to 0,$$
we have the Koszul filtration on $\bigwedge^kV_{a,b}$  whose associated gradeds
are $\bigwedge^kV_{a-1,b-1}$, $V_{1,1} \otimes \bigwedge^{k-1}V_{a-1,b-1}$, and
$\bigwedge^{k-2}V_{a-1,b-1}$. Since each graded piece has degree zero,
$\bigwedge^kV_{a,b}$ is unstable provided that one of the graded pieces is
unstable. Furthermore, if  $\bigwedge^{k-1}V_{a-1,b-1}$ is unstable, then so is
$V_{1,1} \otimes \bigwedge^{k-1}V_{a-1,b-1}$. Thus we may assume that
$\bigwedge^{k-1}V_{a-1,b-1}$ is semistable. Note that the rank of $V_{a-1,b-1}$
is $n-2$ and that $1\leq k-1 \leq n-3$. Thus, since $\iota^*V \cong
\scrO_{\Pee^1}(1) \oplus \scrO_{\Pee^1}^{n-2} \oplus
\scrO_{\Pee^1}(-1)$, $\iota^*\bigwedge^{k-1}V_{a-1,b-1}$ is not the trivial
vector bundle on $\widetilde C$. Since by assumption $\bigwedge^{k-1}V_{a-1,b-1}$
is semistable, it must have some support at the singular point $\{0\}$, and in
particular there is a nonzero homomorphism $\mathcal{F} \to
\bigwedge^{k-1}V_{a-1,b-1}$. Identifying $V_{1,1}$ with its dual gives an
identification $V_{1,1} \otimes
\bigwedge^{k-1}V_{a-1,b-1}\cong Hom(V_{1,1},
\bigwedge^{k-1}V_{a-1,b-1})$. The surjection $V_{1,1}\to \mathcal{F}$ then
induces an inclusion $Hom (\mathcal{F},\mathcal{F}) \to V_{1,1} \otimes
\bigwedge^{k-1}V_{a-1,b-1}$. But $Hom (\mathcal{F},\mathcal{F})\cong
\iota_*\scrO_{\widetilde C}$ has degree one, and so $V_{1,1} \otimes
\bigwedge^{k-1}V_{a-1,b-1}$ is unstable. It follows that $\bigwedge^kV_{a,b}$ is
unstable as well.
\end{proof}

We shall see below that, in contrast, $\bigwedge^kV_{1,n-1}$ and
$\bigwedge^kV_{n-1,1}$ are always semistable.

\subsection{Vector bundles via \'etale covers}

In Proposition~\ref{nodegluing}, we described bundles on a nodal curve $C$ via
bundles on the normalization together with gluing maps.  In this section, we
describe another method for constructing   vector bundles
over a nodal curve, using instead the nontrivial
\'etale covers of degree
$n$. (See Teodorescu \cite{Teo} for the case $n=2$.) In fact, this is a
special case of Proposition~\ref{nodegluing}, where the bundle on the
normalization is a direct sum of line bundles and the gluing maps are of a
special type, but it will be convenient to have this second description as
well. 

\begin{theorem} Let $C$ be a nodal curve, and let $\pi\colon C_n \to C$ be the
unique connected \'etale cover of $C$ of degree $n$. Let $q\in C_n$ be a smooth
point. Then $\pi_*\scrO_{C_n}(q)$ is a stable vector bundle on $C$ of rank $n$
and degree one, and moreover every stable vector bundle on $C$ of rank $n$
and degree one arises in this way. Likewise, $\pi_*\scrO_{C_n}(-q)$ is stable of
degree $-1$, and every stable  rank $n$ bundle on $C$ of degree $-1$ arises in
this way.
\end{theorem}
\begin{proof} Since $\pi$ is \'etale, it is a finite
flat morphism of degree $n$, and thus
$W=\pi_*\scrO_{C_n}(q)$ is a vector bundle on
$C$ of rank
$n$ and degree one. Moreover, for every line bundle $\lambda$ of degree
zero on $C$, $W\otimes
\lambda=
\pi_*(\scrO_{C_n}(q)
\otimes
\pi^*\lambda) = \pi_*\scrO_{C_n}(q')$ for some smooth point $q'$. Since every
stable vector bundle on $C$ of rank $n$ and degree one is of the form $W_n\otimes
\lambda$ for some line bundle $\lambda$ of degree zero, it will suffice to prove
that $\pi_*\scrO_{C_n}(q)$ is stable.

Applying $\pi_*$ to the inclusion $\scrO_{C_n} \to \scrO_{C_n}(q)$ gives an exact
sequence
$$0 \to \bigoplus_{i=0}^{n-1}\theta^i \to W \to \Cee_p\to 0,$$
where $p=\pi(q)$ and $\theta$ is the $n$-torsion line bundle in $\Pic C$ defining
$C_n$. Suppose that $W$ is unstable. Then there is a quotient $Q$ of $W$ of
degree at most zero. The image of $\bigoplus_{i=0}^{n-1}\theta^i$ in $Q$ is 
a subsheaf of $Q$ with torsion cokernel $\tau$ and hence its degree is $\deg
Q-\ell(\tau) \leq 0$, with equality if and only if $\tau=0$ and $\deg Q=0$.
Since $\bigoplus_{i=0}^{n-1}\theta^i$ is semistable of degree zero, it follows
then that $\deg Q=0$ and that the induced map $\bigoplus_{i=0}^{n-1}\theta^i
\to Q$ is surjective. Since $\bigoplus_{i=0}^{n-1}\theta^i$ is a direct sum of
line bundles of degree zero, this is only possible if $Q$ is a direct
summand of $\bigoplus_{i=0}^{n-1}\theta^i$. Hence the surjection $W\to Q$ also
splits, so that $W =Q\oplus W'$, where $\deg W' =1$. In particular, by the
Riemann-Roch theorem, $h^0(C;W'\otimes \lambda)\geq 1$ for every line bundle
$\lambda$ of degree zero. Clearly, by choosing $\lambda = \theta^j$ for an
appropriate power of $\theta$, we can arrange that $h^0(C;Q\otimes \lambda) \neq
0$. Thus $h^0(C; W\otimes \lambda) \geq 2$. On the other hand,
$$h^0(C; W\otimes \lambda) = h^0(C; \pi_*(\scrO_{C_n}(q) \otimes
\pi^*\lambda)) = h^0(C_n; \scrO_{C_n}(q) \otimes
\pi^*\lambda) =1,$$
as one checks by a straightforward direct calculation on $C_n$. This is a
contradiction, so that $W$ is stable. The  statements about $W_n\spcheck$
follow from a very similar argument, or from relative duality applied to the
\'etale morphism $\pi$.
\end{proof} 

\begin{corollary}\label{tensor}  
Let $f\colon C_{n_1}\times _CC_{n_2} \times _C\cdots \times _CC_{n_k} \to C$,
$\pi_i\colon C_{n_i}\to C$, and $p_i \colon C_{n_1}\times _CC_{n_2} \times
_C\cdots
\times _CC_{n_k}
\to C_{n_i}$  be  the natural projections.  Then there exists smooth points $q_i
\in C_{n_i}$ such that the bundle
$$W_{n_1} \otimes \cdots\otimes W_{n_a} \otimes W_{n_{a+1}}\spcheck \otimes
\cdots
\otimes W_{n_k}\spcheck $$
is isomorphic to 
$$f_* \left(p_1^*\scrO_{C_{n_1}}(q_1)\otimes \cdots \otimes
p_a^*\scrO_{C_{n_a}}(q_a)\otimes
p_{a+1}^*\scrO_{C_{n_{a+1}}}(-q_{a+1})\otimes \cdots \otimes
p_k^*\scrO_{C_{n_k}}(-q_k)\right).$$
\end{corollary}
\begin{proof} An easy exercise shows that, in the case of two factors 
$C_{n_1}$ and $C_{n_2}$, if $\mathcal{S}_i$ is a coherent sheaf on $C_{n_i}$,
then 
$f_*\left(p_1^*\mathcal{S}_1\otimes p_2^*\mathcal{S}_2\right) =
(\pi_1)_*\mathcal{S}_1\otimes (\pi_2)_*\mathcal{S}_2$. The general case
follows by induction.
\end{proof}

We now use the  \'etale cover $\pi\colon C_n \to C$ to construct semistable
vector bundles on $C$ of degree zero. Let us begin with some notation. Order
the components of $C_n$ as $D_0, D_1, \dots, D_{n-1}$, where $D_i\cap D_{i+1}
\neq \emptyset$. Let $d_0, \dots, d_{n-1}$ be integers. We shall view the
subscripts of the $D_i$ or the $d_i$ as integers modulo $n$. Denote by
$\scrO_{C_n}(d_0,
\dots, d_{n-1})$ a line bundle of multidegree
$(d_0,
\dots, d_{n-1})$ on $C_n$, i.e.\ a line bundle whose restriction to each
component $D_i$ has degree $d_i$. This does not specify
$\scrO_{C_n}(d_0, \dots, d_{n-1})$ uniquely, but it is specified up to 
tensoring by a line bundle of the form $\pi^*\lambda$, where $\lambda$ is a
line bundle on
$C$ of degree zero. Thus, the rank $n$ vector bundle $\pi_*\scrO_{C_n}(d_0,
\dots, d_{n-1})$ is well-defined up to tensoring with a line bundle  $\lambda$   
of degree zero. 

Assume now that  $\sum _id_i = 0$. Then $\deg \pi_*\scrO_{C_n}(d_0,
\dots, d_{n-1}) =0$, and we wish to decide when it is semistable. If all of the
$d_i$ are zero, then, up to tensoring by a line bundle of degree zero, $\pi_*\scrO_{C_n}(d_0,
\dots, d_{n-1}) =\bigoplus_{i=0}^{n-1}\theta^i$, where $\theta$ is a line
bundle of degree zero, and thus $\pi_*\scrO_{C_n}(d_0,
\dots, d_{n-1})$ is semistable. If $d_i\geq 2$ for some $i$, then it is easy to
see that $\pi_*\scrO_{C_n}(d_0,
\dots, d_{n-1})$ is unstable, and by relative duality or a direct argument the
same is true if $d_i\leq -2$ for some $i$, although we shall not use this fact.
Thus we may assume that all of the $d_i$ are $0$ or $\pm 1$. The next lemma tells
us when $\pi_*\scrO_{C_n}(d_0,
\dots, d_{n-1})$ is semistable:

\begin{lemma}\label{ssetale} With notation as above, suppose that, for every
$i$,
$d_i$ is either  $0$ or $\pm 1$ and that $\sum _id_i=0$. Then 
$\pi_*\scrO_{C_n}(d_0,
\dots, d_{n-1})$ is semistable if and only if, for every $i$, if $d_i=1$,
$d_{i+1} = \cdots = d_{i+k-1} = 0$ and $d_{i+k}$ is nonzero, then $d_{i+k} = -1$.
In other words, the nonzero multidegrees must alternate in sign. In this case, 
$\pi_*\scrO_{C_n}(d_0,
\dots, d_{n-1})$ is supported at the singular point $0$. Finally, 
$\pi_*\scrO_{C_n}(d_0,\dots, d_{n-1})$ is strongly indecomposable if and only if
exactly one of the $d_i$ is equal to $1$ and one is equal to $-1$.
\end{lemma}
\begin{proof} A vector bundle $V$ of degree zero on $C$ is unstable if
and only if, for every line bundle $\lambda$ of degree zero on $C$, $h^0(C;
V\otimes \lambda) \neq 0$, and $V$ is semistable with support equal to $\{0\}$ if
and only if, for every line bundle $\lambda$ of degree zero on $C$, $h^0(C;
V\otimes \lambda) = 0$. Now
$$H^0(C; \pi_*\scrO_{C_n}(d_0, \dots, d_{n-1})\otimes \lambda) = H^0(C_n;
\scrO_{C_n}(d_0, \dots, d_{n-1})\otimes \pi^*\lambda).$$
Moreover $\scrO_{C_n}(d_0, \dots, d_{n-1})\otimes \pi^*\lambda$ has the same
multidegree as $\scrO_{C_n}(d_0, \dots, d_{n-1})$, namely $(d_0, \dots,
d_{n-1})$. Let $\lambda$ be a line bundle  of degree zero on $C$. Direct
computation then shows that $H^0(C_n;
\scrO_{C_n}(d_0, \dots, d_{n-1})\otimes \pi^*\lambda) =0$ if the nonzero $d_i$
alternate in sign and that $H^0(C_n;
\scrO_{C_n}(d_0, \dots, d_{n-1})\otimes \pi^*\lambda)\neq 0$ if there are two
$d_i$ equal to $1$ which are separated only by zeroes. This proves the first two
statements.

To see the remaining statement, it suffices to show, for the unique
non-locally free torsion free rank one sheaf $\mathcal{F}$ of degree zero,
that
$$\dim
\Hom (\mathcal{F},
\pi_*\scrO_{C_n}(d_0,\dots, d_{n-1})) =1$$ if exactly one of the $d_i$ is equal to
$1$ and one is equal to $-1$, and that the dimension is greater than $1$
otherwise. To do so, note that $\mathcal{F} = \iota_*\scrO_{\widetilde C}\otimes
\scrO_C(-p_0)$. Moreover,
$$\Hom (\mathcal{F},\pi_*\scrO_{C_n}(d_0,\dots, d_{n-1}))=\Hom
(\pi^*\mathcal{F},\scrO_{C_n}(d_0,\dots, d_{n-1})).$$
But $\pi^*\mathcal{F} \cong  (\iota_n)_*\scrO_{\widetilde C_n}\otimes
\pi^*\scrO_C(-p_0)$, where $\iota_n \colon \widetilde C_n \to C_n$ is the
normalization map. There is a canonical isomorphism
$Hom((\iota_n)_*\scrO_{\widetilde C_n}, \scrO_{C_n})\cong I$, where $I$ is
the ideal of the reduced singular locus of $C_n$, i.e.\ $I=\pi^*\frak m_0$. Thus 
\begin{align*}
\Hom(\pi^*\mathcal{F},\scrO_{C_n}(d_0,\dots, d_{n-1})) &= H^0(C_n;
 Hom ((\iota_n)_*\scrO_{\widetilde C_n}\otimes \pi^*\scrO_C(-p_0),
\scrO_{C_n}(d_0,\dots, d_{n-1})) \\
&=H^0(C_n; I\otimes \scrO_{C_n}(d_0+1,\dots,
d_{n-1}+1)).
\end{align*}
In other words, on each component $D_i\cong \Pee^1$, we look for sections of
$\scrO_{\Pee^1}(d_i+1)$ vanishing at the two points of $D_i$ lying in the singular
locus. Clearly, the only contributions are from the components $D_i$ where $d_i+1
= 2$, i.e.\ $d_i =1$, and each such component contributes an independent section.
Thus, $\dim \Hom (\mathcal{F},
\pi_*\scrO_{C_n}(d_0,\dots, d_{n-1})) =1$ if exactly one of the $d_i$ is $1$ and
it is greater than one in all other cases.
\end{proof}

\begin{remark} One can show that, if $(d_0, \dots, d_{n-1})$ satisfies: $d_0
= 1, d_{a} = -1$, and $d_i =0$ for $i\neq 0,a$, then
$\pi_*\scrO_{C_n}(d_0,\dots, d_{n-1})$ is isomorphic to $V_{a, n-a;\lambda}$,
where  $\lambda$ depends on the actual gluing maps. We shall prove the special
case of this where $a=1$ below.
\end{remark}

\begin{proposition}\label{wedgess} In the above notation, suppose that $d_0
=1$,
$d_1=-1$ and all of the remaining $d_i$ are $0$. Then  $\bigwedge
^k\pi_*\scrO_{C_n}(d_0,\dots, d_{n-1})$ is semistable and supported at the
singular point.
\end{proposition}
\begin{proof} We may assume that $k< n$. The vector bundle $\bigwedge
^k\pi_*\scrO_{C_n}(d_0,\dots, d_{n-1})$ is obtained as follows. For each
$k$-element subset  $\{a_1, \dots, a_k\}$ of $\Zee/n\Zee$, we have a component
$D_{\{a_1, \dots, a_k\}}\cong \Pee^1$ and a line bundle on $D_{\{a_1, \dots,
a_k\}}$ of degree
$d_{a_1} + 
\cdots + d_{a_k}$. The component $D_{\{a_1, \dots, a_k\}}$ is glued to
$D_{\{a_1+1, \dots, a_k+1\}}$, where the addition is in $\Zee/n\Zee$, and the
line bundles are glued in an appropriate way (which we shall not need to
describe). This gives an \'etale cover $C'$ of $C$, with several components, and a
line bundle on $C'$ whose push-forward is $\bigwedge ^k\pi_*\scrO_{C_n}(d_0,\dots,
d_{n-1})$.

Clearly, $d_{a_1} +  \cdots + d_{a_k}$ is either $0$, $1$, or $-1$, and we must
show that the nonzero degrees alternate in a suitable sense. Suppose that
$D_{\{a_1, \dots, a_k\}}$ is a component where the line bundle has degree one.
Then we may assume that $a_1 =0$ and that we have ordered the remaining $a_i$ so
that $2\leq a_2 < a_3< \cdots < a_k$. Then $D_{\{a_1, \dots, a_k\}}$ is glued to
$D_{\{a_1+1, \dots, a_k+1\}}$, where $a_1+1 = 1$ and the $a_i+1$ are increasing
unless $a_k = n-1$, in which case $a_k+1 = 0$. If we are not in this last case,
the degree on $D_{\{a_1+1, \dots, a_k+1\}}$ is $-1$, and otherwise it is zero.
Assuming that the degree is zero, then the next component is $D_{\{a_1+2, \dots,
a_k+2\}}$, $a_k+2 =1$, $a_1+2=2$, and the integers $a_i+2$, $2\leq i\leq k-1$
are in increasing order unless $a_{k-1} =n-2$, in which case $a_{k-1}+2 =0$.
Thus either the degree is $-1$ or it is zero and $\{a_1+2, \dots,
a_k+2\} =\{0,1, a_2+2, \dots , a_{k-2}+2\}$. Continuing in this way, either we
reach degree $-1$ at some stage or all the degrees are zero until at stage $k$ we
reach $\{0, 1, \dots, k-1\}$. At the next stage, we get $\{1, \dots, k\}$. where
by assumption $k<n$, and the corresponding degree is then $-1$.
\end{proof}

\begin{corollary}\label{wedgesemistable} If $d_0 =1$, $d_1=-1$ and
all of the remaining $d_i$ are $0$, then the bundle $\pi_*\scrO_{C_n}(d_0,\dots,
d_{n-1})$ is either $V_{1,n-1}\otimes \lambda$ or $V_{n-1,1}\otimes \lambda$,
where $\lambda $ is a line bundle of degree zero. Hence $\bigwedge ^kV_{1,n-1}$
and
$\bigwedge ^kV_{n-1,1}$ are semistable for all $k$.
\end{corollary}
\begin{proof} Under our assumptions, $\pi_*\scrO_{C_n}(d_0,\dots,
d_{n-1})$ is semistable, supported at the singular point, and strongly
indecomposable. By Corollary~\ref{wedgeunstable}, $\pi_*\scrO_{C_n}(d_0,\dots,
d_{n-1})$ cannot be isomorphic to $V_{a,b}\otimes \lambda$ if both $a$ and $b$
are at least
$2$. Thus $\pi_*\scrO_{C_n}(d_0,\dots,
d_{n-1})$ is either $V_{1,n-1}\otimes \lambda$ or $V_{n-1,1}\otimes \lambda$.
Hence at least one of
$\bigwedge ^kV_{1,n-1}$ or
$\bigwedge ^kV_{n-1,1}$ is semistable for all $k$, and by duality so is the
other.
\end{proof}

\section{$G$-bundles over singular curves}

\subsection{$G$-bundles over cuspidal curves}

As in the case of vector bundles, we begin by assuming that $C$ is cuspidal.
Let
$G$ be a reductive linear algebraic group of rank $r$. We want to generalize
our results on vector bundles to
$G$-bundles.  Our first goal is to prove:

\begin{theorem}\label{cuspbundles} Let $\xi$ be a principal
$G$-bundle over $C$, let $C_0$ be a Zariski open subset of $C$ containing the
singular point,  and let
$\tau\colon
\iota^*(\xi|C_0)
\to \widetilde C_0 \times G$ be an   isomorphism  of
principal
$G$-bundles. Then:
\begin{enumerate}
\item[\rm (i)] There exists a unique right invariant vector field $X$ on $G$
such that  the local sections
of $\xi|C_0$ are  the local sections of $\iota^*\xi|\widetilde C_0$ which,
when viewed via the isomorphism $\tau$ as local   functions
$f\colon \widetilde C_0 \to G$,  satisfy 
$$f_*\left(\frac{\partial }{\partial t}\right) = -X_{f(0)}\in T_{f(0)}G.$$
\item[\rm (ii)]  For $C_0=C$, this correspondence sets up a bijection between
pairs
$(\xi,\tau)$ as above and right invariant vector fields on $G$, which we can
identify with $\frak g$.
\item[\rm (iii)] In the situation of {\rm (ii)}, if the pair  $(\xi,\tau)$ 
corresponds to
$X\in
\frak g$, then $\Aut \xi \cong \{g\in G: \ad g(X) = X\}$. In particular,
$\dim \Aut
\xi
\geq r$, with equality if and only if $X$ is a regular element of $\frak g$.
\end{enumerate}
\end{theorem}
\begin{proof} We shall just write out the case $C_0=\Spec R$, where $R=
\Cee[t^2, t^3]$ and   $\widetilde C_0=\Spec
\widetilde R$, where
$\widetilde R =\Cee[t]$; the general case follows by localizing the argument.
We start with a more general situation: Let
$M=\Spec S$ be a smooth affine scheme, and consider smooth morphisms $\xi_0
\to C_0$, locally trivial in the \'etale topology, together with an
isomorphism
$\iota^*\xi_0 \cong \widetilde C_0
\times M$.  Let
$\frak m = (t^2, t^3)R$ be the maximal ideal at the origin of $R$. We are
looking for finite $C_0$-morphisms
$\widetilde C_0
\times M \to \xi_0$ which identify $\widetilde C_0
\times M$ with $\iota^*\xi_0$. By a theorem of Chevalley \cite[ex.\ 4.2, p.\
222]{Hart}, $\xi_0$ is affine, and hence $\xi_0   = \Spec A$ for a subring
$A$ of $\widetilde R\otimes _\Cee S=S[t]$, which is an algebra over $R$. The
assumptions that $\xi_0 \to C_0$ is locally
trivial with fibers isomorphic to $M$ and that $\iota^*\xi_0\cong  \widetilde
C_0 \times M$ imply $A$ is flat over $R$ and that the scheme-theoretic fiber
of $\xi_0$ over $0\in C_0$ is $M$. These requirements easily imply that $A$
has the following two  properties:
\begin{enumerate}
\item $A/\frak m A\cong S[t]/tS[t] = S$ under the natural homomorphism;
\item The $R$-algebra $(\widetilde R\otimes _\Cee S)/A$ is annihilated by
$\frak m$.
\end{enumerate} The second property implies that $t^2S\subseteq A$. Using
this and the first, it follows that there exists a function $X\colon S\to S$
such that every element of $A$ has the form $s_0 + X(s_0)t +\sum _{i\geq
2}s_it^i$, and conversely every such element of $S[t]$ lies in $A$. The fact
that $A$ is a subring implies that $X$ is a derivation. A section of
$\widetilde C_0
\times M$ is of the form $x\mapsto (x, f(x))$, where $f$ is a morphism from
$\widetilde C_0$ to $\Spec S$ and thus corresponds to the homomorphism from
$S$ to $\Cee[t]$ given by $s\mapsto s\circ f$. In this case, the
corresponding homomorphism
$S[t]\to
\Cee[t]$ is the natural extension of the above homomorphism which sends $t$ to
$t$. A local calculation shows the following: suppose that the section
$x\mapsto (x, f(x))$ of
$\widetilde C_0
\times M$ has the property that it induces a section $C_0 \to \Spec A$, in
other words that the image of $A\subseteq S[t]$ under the above homomorphism
is contained in $\Cee[t^2,t^3]$. Then the function
$f\colon
\widetilde C_0
\to M$ satisfies
$$f_*\left(\left.\frac{\partial }{\partial t}\right|_{t=0}\right) =
-X_{f(0)}\in T_{f(0)}M,$$ and conversely every such function defines a
section. Clearly the $R$-automorphisms of $A$ are identified with the
$\widetilde R$-automorphisms of $\widetilde R\otimes _\Cee S$ fixing $X$.

Conversely, given a derivation $X$ of $S$, we can define the subring $A$ via
$$A=\{s_0 + X(s_0)t +\sum _{i\geq 2}s_it^i: s_i\in S\}.$$ Clearly $A$ is an
algebra over $R$ and $A\otimes _R\widetilde R \to S[t]$ is an isomorphism.
Let $\xi_0 = \Spec A$, so that the morphism $\widetilde C_0\times M \to
\xi_0$ is a bijection. We next show that the induced morphism
$\xi_0 \to C_0$ is smooth, in fact there are  \'etale local  cross sections.
Given $p\in M$ and $z_1, \dots, z_n$  coordinates at $p$, let
$z_i'=z_i + X(z_i)t$. Then $z_i'\in A$, and it is easy to check that the
maximal ideal of $A$ corresponding to $(0,p)$ is generated by $t^2, t^3,
z_1', \dots, z_n'$. The functions $t^2, t^3, z_1', \dots, z_n'$ define a  
morphism
$\xi_0\to C_0\times \Aff^n\subseteq \Aff^{n+2}$, sending $(0,p)$ to $(0,0)$,
say. The induced homomorphism on the completed local rings
$\widehat{\scrO}_{C_0\times\Aff^n,(0,0)} \to \widehat{\scrO}_{\xi_0, (0,p)}$
is then surjective. Since both rings have the same dimension and 
$\widehat{\scrO}_{ C_0\times
\Aff^n,(0,0)}$ is an integral domain, it follows that the homomorphism is an
isomorphism, and hence the map $\xi_0 \to C_0\times \Aff^n$ is \'etale at
$(p,0)$. In particular, some open subset of the preimage of $C_0\times \{0\}$
defines an \'etale cross section.

In case $M=G$, it is easy to see that, given a derivation $X$ of the affine
coordinate ring of $G$, the condition that the right multiplication action of
$G$ on $\widetilde C \times G$ induces a right action of $G$ on $\xi$ is
exactly that $X$ is right invariant. Thus, identifying the right invariant
vector fields on $G$ with $\frak g$ gives Part (i) of the  theorem. The proofs
of Parts (ii) and (iii) are straightforward and left to the reader.
\end{proof}

It is easy to check that, in case $G = GL_n(\Cee)$, the element $-X_{f(0)}\in
\frak g$ is just the matrix $A$ of Lemma~\ref{cusp1}.

Later, we shall need the following:

\begin{corollary}\label{muinverse}  Suppose that $\xi$ is a $G$-bundle over
$C$ and that
$\tau\colon \iota^*\xi \to \widetilde C \times G$ is a global trivialization
of $\iota^*\xi$. Let $X$ be the corresponding right invariant vector
field of $G$, identified with an element of $\frak g$. Under the
$\Cee^*$-action of Definition~\ref{action}, the bundle $\mu^*\xi$
corresponds to the vector field $\mu^{-1}\cdot X$.
\end{corollary}
\begin{proof} This follows from the construction  of
Theorem~\ref{cuspbundles} and the fact that $\mu^*t=\mu\cdot t$.
\end{proof}

In case the $G$-bundle $\iota^*\xi$ is trivialized on $\widetilde C$, we use
the local trivialization induced by the global one. Changing the global
trivialization by an element $g\in G$ has the effect of replacing $X\in \frak
g$ by $\ad (g)(X)$. Thus:

\begin{corollary} The correspondence of Theorem~\ref{cuspbundles}
sets up a bijection between isomorphism classes of
$G$-bundles $\xi$ over $C$ such that
$\iota^*\xi$ is isomorphic to the trivial bundle and  $\ad
G$-orbits in $\frak g$.
\qed
\end{corollary}

Let $\xi$ be a bundle such that  $\iota^*\xi$ is isomorphic to the
trivial bundle. We can define S-equivalence for the restricted class
of such bundles, in the usual way.  It is easy to see that this definition
is equivalent to the definition of S-equivalence, in the sense of geometric
invariant theory, for the adjoint action of $G$ on $\frak g$.  We shall call
such a bundle
$\xi$  \textsl{regular} if 
$\dim \Aut \xi = r$. Thus, every S-equivalence class of bundles which become 
isomorphic to the trivial bundle on
$\widetilde C$ contains a unique
regular representative. Furthermore, the adjoint quotient of $\frak g$, 
namely $\frak h/W$, where
$\frak h$ is a Cartan subalgebra of $\frak g$ and $W$ is the corresponding
Weyl group, describes the moduli space of S-equivalence classes of bundles
which become  isomorphic to the trivial bundle on
$\widetilde C$.

Another corollary of Theorem~\ref{cuspbundles} is:

\begin{corollary} Let $\xi$ be a holomorphic, topologically trivial $G$-bundle
over $C$. Then
$\dim H^1(C; \ad \xi) \geq r$.
\end{corollary}
\begin{proof} There is a small deformation of $\iota^*\xi$ to a semistable
bundle over $\widetilde C$, which is necessarily the trivial bundle. It is
easy to see that we can deform the choice of $X\in \frak g$ to obtain a small
deformation of $\xi$ to a $G$-bundle $\xi'$ whose pullback to $\widetilde C$
is trivial, and thus by  Theorem~\ref{cuspbundles} satisfies $h^0(C;
\ad \xi') = h^1(C; \ad \xi') \geq r$.
By semicontinuity, the same must hold for $\xi$.
\end{proof}

We also need a variant of Theorem~\ref{cuspbundles} for families.

\begin{theorem}\label{cuspfamilies} Let $B$ be a scheme. The isomorphism
classes of principal
$G$-bundles $\Xi\to C\times B$, together with $G$-isomorphisms $\iota^*\Xi
\cong \widetilde C\times B \times G$, are classified by morphisms $f\colon B
\to
\frak g$. Thus isomorphism classes of principal
$G$-bundles $\Xi\to C\times B$ such that $\iota^*\Xi
\cong \widetilde C\times B \times G$  are classified by morphisms $f\colon B
\to \frak g$ modulo the adjoint action of the group of morphisms from $B$ to
$G$.
\end{theorem}
\begin{proof} The arguments of the previous theorem show that it suffices to
classify global sections of the sheaf of derivations of $B\times G$, which
are linear over $\scrO_B$ and   invariant under the right action of $G$.
These may be identified with an element of $\scrO_B\otimes \frak g$, i.e.\ a
morphism from $B$ to $\frak g$.
\end{proof}

\subsection{The nodal case}

We turn now to $G$-bundles over nodal curves. As before, we assume that
$\{x,y\}\subseteq  \widetilde C$ is the preimage of the singular point. There
is the following analogue of Theorem~\ref{cuspbundles}, whose proof however is
simpler:

\begin{theorem}\label{nodebundles} Let $\xi$ be a principal
$G$-bundle over $C$, let $C_0$ be a Zariski open subset of $C$ containing the
singular point,  and let $\tau\colon \iota^*(\xi|C_0)
\to \widetilde C_0 \times G$ be an   isomorphism  of
principal $G$-bundles. Then:
\begin{enumerate}
\item[\rm (i)] There exists a unique $g\in G$
such that  the local sections
of $\xi|C_0$ are  the local sections of $\iota^*\xi|\widetilde C_0$ which,
when viewed via the isomorphism $\tau$ as local   functions
$f\colon \widetilde C \to G$,  satisfy 
$$f(x) =gf(y).$$
\item[\rm (ii)] For $C_0=C$, this correspondence sets up a bijection between
pairs
$(\xi,\tau)$ as above and $g\in G$.
\item[\rm (iii)] In the situation of {\rm (ii)}, if the pair  $(\xi,\tau)$ 
corresponds to
$g\in G$, then $\Aut \xi$ is isomorphic to the centralizer of $g$ in $G$. In
particular,
$\dim \Aut \xi \geq r$, with equality if and only if $g$ is a regular element
of $G$.
\end{enumerate}
\end{theorem}
\begin{proof} For simplicity, we assume
that
$\widetilde C_0 =\Spec \widetilde R$, where $\widetilde R = \Cee[t]$, and that
$C_0 = \Spec R$, where $R= \{(f_1, f_2)\in \widetilde R: f_1(x) = f_2(y)\}$.
We again consider the more general situation, where
$M=\Spec S$ is an affine scheme, not necessarily smooth, and $\xi_0\to C_0$
is a morphism such that $\iota^*\xi_0 \cong \widetilde C_0 \times M$.  Let
$\frak m$ be the maximal ideal of
$R$. As in the proof of  Theorem~\ref{cuspbundles}, we are looking for a
subring
$A$ of
$\widetilde R\otimes _\Cee S = S[t]$ which is an
$R$-subalgebra such that $A/\frak m A\cong S$ and such that
$\widetilde R\otimes _\Cee S/A$ is annihilated by $\frak m$. This implies
that there exists an automorphism $\varphi\colon M\to M$ such that
$$A=\Big\{s=\sum_{i\geq 0}s_it^i\in S[t]: s(x)=\varphi^*(s(y))\Big\}.$$
Clearly, as a ringed space, in the obvious categorical sense,
$$\xi_0 =\Spec A = \widetilde C_0\times M/(x, p)\sim (y,\varphi(p)),$$ and the
local sections of $\xi_0$ correspond to those morphisms $f\colon
\widetilde C_0 \to M$ such that $f(x) = \varphi(f(y))$.  Moreover, the
automorphisms of $\xi_0$ covering the identity on $C_0$ are exactly those
automorphisms $\psi$ of $\widetilde C_0 \times M$ covering the identity on
$\widetilde C_0$ and such that
$$\psi(x,p) = \psi (y, \varphi(p)).$$

Conversely, given
$A$, we shall show that there is an \'etale base change $C'_0\to C_0$ such
that
$\xi_0\times_{C_0}C'_0 \cong C'_0\times M$ as spaces over $C'_0$. In fact, let
$C_1=C_0-\{x\}$ and $C_2=C_0-\{y\}$, and take the Wirtinger double cover
$$C'_0 = C_1\amalg C_2/x\in C_2\sim y\in C_1.$$ The natural map $C'_0\to C_0$
is then an \'etale cover of degree $2$, and the preimage of the singular point
in $C_0$ is the unique point $(y,x)$ of $C'_0$. The fiber product
$\xi_0\times_{C_0}C'_0$ is given by
$$(C_1\times M)\amalg (C_2\times M)/(x, p)\in C_2\times M\sim
(y,\varphi(p))\in C_1\times M.$$ We can then define $\sigma \colon C'_0\times
M
\to \xi_0\times_{C_0}C'_0$ by
$\sigma|C_2\times M =(\Id,\Id)$ and  $\sigma|C_1\times M =(\Id,\varphi)$. It
is straightforward to check that $\sigma$ defines an isomorphism of ringed
spaces.

In case $M=G$, an automorphism $\varphi\colon G\to G$ commutes with right
multiplication if and only if $\varphi$ is given by left multiplication by
$g\in G$. This proves the first part of the theorem, and the proofs of
the second and third parts are straightforward.
\end{proof}

\begin{corollary} Isomorphism classes of $G$-bundles $\xi$ over $C$ such that
$\iota^*\xi$ is trivial are classified by $\Ad G$-orbits in $G$.\qed
\end{corollary}

Let $\xi$ be a bundle such that  $\iota^*\xi$ is trivial. We can again
define S-equivalence for such bundles, and it is the same as S-equivalence
in the sense of geometric invariant theory for the conjugation action of $G$
on itself. We shall call such a bundle
$\xi$ 
\textsl{regular} if 
$\dim \Aut \xi = r$. Thus, every S-equivalence
class of bundles which become trivial on $\widetilde C$ contains a unique
regular representative. Moreover,   the adjoint quotient of $G$, namely
$H/W$, describes the moduli space of S-equivalence
classes of bundles which become trivial on $\widetilde C$. 

As in the cuspidal case, we have the following:

\begin{corollary} Let $\xi$ be a holomorphic, topologically trivial
$G$-bundle over $C$. Then $\dim H^1(C; \ad \xi) \geq r$. \qed
\end{corollary}

There is also the following for families:

\begin{theorem}\label{nodefamilies} Let $B$ be a scheme. The   isomorphism
classes of principal
$G$-bundles $\Xi\to C\times B$, together with $G$-isomorphisms $\iota^*\Xi
\cong  \widetilde C\times B \times G$, are classified by morphisms $f\colon B
\to G$. Thus isomorphism classes of principal
$G$-bundles $\Xi\to C\times B$ such that $\iota^*\Xi
\cong \widetilde C\times B\times G$  are classified by morphisms $f\colon B
\to G$ modulo the adjoint action of the group of morphisms from $B$ to
$G$.\qed
\end{theorem}

\subsection{Bundles which become trivial on the normalization}

Here we shall let $C$ be either nodal or cuspidal. In this section,  we show
that the  bundles $\xi$ on $C$ such that $\iota^*\xi$ is isomorphic
to the trivial bundle are characterized by a strong semistability property: 

\begin{theorem}\label{unstable} Let $C$ be nodal or cuspidal 
and let $\xi$ be a $G$-bundle.  Then
$\iota^*\xi$ is trivial if and only if, for every representation $\rho\colon
G \to GL_N(\Cee)$, the associated vector bundle $\xi\times _G\Cee^N$ is
semistable.
\end{theorem}
\begin{proof}  It is easy to see that, if $V$ is an unstable vector bundle on
$C$, then $\iota^*V$ is unstable. Thus if $\iota^*V$ is semistable, then so
is $V$. Applying this remark to the vector bundle $\xi\times _G\Cee^N$
associated to a representation $\rho$, we see that, if
$\iota^*\xi$ is trivial, then so is $\iota^*\xi\times _G\Cee^N$, and hence
$\xi\times _G\Cee^N$ is semistable.

To prove the converse, we begin with the following notation. If
$\xi$ is a
$G$-bundle on $C$, then $\iota^*\xi$ is a $G$-bundle on $\widetilde C \cong
\Pee^1$.  By Grothendieck's theorem, every $G$-bundle on $\Pee^1$ reduces to
an
$H$-bundle, unique up to the action of the Weyl group. Denote by
$\Lambda$ the fundamental group of $H$, i.e.\  the coroot lattice   if $G$ is
simply connected. Using the exponential sheaf sequence
$$0\to \Lambda \to \frak h\otimes _\Cee\scrO_{\Pee^1} \to \underline{H}\to
1,$$ it follows that $H$-bundles over $\Pee^1$ are classified by
$\lambda \in\Lambda$. 

\begin{defn}\label{Pee} Given $\lambda \in \Lambda$, denote by
$\gamma_\lambda$ the corresponding $H$-bundle. 
\end{defn}

Next, there is the following lemma concerning the instability of vector
bundles on $C$.

\begin{lemma}\label{whenunstable} Let $V$ be a vector bundle of degree zero on
$C$ such that
$\iota^*V\cong
\bigoplus _i\scrO_{\Pee^1}(a_i)$. If $a_i \geq 2$ for some $i$, then $V$ is
unstable.
\end{lemma}
\begin{proof} Suppose for example that $C$ is nodal and let $x,y$ be the
preimages of the singular points. If
$a_i\geq 2$, then there is a section of $\scrO_{\Pee^1}(a_i)$ vanishing at
$x$ and $y$. It is easy to check that this defines an inclusion of
$\iota_*\scrO_{\Pee^1}$ in $V$. But $\deg \iota_*\scrO_{\Pee^1} = 1$, so that
$V$ is unstable. A similar argument works in the cuspidal case, where we
choose instead a section vanishing to order 2 at the preimage of the cusp
point.
\end{proof}

Returning to the proof of Theorem~\ref{unstable}, suppose that $\iota^*\xi$
is nontrivial. By the remarks at the beginning of the section, $\iota^*\xi =
\gamma_\lambda\times _HG$ for some $\lambda \in \Lambda$. If $\rho\colon G
\to \Cee^N$ is a representation with weights $\mu$, then 
$$\iota^*\xi\times _G\Cee^N = \bigoplus _\mu\scrO_{\Pee^1}(\mu(\lambda)).$$
Clearly, if $\lambda\neq 0$, then  by choosing $\rho$ correctly, we can make
some of the integers
$\mu(\lambda)$ arbitrarily large, and in particular at least 2. Thus, by the
preceding lemma, $\iota^*\xi\times _G\Cee^N$ is unstable.
\end{proof}

\subsection{A partial moduli space}\label{partialmodulispace}

If $C$ is a smooth elliptic curve, we can define the moduli space of
semistable holomorphic $G$-bundles $\mathcal{M}_C(G)$. When $C$ is singular,
there is no longer a good notion of semistability, defined for every group
$G$, which leads to a compact moduli space. However, we do have the following,
which can be proved for example by examining the arguments of Balaji
and Seshadri \cite{BS}, especially those in Section 8 and the proof of
Proposition~2.8:

\begin{theorem} Let $C$ be nodal or cuspidal. Then there is a coarse moduli
space $\mathcal{M}_C^0(G)$ whose points correspond to S-equivalence
classes of $G$-bundles
$\xi$ which pull back to the trivial bundle on the normalization, or
equivalently such that $\xi \times _G\Cee^N$ is a semistable vector bundle
for every representation $\rho\colon G \to GL_N(\Cee)$. Likewise, given a
Weierstrass fibration $\pi\colon Z\to B$, there is a relative moduli space
$\mathcal{M}_{Z/B}^0(G)\to B$. \qed
\end{theorem}

Here the method of proof of \cite{BS} shows that the equivalence used in
defining the coarse moduli space is the usual S-equivalence of bundles, where
all bundles satisfy the hypotheses of the theorem. Similarly, if $C$ is any
reduced irreducible curve, there is a coarse moduli
space $\mathcal{M}_C^0(G)$ whose points correspond to $G$-bundles $\xi$ which
pull back to semistable bundles on the normalization $\widetilde C$, modulo
S-equivalence.

We can describe the structure of the relative moduli space
$\mathcal{M}_{Z/B}^0(G)$ quite concretely. Let $Z_{\textrm{reg}}\to B$ be the
commutative group scheme defined by the regular points of the morphism $\pi$,
and let $\underline{\Lambda}$ be the product group scheme $B\times \Lambda$.
Then we can define a tensor product group scheme, which we denote by
$Z_{\textrm{reg}}\otimes \underline{\Lambda}$, and the Weyl group $W$ acts on
the tensor product by its action on $\Lambda$. Standard constructions then
prove the following:

\begin{theorem}\label{Eotomes} There is an $H$-bundle over $Z\times
_B(Z_{\rm reg}\otimes \underline{\Lambda})$ whose associated $G$-bundle
becomes trivial on the normalization of every singular fiber. The resulting
morphism $Z_{\rm reg}\otimes \underline{\Lambda} \to \mathcal{M}_{Z/B}^0(G)$
is $W$-invariant and identifies $\mathcal{M}_{Z/B}^0(G)$ with $(Z_{\rm
reg}\otimes \underline{\Lambda})/W$ as schemes over $B$. \qed
\end{theorem}

\section{The parabolic construction for singular curves and families}

\subsection{Cohomology computations for singular curves}\label{functor}

In this section $C$ will be either nodal or cuspidal, and $\iota\colon
\widetilde C \to C$ will denote the normalization. For a simply connected
group $G$,
$\alpha$ denotes the special root, $P$ is the corresponding maximal parabolic
subgroup, $U$ is the unipotent radical of $P$, and $L$ is a Levi factor of
$P$. Thus 
$$L \cong \left\{(A_1, \dots, A_t) \in \prod
  _{i=1}^tGL_{n_i}(\Cee): \det A_1= \cdots = \det A_t\right\},$$ where $t\leq
3$.  An
$L$-bundle   over $C$ thus corresponds to $t$ vector bundles $V_1, \dots,
V_t$, where $V_i$ has rank $n_i$, such that $\det V_1 \cdots =\det V_t$. We
define
$\eta_0$ to be the unique $L$-bundle corresponding to the bundles
$W_{n_1}\spcheck,
\dots, W_{n_t}\spcheck$. Its determinant, i.e.\ the line bundle associated to
the primitive dominant character of $L$, is $\scrO_C(-p_0)$.

We now  examine the parabolic construction for $\eta_0$, in other words
the set of liftings of $\eta_0$ to a $P$-bundle. The set of isomorphism
classes of pairs
$(\xi_P,\varphi)$ where is a $P$-bundle over $C$ and $\varphi$ is an
isomorphism from $\xi_P/U$ to $\eta_0$ is as usual the cohomology set $H^1(C;
U(\eta_0))$. There is also a corresponding functor from schemes over $\Cee$
to sets defined in \cite[4.1.1]{FMII}: if $S$ is a scheme over $\Cee$, let
$\mathbf{F}(S)$ be the set of isomorphism classes of pairs $(\Xi_P, \Phi)$,
where $\Xi_P$ is a 
$P$-bundle over $C\times S$ and $\Phi\colon \Xi_P/U \to \pi_1^*\eta_0$ is an
isomorphism of $L$-bundles. In this paper, in contrast to the notational
conventions of \cite{FMII}, since our main interest is in
$G$-bundles, we shall use the letter $\xi$ to denote the $G$-bundle
$\xi_P\times_PG$, and similarly for $\Xi$. We wish to study
$H^1(C; U(\eta_0))$ and its linearized analogue  $H^1(C;
\frak u(\eta_0))=\bigoplus _{k> 0}H^1(C; \frak u^k(\eta_0))$. The center of
$L$ contains a subgroup isomorphic to $\Cee^*$ which acts on $H^1(C;
U(\eta_0))$ and $H^1(C;
\frak u(\eta_0))$. 
Let us begin by recalling the situation for a smooth elliptic curve
$E$
\cite[Corollary 2.1.7 and Theorem 3.3.1]{FMII}. In this case, the
semistability of
$\eta_0$ and the fact that $\frak u^k$ is an irreducible $L$-module imply
that $\frak u^k(\eta_0)$ is a semistable vector bundle over $E$ of negative
degree. Thus, $H^0(E;\frak u^k(\eta_0))=0$ for all $k\geq 1$, and so by the
Riemann-Roch theorem the dimension of $H^1(E;\frak u (\eta_0))$ is simply
given by the negative of the degree of $\frak u (\eta_0)$, which can be
computed to be $r+1$.

In case
$C$ is singular, it is rarely the case that the vector bundles $\frak
u^k(\eta_0)$ are semistable. However, we can directly enumerate them
in Table~\ref{Table1}. As a guide to working out the bundle $\frak
u^k(\eta_0)$, note that $\frak u^k$ is an irreducible $L$-module with highest
weight
$\lambda_k(\alpha)$, the highest root such that $\alpha$ occurs with
coefficient $k$. Moreover, the determinant of
$\frak u^k(\eta_0)$ as a vector bundle is $\scrO_C(-i(k)p_0)$, where $i(k)$ is
the number of $\beta \in \widetilde \Delta$ such that $g_\beta = k$.

\begin{table}
{\begin{tabular}{||c|c|c|c||}
\hline
$G$ &  $k$ & $i(k)$ & $\frak u^k(\eta_0)$ 
\\ 
\hline\hline
$A_{n-1}$ &  $1$  & $n$ & $W_k\spcheck\otimes W_{n-k}\spcheck$      \\
\hline  
$B_{n}$ & $1$ & $3$ & $W_{n-1}\spcheck\otimes \ad W_2$ \\ \hline 
$B_{n}$ &  $2$ & $n-2$  & $\bigwedge ^2W_{n-1}\spcheck$ \\  \hline
$C_{n}$ &  $1$ & $n+1$  & $\Sym ^2W_{n-1}\spcheck$ \\ \hline
$D_{n}$ & $1$ & $4$ & $W_{n-2}\spcheck\otimes ( W_2\otimes W_2\spcheck)$ \\
\hline 
$D_{n}$ &  $2$ & $n-3$  & $\bigwedge ^2W_{n-2}\spcheck$ \\ \hline
$E_6$ &  $1$ & $3$  & $W_2\spcheck\otimes W_3\spcheck\otimes
W_3\spcheck\otimes
\scrO_C(p_0) = W_2\otimes W_3\spcheck\otimes W_3\spcheck$ \\
\hline
$E_6$ &  $2$ & $3$  & $\bigwedge^2W_3\spcheck\otimes
\bigwedge^2W_3\spcheck\otimes  
\scrO_C(p_0) = W_3\otimes W_3\otimes  \scrO_C(-p_0)$
\\
\hline 
$E_6$ &  $3$ & $1$  & $W_2\spcheck$ \\ \hline 
$E_7$ &  $1$ & $2$  & $W_2\spcheck\otimes W_3\spcheck\otimes
W_4\spcheck\otimes
\scrO_C(p_0) = W_2\otimes W_3\spcheck\otimes W_4\spcheck$ \\
\hline
$E_7$ &  $2$ & $3$  & $\bigwedge^2W_3\spcheck\otimes
\bigwedge^2W_4\spcheck\otimes  
\scrO_C(p_0) = W_3\otimes \bigwedge^2W_4\otimes  \scrO_C(-p_0)$ \\ \hline
$E_7$ &  $3$ & $2$  &  $\bigwedge^3W_4\spcheck\otimes W_2\spcheck\otimes  
\scrO_C(p_0) =W_4\otimes W_2\spcheck $
\\ \hline
$E_7$ &  $4$ & $1$  & $W_3\spcheck$ \\ \hline 
$E_8$ &  $1$ & $1$  & $W_2\spcheck\otimes W_3\spcheck\otimes
W_5\spcheck\otimes
\scrO_C(p_0) = W_2\otimes W_3\spcheck\otimes W_5\spcheck$ \\  \hline
$E_8$ &  $2$ & $2$  &   $\bigwedge^2W_3\spcheck\otimes
\bigwedge^2W_5\spcheck\otimes  
\scrO_C(p_0) =\bigwedge ^2W_5\spcheck\otimes W_3$ \\ \hline 
$E_8$ &  $3$ & $2$  & $W_2\spcheck\otimes \bigwedge^3W_5\spcheck\otimes  
\scrO_C(p_0)=\bigwedge ^2W_5\otimes W_2\spcheck$ \\ \hline 
$E_8$ &  $4$ & $2$  & $W_3\spcheck\otimes \bigwedge^4W_5\spcheck\otimes  
\scrO_C(p_0)= W_5\otimes W_3\spcheck$ \\ \hline 
$E_8$ &  $5$ & $1$  & $W_2\spcheck\otimes \bigwedge^2W_3\spcheck\otimes  
\scrO_C(p_0)=W_2\spcheck\otimes W_3$ \\
\hline 
$E_8$ &  $6$ & $1$  & $W_5\spcheck$ \\ \hline
$F_4$ &  $1$ & $2$  & $W_2\spcheck\otimes \Sym^2W_3\spcheck\otimes
\scrO_C(p_0)$ \\
\hline 
$F_4$ &  $2$ & $2$  & $\Sym^2W_3\otimes \scrO_C(-p_0)$\\ \hline 
$F_4$ &  $3$ & $1$  & $W_2\spcheck$  \\ \hline 
$G_2$ &  $1$ & $2$  & $\Sym^3W_2\spcheck\otimes \scrO_C(p_0)$ \\ \hline  
$G_2$ &  $2$ & $1$  &  $\scrO_C(-p_0)$\\

\hline\hline
\end{tabular}}
\caption{The vector bundles $\frak u^k(\eta_0)$}\label{Table1}
\end{table}

The meaning of the following theorem is that, although the bundles $\frak
u^k(\eta_0)$ are typically not semistable, their cohomology almost always
looks like that of a semistable bundle.

\begin{theorem}\label{casebycase} Let $G$ be a simple and simply connected
group. If $C$ is nodal, then for  all $k$, $H^0(C;
\frak u^k(\eta_0)) = 0$. If $C$ is cuspidal, then 
$H^0(C; \frak u^k(\eta_0)) = 0$ except for the case where $G$ is of type $E_8$
and $k=1$, where $h^0(C; \frak u^1(\eta_0)) = 1$. 
\end{theorem}
\proof  The proof boils down to a very lengthy computation. First, by
the stability of
$W_k\spcheck$,
$H^0(C; W_k\spcheck)= 0$, which handles the case of $E_6, k=3$,  $E_7, k=4$,
$E_8, k=6$, $F_4, k=3$, and $G_2, k=2$. To handle the case of $A_{n-1}$ as
well as
$B_n, k=2$, $C_n$, and $D_n, k=2$, we use the fact that
$H^0(C; W_k\spcheck\otimes W_{n-k}\spcheck)=0$ by \cite[3.3]{FMW}. By the
stability of $W_a$ and $W_b$,
$H^0(C;W_a\spcheck\otimes W_b) = \Hom (W_a, W_b) = 0$ as long as $a<b$. This
handles
$E_7, k=3$, $E_8, k=4,5$. Another easy case along these lines is $E_6, k=2$,
where
$H^0(C; W_3\otimes W_3\otimes  \scrO_C(-p_0)) = \Hom(W_3\spcheck, W_3\otimes 
\scrO_C(-p_0)) = 0$ since both $W_3\spcheck$ and $W_3\otimes  \scrO_C(-p_0)$
are stable and the slope of $W_3\otimes  \scrO_C(-p_0)$ is $-2/3<
\mu(W_3\spcheck )=-1/3$. This also handles the case of $F_4, k=2$. The cases
of $B_n, k=1$,
$D_n, k=1$, and $G_2, k=1$ are handled by the following lemma:

\begin{lemma}\label{Endss} For all $k\geq 1$, the bundle $W_k\otimes
W_k\spcheck$ is semistable. If $V$ is a semistable sheaf of degree zero, then
$h^0(C; W_n\spcheck
\otimes V) =0$. Thus, for all positive $n$ and $k$, $h^0(C; W_k\otimes
W_k\spcheck\otimes W_n\spcheck) =0$.
\end{lemma}
\begin{proof} First, we may write $W_k\otimes W_k\spcheck=\scrO_C\oplus \ad
W_k$, where $\ad W_k$ has degree zero and $h^0(\ad W_k) =0$ since $W_k$ is
simple. But then $\ad W_k$ is semistable, for if it had a torsion free
subsheaf $S$ of positive degree, then we would have $h^0(\ad W_k) \geq
h^0(S)> 0$. To prove the second statement, write $V=V_0\oplus V_1$, where
$V_0$ is supported at $p_0\in C$ and $V_1$ has support not containing $p_0$.
Since $V_0$ has a filtration whose associated gradeds are isomorphic to
$\scrO_C$ and $H^0(C; W_n\spcheck)=0$, it follows that $H^0(C; W_n\spcheck
\otimes V_0) =0$ as well. Thus it suffices to show that $H^0(C; W_n\spcheck
\otimes V_1) =0$, which follows inductively from the fact that $h^0(C;V_1) =
h^0(C; V_1\otimes
\scrO_C(-p_0)) = 0$. The final statement follows from the first two.
\end{proof}

Returning to the proof of the theorem, we now handle the case of $E_6, E_7,
E_8, F_4$ and $k=1$. From the inclusion $W_2\otimes W_3\spcheck\otimes
W_{n-1}\spcheck
\subseteq W_2\otimes W_3\spcheck\otimes W_n\spcheck$, we see that it suffices
to prove:

\begin{lemma}\label{E8comps} If $C$ is nodal, then $h^0(C; W_2\otimes
W_3\spcheck\otimes W_5\spcheck) = 0$. If $C$ is cuspidal, then $h^0(C;
W_2\otimes W_3\spcheck\otimes W_4\spcheck) = 0$ and $h^0(C; W_2\otimes
W_3\spcheck\otimes W_5\spcheck) =1$.
\end{lemma}
\begin{proof} In the nodal case, a computation based on 
Corollary~\ref{tensor} shows that $h^0(C; W_2\otimes W_3\spcheck\otimes
W_5\spcheck) = 0$. 

In the cuspidal case, first consider the bundle $W_2\otimes W_3\spcheck\otimes
W_4\spcheck$ of rank $24$. It is easy to see that $H^0(\widetilde C;
\iota^*(W_2\otimes W_3\spcheck\otimes W_4\spcheck))$ has dimension
$23$. By Corollary~\ref{sections}, there are $24$ linear conditions that
these sections must satisfy in order to give  sections of $W_2\otimes
W_3\spcheck\otimes W_4\spcheck$. One can check using Lemma~\ref{cuspcomps} and
Corollary~\ref{linalg} that there is no nonzero section satisfying these
conditions.  In the case of
$W_2\otimes W_3\spcheck\otimes W_5\spcheck$, it is likewise easy to check that
$H^0(\widetilde C; \iota^*(W_2\otimes W_3\spcheck\otimes W_5\spcheck))$ has
dimension
$30$, which is the rank of $W_2\otimes W_3\spcheck\otimes W_5\spcheck$. In
this case, although there are $30$ equations in $30$ unknowns, it turns out
that there is a $1$-dimensional space of solutions.  We leave these involved
computations to the reader.
\end{proof}

We will give a different proof of Lemma~\ref{E8comps} in the next section. 

There remain the cases of $E_7, k=2$ and $E_8, k=2,3$. These can be computed
on an \emph{ad hoc} basis, again by using Lemma~\ref{cuspcomps} and
Corollary~\ref{linalg}. The result is summarized in the following lemma:

\begin{lemma} Let $C$ be either nodal or cuspidal. Then  
$$H^0(C; W_3\otimes \bigwedge^2W_4\otimes  \scrO_C(-p_0)) = H^0(C;\bigwedge
^2W_5\spcheck\otimes W_3) = H^0(C; \bigwedge ^2W_5\otimes W_2\spcheck)
=0.\qed$$
\end{lemma}

This completes the proof of Theorem~\ref{casebycase}. \qed

\begin{corollary}\label{KSisom} Let $C$ be a   Weierstrass
cubic. Suppose either that 
$C$ is not cuspidal or that
$G$ is not of type $E_8$. Then  the functor $\mathbf{F}$ corresponding to
$H^1(C; U(\eta_0))$ defined at the beginning of this section is representable
by an affine space
$\Aff^{r+1}$ and a pair $(\Xi_P, \Phi)$. There is a marked point  
$0\in \Aff^{r+1}$ corresponding to the trivial lift
$\eta_0\times_LP$.The embedding of $\Cee^*$ as the identity component of the
center of $L$ induces a $\Cee^*$-action on $\Aff^{r+1}$ fixing $0$ and a
compatible action  on the pair
$(\Xi_P, \Phi)$.  Let $\Xi
\to C\times
\Aff^{r+1}$ be the 
$G$-bundle
$\Xi_P\times_PG$.  If $\xi_x$ denotes the restriction of
$\Xi$ to $C\times \{x\}$, then, for all
$x\neq 0$,  $h^0(C;
\ad \xi_x) = r$ and the Kodaira-Spencer homomorphism from the tangent space
$T_x$ to
$\Aff^{r+1}$ at $x$ to $H^1(C; \ad \xi_x)$ induces an isomorphism from the
quotient of $T_x$ by the tangent space to the $\Cee^*$-orbit through $x$ to
$H^1(C; \ad \xi_x)$. 
\end{corollary}
\begin{proof} The first two statements follow from
Theorem~\ref{casebycase} and  \cite[Theorem A.2.2]{FMII}.
The final  statement follows by checking directly that the method of proof of
\cite[Corollary 4.5.2]{FMII} works in the singular case as well. 
\end{proof}

\begin{remark} In case $G$ is of type $E_8$ and $C$ is cuspidal, the
dimension of $H^1(C; \frak u(\eta_0))$ is $10$ instead of $9$. More seriously,
however, the functor $\mathbf{F}$ does not satisfy the representability
criterion of \cite[Theorem A.2.2]{FMII}, and in fact one can show that it is
not representable.
\end{remark}

By applying  \cite[Theorem A.2.6 and Remark A.2.7]{FMII} to the
conformal bundle $\hat \eta_0$ constructed in \S\ref{conformbund}, we have:

\begin{corollary}\label{affinebundle} Suppose that $\pi\colon Z \to B$ is a 
Weierstrass fibration and either that $G\neq E_8$ or that $\pi$ has no
cuspidal fibers. Let $\hat \eta_0$ be the $\hat L$-bundle of
\S\ref{conformbund}. Then there is a bundle of affine spaces $\mathcal{A}\to
B$ whose fibers have dimension $r+1$ and a pair $(\widehat \Xi_{\hat P},
\Phi)$ consisting of a
$\hat P$-bundle $\widehat \Xi_{\hat P}$ over $Z\times_B\mathcal{A}$ and an
isomorphism
$\Phi\colon
\widehat \Xi_{\hat P}/U \to \pi_1^*\hat \eta_0$, with the property that the
set
$H^1(Z; U(\hat \eta_0))$ of liftings of $\eta_0$ to a $\hat P$-bundle is
identified with the set of sections of $\mathcal{A}$. The   bundle
$\mathcal{A}$ has a zero-section corresponding to the trivial lift
$\pi_1^*\hat
\eta_0\times _{\hat L}\hat P$. The group
$\Cee^*$ acts on
$\mathcal{A}$, covering the identity on $B$ and fixing the zero section, and
the quotient of the complement of the zero section is a bundle of weighted
projective spaces
$\mathcal{W}\mathcal{P}(G)$. If
$B=\Spec R$ is an affine space $\Aff^k$, then $\mathcal{A} \cong \Spec R[z_0,
\dots, z_r]$ as spaces over $\Spec R$. Finally, in case
$G\neq E_8$, the bundles $\mathcal{A}\to B$ form a space over the
moduli stack. \qed
\end{corollary}

There is also a universal conformal bundle $\mathcal{W}\mathcal{P}(G)$  in an
appropriate sense. For each group $G$, we have defined a subgroup $F$ of the
center of $G$ in Theorem~\ref{conf1}. We also have the subgroup
$\Zee/n_\alpha\Zee$ defined in
\cite[Lemma 1.2.4]{FMII}. The group $F\cdot \Zee/n_\alpha\Zee$
is always cyclic, since in case $G$ is of type $D_{2n}$, $n_\alpha =1$ by 
\cite[Lemma 3.2.2]{FMII}. Thus there is an embedding of $F\cdot
\Zee/n_\alpha\Zee$ in $\Cee^*$ and we can form the group
$\widetilde G = G\times _{F\cdot \Zee/n_\alpha\Zee}\Cee^*$. We then have:

\begin{theorem}\label{families} Let $Z\to B$ be a Weierstrass fibration and
suppose either that $G\neq E_8$ or that $\pi$ has no cuspidal fibers. For
each group
$G$, let $\mathcal{W}\mathcal{P}(G)$ be the bundle of weighted projective spaces constructed
above and let $\mathcal{W}\mathcal{P}(G)_{\rm reg}$ be the locus where the projection 
$\mathcal{W}\mathcal{P}(G)\to B$ is a smooth morphism. Then there  is a
$\widetilde G$-bundle $\widetilde \Xi
\to Z\times _B\mathcal{W}\mathcal{P}(G)_{\rm reg}$ which reduces on every fiber $C$ to the
universal conformal bundle over $C\times\WP_{\rm reg}$ constructed in
\rm{\cite[Prop.\ 4.3.2]{FMII}}.\qed
\end{theorem}

Thus, sections of $\mathcal{W}\mathcal{P}(G)_{\rm reg}\to B$ define 
$\widetilde G$-bundles over
$Z$, which can be reduced to $G$-bundles over every fiber. However, the
bundles over $Z$ do not necessarily reduce to  $G$-bundles.

\subsection{The case of $E_8$ revisited}\label{E8revisited}

We will give another argument for Lemma~\ref{E8comps}, which will be used
later, and will also handle the weight one cohomology for $E_6$ and $E_7$.
The idea of the proof will be to work on the moduli stack, or equivalently
over $\Aff^2$ with the universal family $\mathcal{E}$, and to give an
inductive computation of the relative cohomology of $\mathcal{W}_2\otimes 
\mathcal{W}_3\spcheck
\otimes\mathcal{W}_n\spcheck$ for $n\leq 5$. We will stick to the
case of a general
$\pi\colon Z
\to B$, using the universal family as needed. The
$L$-bundle $\eta_0$ corresponds to the three vector bundles
$(\mathcal{W}_2\otimes
\pi^*\mathcal{L}^7)\spcheck$, $(\mathcal{W}_3\otimes
\pi^*\mathcal{L}^4)\spcheck$, $(\mathcal{W}_5\otimes
\pi^*\mathcal{L})\spcheck$. We wish to analyze the sheaves
$R^i\pi_*(\mathcal{W}_2\otimes \mathcal{W}_3\spcheck
\otimes\mathcal{W}_5\spcheck\otimes
\pi^*\mathcal{L}^{2})$
for $i=0,1$. The factor $\pi^*\mathcal{L}^{2}$ will not however be important.
We  use the exact sequence
$$0 \to \mathcal{W}_4\spcheck \to \mathcal{W}_5\spcheck \to
\pi^*\mathcal{L}^{-4}
\to 0.$$ Thus there is an exact sequence
$$0 \to \mathcal{W}_2\otimes \mathcal{W}_3\spcheck
\otimes\mathcal{W}_4\spcheck\otimes
\pi^*\mathcal{L}^{2} \to \mathcal{W}_2\otimes \mathcal{W}_3\spcheck
\otimes\mathcal{W}_5\spcheck\otimes
\pi^*\mathcal{L}^{2}\to \mathcal{W}_2\otimes \mathcal{W}_3\spcheck\otimes
\pi^*\mathcal{L}^{-2}\to 0.$$
 The idea will be to apply $R^i\pi_*$ to the above exact
sequence and analyze the coboundary map. We begin with a series of lemmas.

\begin{lemma}\label{somemorecomps} $R^0\pi_*(\mathcal{W}_n\otimes 
\mathcal{W}_{n+1}\spcheck)
\cong \scrO_B$ and $R^1\pi_*(\mathcal{W}_n\otimes 
\mathcal{W}_{n+1}\spcheck)=0$.
\end{lemma}
\begin{proof} By
stability, $\Hom (W_n, W_{n+1}) =0$ and hence $H^1(C; W_n\otimes W_{n+1}) =0$
as well. Thus, there is a unique nonzero homomorphism $W_{n+1} \to W_n$ mod
scalars, and clearly it is the surjection of the inductive construction
realizing $W_{n+1}$ as an extension of $W_n$ by $\scrO_C$. Thus the
surjection $\mathcal{W}_{n+1} \to \mathcal{W}_n$ defines an everywhere
generating section of $R^0\pi_*Hom (\mathcal{W}_{n+1}, \mathcal{W}_n) =
R^0\pi_*(\mathcal{W}_n\otimes  \mathcal{W}_{n+1}\spcheck)$. Hence
$R^0\pi_*(\mathcal{W}_n\otimes  \mathcal{W}_{n+1}\spcheck) \cong \scrO_B$.
\end{proof}

\begin{lemma}\label{stillmorecomps} $R^0\pi_*(\mathcal{W}_2\otimes
\mathcal{W}_3\spcheck \otimes \mathcal{W}_2\spcheck) = 0$ and
$$R^1\pi_*(\mathcal{W}_2\otimes
\mathcal{W}_3\spcheck \otimes \mathcal{W}_2\spcheck) \cong
\mathcal{L}^2\oplus \scrO_B\oplus \mathcal{L}^{-2}\oplus \mathcal{L}^{-3}.$$
\end{lemma}
\begin{proof} We shall apply $R^i\pi_*$ to the exact sequence
$$0 \to \mathcal{W}_2\otimes \mathcal{W}_3\spcheck \otimes \scrO_Z(-\sigma)
\to \mathcal{W}_2\otimes \mathcal{W}_3\spcheck \otimes \mathcal{W}_2\spcheck
\to \mathcal{W}_2\otimes \mathcal{W}_3\spcheck \otimes \pi^*\mathcal{L}^{-1}
\to 0.$$
Let $\delta \colon R^0\pi_*(\mathcal{W}_2\otimes \mathcal{W}_3\spcheck \otimes
\pi^*\mathcal{L}^{-1}) \to R^1\pi_*(\mathcal{W}_2\otimes \mathcal{W}_3\spcheck \otimes
\scrO_Z(-\sigma))$ be the connecting homomorphism.
Using the fact that $\mathcal{W}_2\otimes\scrO_Z(-\sigma)\cong
\mathcal{W}_2\spcheck\otimes \mathcal{L}$ and \cite[Corollary 4.7]{FMW}, we
see that
\begin{align*}
R^0\pi_*(\mathcal{W}_2\otimes \mathcal{W}_3\spcheck \otimes
\scrO_Z(-\sigma))&=0;\\
R^1\pi_*(\mathcal{W}_2\otimes \mathcal{W}_3\spcheck \otimes
\scrO_Z(-\sigma))&= \mathcal{L}^2\oplus \scrO_B\oplus \mathcal{L}^{-1} \oplus
\mathcal{L}^{-2} \oplus \mathcal{L}^{-3}.
\end{align*}
Likewise, from Lemma~\ref{somemorecomps}, it follows that
$R^0\pi_*(\mathcal{W}_2\otimes \mathcal{W}_3\spcheck \otimes
\pi^*\mathcal{L}^{-1}) = \mathcal{L}^{-1}$ and that
$R^1\pi_*(\mathcal{W}_2\otimes \mathcal{W}_3\spcheck \otimes
\pi^*\mathcal{L}^{-1}) =0$. Thus  $R^0\pi_*(\mathcal{W}_2\otimes
\mathcal{W}_3\spcheck \otimes \mathcal{W}_2\spcheck) =\Ker\delta$ and
$R^1\pi_*(\mathcal{W}_2\otimes
\mathcal{W}_3\spcheck \otimes \mathcal{W}_2\spcheck)= \Coker\delta$. On the
other hand, $\delta$ defines a section of $H^0(B; \mathcal{L}^{3}\oplus
\mathcal{L} \oplus \scrO_B \oplus \mathcal{L}^{-1}\oplus \mathcal{L}^{-2})$,
and clearly this section is well-defined over the moduli stack.  By
Proposition~\ref{sectionsoverstack}, the only possibility is that the
component of $\delta$ in $H^0(B;\scrO_B)$ is a constant, possibly zero, and
that all other components are zero. If $\delta =0$, then 
$R^1\pi_*(\mathcal{W}_2\otimes
\mathcal{W}_3\spcheck \otimes \mathcal{W}_2\spcheck)$ is a rank four vector
bundle, implying that for every $C$, $h^1(C; W_2\otimes W_3\spcheck
\otimes W_2\spcheck) =4$ and thus that $h^0(C; W_2\otimes W_3\spcheck
\otimes W_2\spcheck) =1$ by Riemann-Roch. On the other hand, for a smooth
elliptic curve
$C$, we know that
$H^0(C; W_2\otimes W_3\spcheck
\otimes W_2\spcheck) =0$ by semistability. Thus $\delta$ defines an
isomorphism from
$\mathcal{L}^{-1}$ to the corresponding factor of $\mathcal{L}^2\oplus 
\scrO_B\oplus \mathcal{L}^{-1} \oplus
\mathcal{L}^{-2} \oplus \mathcal{L}^{-3}$. This completes the proof of
Lemma~\ref{stillmorecomps}. 
\end{proof}

The next lemma deals with the weight one cohomology for $E_6$:

\begin{lemma}\label{E6} For all Weierstrass cubics $C$, $H^0(C; W_2\otimes
W_3\spcheck\otimes W_3\spcheck) =0$, and
\begin{align*}
R^0\pi_*(\mathcal{W}_2\otimes \mathcal{W}_3\spcheck
\otimes\mathcal{W}_3\spcheck) &=0;\\
R^1\pi_*(\mathcal{W}_2\otimes \mathcal{W}_3\spcheck
\otimes\mathcal{W}_3\spcheck) &\cong \mathcal{L}^{2}\oplus \scrO_B\oplus
\mathcal{L}^{-3}.
\end{align*}
\end{lemma}
\begin{proof} Applying $R^i\pi_*$ to the exact sequence
$$0 \to \mathcal{W}_2\otimes \mathcal{W}_3\spcheck
\otimes\mathcal{W}_2\spcheck \to \mathcal{W}_2\otimes \mathcal{W}_3\spcheck
\otimes\mathcal{W}_3\spcheck \to \mathcal{W}_2\otimes \mathcal{W}_3\spcheck
\otimes\pi^*\mathcal{L}^{-2} \to 0$$
and using Lemma~\ref{stillmorecomps} gives a connecting homomorphism
$\delta\colon \mathcal{L}^{-2}\to \mathcal{L}^{2} \oplus \scrO_B \oplus
\mathcal{L}^{-2} \oplus \mathcal{L}^{-3}$ such that
$\Ker \delta = R^0\pi_*(\mathcal{W}_2\otimes
\mathcal{W}_3\spcheck \otimes\mathcal{W}_3\spcheck)$ and $\Coker
\delta = R^1\pi_*(\mathcal{W}_2\otimes \mathcal{W}_3\spcheck
\otimes\mathcal{W}_3\spcheck)$. Again using the fact that $\delta$ is defined
over the moduli stack, it follows from Proposition~\ref{sectionsoverstack}
that the induced homomorphism $\mathcal{L}^{-2} \to \mathcal{L}^{2}$ is a
constant multiple of $G_2$, the homomorphism $\mathcal{L}^{-2} \to
\mathcal{L}^{-2}$ is given by a constant, possibly zero, and all other
components are zero. Arguments as in the proof of the previous lemma show
that $\delta \neq 0$. If the induced map
$\mathcal{L}^{-2}
\to
\mathcal{L}^{-2}$ were zero, then it would follow that 
$R^1\pi_*(\mathcal{W}_2\otimes \mathcal{W}_3\spcheck
\otimes\mathcal{W}_3\spcheck)$ is the direct sum of a rank three vector
bundle plus a torsion sheaf   supported along the divisor where $G_2=0$, and
of rank one along its support. By standard cohomology and base change result,
for
$C$ a smooth elliptic curve with $g_2=0$, we would have $h^1(C; W_2\otimes
W_3\spcheck\otimes W_3\spcheck)$ of rank at least four, which contradicts the
semistability of 
$W_2\otimes W_3\spcheck\otimes W_3\spcheck$. Thus the induced homomorphism 
$\mathcal{L}^{-2} \to
\mathcal{L}^{-2}$ is an isomorphism, and so the quotient is 
$\mathcal{L}^{2}\oplus \scrO_B\oplus
\mathcal{L}^{-3}$. In particular, $R^1\pi_*(\mathcal{W}_2\otimes 
\mathcal{W}_3\spcheck
\otimes\mathcal{W}_3\spcheck)$ is locally free of rank three, so that,
again by base change, $h^1(C; W_2\otimes
W_3\spcheck\otimes W_3\spcheck) = 3$. Thus, by Riemann-Roch, $h^0(C;
W_2\otimes W_3\spcheck\otimes W_3\spcheck) = 0$. This completes
the proof.
\end{proof} 

Next we deal with the weight one cohomology for $E_7$:

\begin{lemma}\label{E7} For all Weierstrass cubics $C$, $H^0(C; W_2\otimes
W_3\spcheck\otimes W_4\spcheck) =0$, and
\begin{align*} R^0\pi_*(\mathcal{W}_2\otimes \mathcal{W}_3\spcheck
\otimes\mathcal{W}_4\spcheck) &=0;\\
R^1\pi_*(\mathcal{W}_2\otimes \mathcal{W}_3\spcheck
\otimes\mathcal{W}_4\spcheck) &\cong \mathcal{L}^{2}\oplus \scrO_B.
\end{align*}
\end{lemma}
\begin{proof} In this case, we apply $R^i\pi_*$ to the exact sequence
$$0 \to \mathcal{W}_2\otimes \mathcal{W}_3\spcheck
\otimes\mathcal{W}_3\spcheck \to \mathcal{W}_2\otimes \mathcal{W}_3\spcheck
\otimes\mathcal{W}_4\spcheck \to \mathcal{W}_2\otimes \mathcal{W}_3\spcheck
\otimes\pi^*\mathcal{L}^{-3} \to 0.$$
The result is a connecting homomorphism $\delta\colon \mathcal{L}^{-3} \to 
\mathcal{L}^{2}\oplus \scrO_B\oplus
\mathcal{L}^{-3}$. The only possible nonzero term is the induced map
$\mathcal{L}^{-3} \to \mathcal{L}^{-3}$, and we can argue as before.
\end{proof} 

Finally we turn to the case of $E_8$, and consider the exact sequence
$$0 \to \mathcal{W}_2\otimes \mathcal{W}_3\spcheck
\otimes\mathcal{W}_4\spcheck\to \mathcal{W}_2\otimes \mathcal{W}_3\spcheck
\otimes\mathcal{W}_5\spcheck\to \mathcal{W}_2\otimes
\mathcal{W}_3\spcheck\otimes \pi^*\mathcal{L}^{-4}\to 0.$$

\begin{lemma}\label{E8} Let $\delta\colon R^0\pi_*(\mathcal{W}_2\otimes
\mathcal{W}_3\spcheck\otimes \pi^*\mathcal{L}^{-4}) \to R^1\pi_*
(\mathcal{W}_2\otimes \mathcal{W}_3\spcheck
\otimes\mathcal{W}_4\spcheck)$ be the connecting homomorphism in the above
exact sequence. Then we may identify $\delta$ with a homomorphism
$\mathcal{L}^{-4}\to \mathcal{L}^2\oplus \scrO_B$. Under this identification,
$\delta = (c_1G_3, c_2G_2)$ with both $c_1$ and $c_2$ nonzero constants. It
follows that if $C$ is either smooth or nodal, then $H^0(C; W_2\otimes
W_3\spcheck\otimes W_5\spcheck) =0$, but if $C$ is cuspidal, then $H^0(C; 
W_2\otimes W_3\spcheck\otimes W_5\spcheck)$ has dimension one.
\end{lemma}
\begin{proof} We work on $\Aff^2$ with the
family $\mathcal{E}$. As before, we can identify $\delta$ with a  homomorphism
$\mathcal{L}^{-4}\to \mathcal{L}^2\oplus \scrO_{\Aff^2}$, which thus has the
form 
$\delta = (c_1g_3, c_2g_2)$, where the $c_i$ are constant. As before, if one
of the $c_i$ is zero, we reach a contradiction by considering the cohomology
for a smooth elliptic curve with either $g_2$ or $g_3$ equal to zero. Thus
both $c_1$ and $c_2$ are nonzero. Hence, 
$R^1\pi_*(\mathcal{W}_2\otimes 
\mathcal{W}_3\spcheck
\otimes\mathcal{W}_5\spcheck)\cong \mathcal{L}^{6}\otimes \frak m_0$, where
$\frak m_0$ is the maximal ideal of $0\in \Aff^2$. Via base change, for
a cuspidal curve $C$, this leads to an isomorphism $\frak m_0/\frak m_0^2 \to
H^1(C; W_2\otimes W_3\spcheck\otimes W_5\spcheck)$, and hence
$h^1(C;W_2\otimes W_3\spcheck\otimes W_5\spcheck) =2$ if $C$ is cuspidal. If
$C$ is nodal or smooth, the same argument shows that $h^1(C;W_2\otimes
W_3\spcheck\otimes W_5\spcheck) =1$. The final statement of Lemma~\ref{E8}
then follows from Riemann-Roch.
\end{proof}

\subsection{Unipotent cohomology on the normalization}

We now analyze the pullback $\iota^*\eta_0$ as a bundle on
$\widetilde C \cong \Pee^1$. Recall that, by Definition~\ref{Pee}, given
$\lambda
\in
\Lambda$, there is the corresponding $H$-bundle
$\gamma_\lambda$ over $\widetilde C$. Then:

\begin{lemma}\label{lambda1} $\iota^*\eta_0 =
\gamma_{-\lambda_1(\alpha)\spcheck}\times _HL$, where 
$\lambda_1(\alpha)\spcheck=\sum _{\beta \in \Delta}\beta\spcheck$ is the
highest coroot such that the coefficient of
$\alpha\spcheck$ in $\lambda_1(\alpha)\spcheck$ is one.
\end{lemma}
\begin{proof} First, as usual we view the maximal torus $D$ of diagonal
matrices of $GL_k(\Cee)$ as the quotient $\Cee^k/\Zee^k$, where $e_1, \dots,
e_k$ is a basis for $\Zee^k$, in such a way that the coroots of $SL_k(\Cee)$
are given by $\alpha _i = e_i-e_{i+1}$. Note that $\langle e_k,
e_{k-1}-e_k\rangle = -1$, so that $e_k$ projects to the element
$-\varpi_{\alpha_k}$. Under the correspondence between
$D$-bundles over
$\Pee^1$ and elements of $\Zee^k$, the element $-e_1$ corresponds to
$\scrO_{\Pee^1}(-1) \oplus \scrO_{\Pee^1}^{k-1}$, and so to
$\iota^*W_k\spcheck$.

Consider the morphism $H\to L \to GL_{n_i}(\Cee)$ induced by projection of
$L$ onto one of its factors. If we label the simple roots corresponding to the
factor $SL_{n_i}(\Cee)$ as $\beta_1, \dots, \beta_{n_i-1}$, so that
$\beta_{n_i-1}$ is the unique simple root not orthogonal to $\alpha$, then
$\alpha\spcheck$ projects to $-\varpi_{\beta_{n_i-1}}\spcheck$, by \cite[Lemma
1.2.3]{FMII}. Thus in the description of the maximal torus of $GL_{n_i}(\Cee)$
above, possibly after reordering the $\beta_i$, we may assume that
$\beta_1\spcheck,
\dots, \beta_{n_i-1}\spcheck, \alpha\spcheck$ correspond to $(e_1-e_2),\dots ,
(e_{n_i-1} - e_{n_i}) ,e_{n_i} $ and so $\beta_1\spcheck + 
\cdots +  \beta_{n_i-1}\spcheck + \alpha\spcheck$ corresponds to
$(e_1-e_2)+\cdots + (e_{n_i-1} - e_{n_i}) + e_{n_i} =e_1$. On the other
hand, $\lambda_1(\alpha)\spcheck$ is the sum of all the simple coroots,
and the only simple coroots which are not annihilated by the roots in
$GL_{n_i}(\Cee)$ under the projection are the coroots $\beta_1\spcheck, 
\dots, \beta_{n_i-1}\spcheck, \alpha\spcheck$. Hence 
$-\lambda_1(\alpha)\spcheck$ projects to the element $-e_1$. It follows that
$\gamma_{-\lambda_1(\alpha)\spcheck}\times _H\Cee^{n_i} =
\iota^*W_{n_i}\spcheck$ for every $i$, which proves the lemma.
\end{proof}

\begin{lemma}\label{more} Let $\beta$ be a positive root. Then $\frak
g^\beta(\gamma_{-\lambda_1(\alpha)\spcheck}) \cong \scrO_{\Pee^1}(-2)$ if and
only if
$\beta =\lambda_1(\alpha)$. In all other cases $\frak
g^\beta(\gamma_{-\lambda_1(\alpha)\spcheck})$ is isomorphic either to
$\scrO_{\Pee^1}$ or to
$\scrO_{\Pee^1}(\pm1)$.
\end{lemma}
\begin{proof} Clearly $\frak g^\beta(\gamma_{-\lambda_1(\alpha)\spcheck})
\cong
\scrO_{\Pee^1}(-\beta(\lambda_1(\alpha)\spcheck) = \scrO_{\Pee^1}(-n(\beta,
\lambda_1(\alpha)))$, where $n(\beta, \gamma)$ denotes the Cartan integer. The
result is then immediate, bearing in mind that $\alpha$ is a long root and
hence so is $\lambda_1(\alpha)$ since, by \cite[Lemma 1.4.6]{FMII}, it is Weyl
conjugate to
$\sigma_1(\alpha) = \alpha$.
\end{proof}

\begin{corollary}\label{onedim} $H^1(\Pee^1; \frak u^k(\iota^*\eta_0)) =0$ for
all
$k>1$, and 
$H^1(\Pee^1; \frak u^1(\iota^*\eta_0)) \cong \Cee$.
\end{corollary}
\begin{proof}
First note that $\frak g^{\lambda_1(\alpha)}\subseteq \frak u^1$. Since
$H^1(\Pee^1;\scrO_{\Pee^1}(a)) = 0$ for all $a\geq -1$ and
$H^1(\Pee^1;\scrO_{\Pee^1}(-2)) =\Cee$, the result is immediate from
Lemma~\ref{lambda1} and Lemma~\ref{more}.
\end{proof}

Let $\widetilde{\mathbf{F}}$ be the functor from schemes over $\Cee$ to
sets corresponding to
$H^1(\Pee^1; U(\iota^*\eta_0))$: for a scheme $S$, $\widetilde{\mathbf{F}}(S)$
is the set of
isomorphism classes of pairs $(\widetilde{\Xi}_P, \Phi)$, where
$\widetilde{\Xi}_P$ is a 
$P$-bundle over $\widetilde C\times S$ and $\Phi\colon \widetilde{\Xi}_P/U \to
\pi_1^*(\iota^*\eta_0)$ is an isomorphism of $L$-bundles. If $\mathbf{F}$
is the corresponding functor for $P$-bundles over $C$ defined at the
beginning of this section, then the map $(\Xi_P, \Phi)\mapsto ((\iota\times
\Id)^*\Xi_P, (\iota\times
\Id)^*\Phi)$ defines  a natural transformation of functors from
$\mathbf{F}$ to $\widetilde{\mathbf{F}}$.

\begin{proposition}\label{reponnorm} The functor $\widetilde{\mathbf{F}}$ is
representable by an affine space  
$\Aff^1$ and a pair over $\widetilde C \times \Aff^1$. The inclusion of
$\Cee^*$ in
$\Aut_L(\iota^*\eta_0)$ defined by the homomorphism  $\varphi_\alpha$ from
$\Cee^*$ to the center of $L$ given in {\rm \cite[Lemma 1.2.3]{FMII}} defines
an action of
$\Cee^*$ on
$\Aff^1$, which is isomorphic to the linear action  with weight
$n_\alpha$. If the functor corresponding to
$H^1(C; U(\eta_0))$ is also representable by an affine space $\Aff^{r+1}$,
then pullback defines an
$\Cee^*$-equivariant morphism $f\colon \Aff^{r+1}\to\Aff^1$.
\end{proposition}
\begin{proof} The representability of the functor by an affine space follows
from
\cite[Theorem A.2.2 and Remark A.2.5]{FMII} and the inclusion of $\Cee^*$ into
the center of
$L$ gives an action on $\Cee^*$ on the corresponding affine space. Arguing as
in
\cite[\S 4.2]{FMII}, there is a
$\Cee^*$-equivariant isomorphism from this affine space to $H^1(\Pee^1; \frak
u(\iota^*\eta_0))$. By the preceding corollary, this last vector space is
one-dimensional, and the action of $\Cee^*$ on it is given by
$\lambda_1(\alpha)\circ
\varphi_\alpha$, where we view $\lambda_1(\alpha)$ as a character on $H$.
Since
$\varphi_\alpha$ is in the kernel of all simple roots except $\alpha$, 
$\lambda_1(\alpha)\circ \varphi_\alpha (t) = \alpha \circ \varphi_\alpha (t) =
t^{n_\alpha}$.
The last statement then follows since both functors are representable.
\end{proof}

\begin{lemma}\label{lemma435} Let  $\widetilde\xi_x$ be the $G$-bundle
corresponding to  
$x\in H^1(\Pee^1; U(\iota^*\eta_0))$. Then $x\neq 0$ if and only if
$\widetilde
\xi_x$ is the trivial $G$-bundle.
\end{lemma}
\begin{proof} By Lemmas~\ref{lambda1} and \ref{more}, the bundle
$\iota^*\eta_0\times _LG$ is not the trivial
$G$-bundle, and hence the origin in
$H^1(\Pee^1; U(\iota^*\eta_0))$ does not correspond to the trivial $G$-bundle.
Since all other points of $H^1(\Pee^1; U(\iota^*\eta_0))$ are
$\Cee^*$-equivalent, it will suffice to show that there is one point $x\in
H^1(\Pee^1; U(\iota^*\eta_0))$ such that the corresponding $P$-bundle
$\widetilde\xi$ satisfies: 
$\widetilde \xi \times _PG$ is trivial. Since all $G$-bundles are
topologically trivial, there is a small deformation of $\iota^*\eta_0\times
_LG$ which is the trivial bundle. Thus it will suffice to show that the
Kodaira-Spencer map of the family corresponding to
$H^1(\Pee^1; U(\iota^*\eta_0))$ is an isomorphism. This map is given by a
homomorphism from
$H^1(\Pee^1; \frak u(\iota^*\eta_0))$ to $H^1(\Pee^1; \ad _G
(\iota^*\eta_0\times _LG))$, which is just the inclusion of $H^1(\Pee^1; \frak
u(\iota^*\eta_0))$ in 
$$H^1(\Pee^1; \ad _G (\iota^*\eta_0\times _LG)) = H^1(\Pee^1; \frak
u(\iota^*\eta_0)) \oplus H^1(\Pee^1; \ad _L (\iota^*\eta_0)) \oplus
H^1(\Pee^1; \frak u(\iota^*\eta_0)\spcheck).$$ Since $\frak
u(\iota^*\eta_0)\spcheck$ is a direct sum of line bundles of degrees $0, \pm
1,$ or
$2$,  $H^1(\Pee^1;
\frak u(\iota^*\eta_0)\spcheck) =0$. Moreover, $H^1(\Pee^1; \ad _L
(\iota^*\eta_0)) =0$ since $\iota^*\eta_0$ is rigid ($\scrO_{\Pee^1}(-1)
\oplus
\scrO_{\Pee^1}^{k-1}$ has no nontrivial deformations). Hence the
Kodaira-Spencer map is an isomorphism.
\end{proof} 

It is easy to give a direct proof of Lemma~\ref{lemma435} by working in the
subgroup of $G$ isomorphic to $SL_2(\Cee)$ whose Lie algebra is $\frak
g^{-\lambda_1(\alpha)} \oplus \Cee \cdot \lambda_1(\alpha)\spcheck \oplus 
\frak g^{\lambda_1(\alpha)}$ and using the reduction of structure of
$\iota^*\eta_0$ to a Cartan subgroup of this $SL_2(\Cee)$.

\begin{corollary}\label{Dinfty} 
Suppose either that $C$ is not cuspidal or that $G$ is not of type $E_8$, so
that the functor corresponding to $H^1(C; U(\eta_0))$ is representable by an
affine space $\Aff^{r+1}$. Let $\WP(G)$ denote the corresponding weighted
projective space $\Aff^{r+1}-\{0\}/\Cee^*$. For $p\in \WP(G)$, let
$\xi_p$ denote the corresponding $G$-bundle, and define
$$D_\infty = \{p\in \WP(G): \iota^*\xi_p \textrm{\ is not isomorphic to the
trivial $G$-bundle}\}.$$ Then
$D_\infty$ is an irreducible Weil divisor on
$\WP(G)$,isomorphic to a weighted projective subspace, and
$\WP(G) -D_\infty$ is an affine space $\Aff^r$.
\end{corollary}
\begin{proof} Let $f\colon \Aff^{r+1} \to \Aff^1$ be the
$\Cee^*$-equivariant morphism of Proposition~\ref{reponnorm}. Since
$\Cee^*$ acts on $\Aff^1$ with weight $n_\alpha$, and $n_\alpha$ divides all
of the weights of the $\Cee^*$-action on $\Aff^{r+1}$, it follows that $f$ is
linear with respect to every linear structure on $\Aff^{r+1}$ for which the
$\Cee^*$-action is linearized. Moreover, $f^{-1}(0)$ is a linear,
$\Cee^*$-invariant subspace, and the quotient $(f^{-1}(0) -\{0\})/\Cee^*$ is
a weighted projective subspace of $\WP(G)$ which is exactly the
hypersurface $D_\infty$. Thus $D_\infty$ is irreducible.

Fixing a nonzero
point $x\in H^1(\widetilde C; U(\iota^*\eta_0))\cong \Aff^1$,   $f^{-1}(x)$
is an affine subspace
$\Aff^r\subseteq \Aff^{r+1}$, and the induced morphism from $\Aff^r$ to
$\WP(G)$ embeds $\Aff^r$ in $\WP(G)$. Clearly, this image is the
complement of $D_\infty$.
\end{proof}

We continue to assume either that $G$ is not of type $E_8$ or that $C$ is not
cuspidal.  Fix once and for all  a nonzero
point $x\in H^1(\widetilde C; U(\iota^*\eta_0))$ and a
trivialization of the corresponding $G$-bundle $\widetilde \xi$,
where $\widetilde \xi$ is the $G$-bundle corresponding to $x$. If $\Xi \to
C\times \Aff^{r+1}$ is the universal bundle of
Corollary~\ref{KSisom}, then restriction defines a universal bundle
$\Xi\to C\times
f^{-1}(x)$ and the trivialization of $\widetilde \xi$ defines a
trivialization of
$(\iota\times \Id)^*\Xi$. By Theorems~\ref{cuspfamilies}
and \ref{nodefamilies}, there are induced morphisms $\sigma\colon
f^{-1}(x) \cong \Aff^r
\to \frak g$, in the cuspidal case, and $\sigma\colon \Aff^r \to G$, in the
nodal case. By Corollary~\ref{KSisom}, the image of $\sigma$ in either case is
contained in the set of regular elements.

\begin{defn}\label{parasectdef} We shall call either   $\sigma\colon \Aff^r
\to
\frak g$ (in the cuspidal case) or $\sigma\colon \Aff^r \to
G$ (in the nodal case)   a
\textsl{parabolic section}. Note that changing the trivialization of
$\widetilde \xi$ corresponds to composing $\sigma$ with the adjoint
action of a fixed $g\in G$.
\end{defn}

\subsection{A morphism to the relative moduli space}

Let $\pi\colon Z\to B$ be a Weierstrass fibration, and suppose either that
$G$ is not of type $E_8$ or that there are no cuspidal fibers of $\pi$. For
each fiber $C$, we have the affine space $\Aff^{r+1}$ and the corresponding
weighted projective space $\WP(G)$. Suppose that $C$ is singular. There is a
Zariski open subset $\WP(G)_0$ of $\WP(G)$ corresponding to $G$-bundles which
are trivial on $\widetilde C$. In the relative situation, we can define the
space $\mathcal{W}\mathcal{P}(G)_0\to B$, whose fibers consist of the full
weighted projective spaces over points corresponding to smooth fibers of
$\pi$, and are the affine spaces $\WP(G)_0$ over those points of $B$
corresponding to singular fibers. By the universal property of the relative
coarse moduli space
$\mathcal{M}^0_{Z/B}(G)$ defined in \S\ref{partialmodulispace}, there is a
morphism $\Psi\colon
\mathcal{W}\mathcal{P}(G)_0\to
\mathcal{M}^0_{Z/B}(G)$. Clearly $\Psi$ is compatible with base change.

\begin{theorem} The morphism $\Psi\colon\mathcal{W}\mathcal{P}(G)_0\to
\mathcal{M}^0_{Z/B}(G)$ is an isomorphism, and hence
$\mathcal{W}\mathcal{P}(G)_0$ is isomorphic to
$(Z_{\rm reg}\otimes \underline{\Lambda})/W$ as schemes over $B$.
\end{theorem}
\begin{proof} First assume that the generic fiber of $\pi$ is smooth. By
Looijenga's theorem, the restriction of $\Psi$ to the fiber above every smooth
point is an isomorphism and the differential of $\Psi$ is an isomorphism by
Corollary~\ref{KSisom}. Thus $\Psi$ is a quasi-finite degree one morphism
between two normal varieties, so that
$\Psi$ is an open embedding by Zariski's main theorem. We must show that
$\Psi$ is surjective. It suffices to check this at the singular fibers. For
such a fiber, the restriction of $\Psi$ is an \'etale injective morphism from
$\Aff^r$ to $\Aff^r$. It is an elementary fact \cite[Theorem (2.1)]{BCW} that
every such morphism is also surjective.

Let $C$ now be a singular curve; if $G$ is of type $E_8$ we also assume that
$C$ is not cuspidal. Choose some Weierstrass fibration
$Z'\to B'$ containing $C$ as a singular fiber and such that the generic fiber
is smooth; if $G$ is of type $E_8$ we also assume that no fiber is
cuspidal. For example, we can use the universal family $\mathcal{E}$ or its
restriction to $\Aff^2-\{0\}$. Applying the above and the compatibility with
base change, we see that $\Psi$ induces an isomorphism over the singular
fibers from $\Aff^r$ to 
$\mathcal{M}_C^0(G)$. The same is now true for every  Weierstrass fibration.
The last statement of the proof follows from  Theorem~\ref{Eotomes}.
\end{proof}

The above theorem implies that the parabolic construction is a
compactification of the space $(Z_{\rm reg}\otimes \underline{\Lambda})/W$.
In the next section, we will compare this compactification with other possible
compactifications.

Applying the above to the case where $Z$ is a point and the corresponding
Weierstrass cubic is singular gives:

\begin{corollary} Suppose that $C$ is nodal. Then $\WP(G)_0 \cong \WP(G)
-D_\infty \cong \Aff^r$, and the induced morphism $\Aff^r\to H/W$ is an
isomorphism. The composition of the inverse isomorphism $H/W \cong \Aff^r$
with the morphism $\sigma \colon \Aff^r\to G$ of
Definition~\ref{parasectdef} is a section of the adjoint quotient morphism
$G \to H/W$, whose image is contained in the set of regular elements. A
similar statement holds if $C$ is cuspidal and $G$ is not of type $E_8$,
provided that we replace $H/W$ by $\frak h/W$ and $G$ by $\frak g$. \qed
\end{corollary}

The content of the above corollary is that the parabolic section functions
much like the Kostant or Steinberg section \cite{Kostant, Steinberg} in the
cuspidal or nodal case. Using the $\Cee^*$-action, we
will give an elementary proof of this fact in the cuspidal case (i.e.\ one
that does not use Looijenga's theorem). It would be nice to have a direct
argument in the nodal case as well.

\subsection{Linearization of the action over the universal elliptic curve}

For the rest of this section, we shall be concerned with constructions over
the moduli stack, including the case of cuspidal fibers. Indeed, the
$\Cee^*$-action on the cuspidal curve will be crucial to the arguments. 
\textbf{Thus, throughout the remainder of this section, we shall assume that
$\boldsymbol{G}$ is not of type $\boldsymbol{E_8}$.}

We now consider the case of the universal elliptic curve $\mathcal{E} \to
\Aff^2 =\Spec \Cee[g_2, g_3]$. There is the $\Cee^*$-action on $\mathcal{E}$
of Definition~\ref{action}. We shall  construct a bundle of affine spaces
$\mathcal{A}_\mathcal{E}$ over $\Aff^2$ and a principal $\hat G$-bundle over
$\mathcal{E}\times_{\Aff^2}\mathcal{A}_\mathcal{E}$ such that the
$\Cee^*$-action lifts to an action on the bundle. There is also a section
$e\colon \Aff^2
\to \mathcal{A}_\mathcal{E}$ corresponding to the origin in every fiber.

\begin{theorem}\label{451}  Let $\lambda\colon \Cee^*\times \mathcal{E}\to
\mathcal{E}$  be the action of Definition~\ref{action}. There is a bundle of
affine spaces
$f\colon \mathcal{A}_\mathcal{E} =\Spec
\Cee[g_2, g_3, z_0, \dots, z_r]\to
\Aff^2=\Spec\Cee[g_2,g_3]$ and a universal $\hat G$-bundle $\widehat\Xi$ over
$\mathcal{E}
\times _{\Aff^2}\mathcal{A}_\mathcal{E}$, whose restriction to each fiber
$C$ is naturally identified with the affine space $H^1(C;U(\hat \eta_0))$
and the universal $\hat G$-bundle $\Xi\times _G\hat G$ over $C\times
H^1(E;U(\hat
\eta_0))$, where $\Xi$ is the $G$-bundle   constructed in
Corollary~\ref{KSisom}. Finally,  there are
$\Cee^*$-actions on
$\mathcal{A}_\mathcal{E}$ and on
$\widehat\Xi$, where the second action commutes with that of $\hat G$, such
that the projections 
$\mathcal{A}_\mathcal{E}
\to
\Aff^2$ and $\widehat\Xi \to \mathcal{E}
\times _{\Aff^2}\mathcal{A}_\mathcal{E}$ are equivariant.
\end{theorem}
\begin{proof} In Theorem~\ref{conf2}, we have
constructed an
$\hat L$-bundle
$\hat \eta_0$ and a
$\Cee^*$-linearization of the bundle
$\hat\eta_0$, which we can view as giving an
isomorphism $\lambda^*\hat\eta_0 \cong \pi_2^*\hat\eta_0$ over $\Cee^*\times
\mathcal{E}$. There is the functor $\mathbf{F}$ which associates to a scheme
$S$ over
$\Aff^2$ the set of isomorphism classes of pairs $(\xi_{\hat P},
\varphi)$, where
$\xi_{\hat P}$ is a $\hat P$-bundle over $\mathcal{E}\times _{\Aff^2}S$ and
$\varphi$ is an isomorphism from $\xi_{\hat P}/U$ to $\pi_1^*\hat\eta_0$,
and, by 
Corollary~\ref{affinebundle}, 
$\mathbf{F}$ is represented by the triple
$(\mathcal{A}_\mathcal{E},\widehat\Xi_{\hat P}, \Phi)$, where
$\mathcal{A}_\mathcal{E}$ is a bundle of affine spaces over $\Aff^2$,
$\widehat\Xi_{\hat P}$ is a $\hat P$-bundle over $\mathcal{E}\times
_{\Aff^2}\mathcal{A}_\mathcal{E}$, and $\Phi\colon \widehat\Xi_{\hat P}/U \to
\pi_1^*\hat\eta_0$ is an isomorphism. It will suffice to find the requisite
$\Cee^*$-actions on $\mathcal{A}_\mathcal{E}$ and on
$\widehat\Xi_{\hat P}$, and then to set $\widehat\Xi = \widehat\Xi_{\hat
P}\times_{\hat P}\hat G$.

Let $\tau \colon \Cee^*\times \mathcal{A}_\mathcal{E} \to \Aff^2$ be defined
by $\tau(\mu, a) = \mu\cdot f(a)$, i.e.\ $\tau =\lambda_0\circ (\Id\times
f)$. Then the morphism 
$$\tilde p\colon \mathcal{E} \times _{\Aff^2}(\Cee^*\times
\mathcal{A}_\mathcal{E},
\tau) \to \mathcal{E}\times \mathcal{A}_\mathcal{E}$$
defined by $\tilde p(e,(\mu,a)) =( \mu^{-1}\cdot e, a)$ has image contained in
$\mathcal{E}\times_{\Aff^2} \mathcal{A}_\mathcal{E}$, and we denote the
induced morphism by $p\colon \mathcal{E} \times _{\Aff^2}(\Cee^*\times
\mathcal{A}_\mathcal{E},
\tau) \to \mathcal{E}\times_{\Aff^2} \mathcal{A}_\mathcal{E}$. Clearly, up to
a permutation of the factors, $\pi_1\circ p$ is the natural projection of 
$\mathcal{E} \times _{\Aff^2}(\Cee^*\times
\mathcal{A}_\mathcal{E})$ to $\Cee^*\times \mathcal{E}$ followed by
$\lambda^{-1}$, the inverse of the action on $\mathcal{E}$. The pullback
bundle $p^*\widehat\Xi_{\hat P}$ has an induced isomorphism 
$$p^*\Phi\colon p^*\widehat\Xi_{\hat P} /U \cong  p^*\pi_1^*\hat\eta_0\cong
(\lambda^{-1})^*\eta_0 \cong \pi_1^*\hat\eta_0,$$
where we have used the given linearization of the action $\lambda$ on
$\hat\eta_0$. By abuse of notation, we denote the induced isomorphism
$p^*\widehat\Xi_{\hat P} /U \to \pi_1^*\hat\eta_0$ by $p^*\Phi$ as well.
The pair $(p^*\widehat\Xi_{\hat P}, p^*\Phi)$ defines an element of
$\mathbf{F}(\Cee^*\times \mathcal{A}_\mathcal{E},\tau)$, and hence there is a
classifying morphism $\lambda'\colon \Cee^*\times \mathcal{A}_\mathcal{E} \to
\mathcal{A}_\mathcal{E}$. Moreover, the following diagram commutes:
$$\begin{CD}
\Cee^*\times \mathcal{A}_\mathcal{E} @>{\lambda'}>> \mathcal{A}_\mathcal{E}\\
@V{\Id\times f}VV @VV{f}V\\
\Cee^*\times \Aff^2 @>{\lambda_0}>> \Aff^2.
\end{CD}$$
It is straightforward to check that $\lambda'$ defines an action of $\Cee^*$
on $\mathcal{A}_\mathcal{E}$. Since both $\lambda$ and $\lambda'$ cover the
action of $\lambda_0$ on $\Aff^2$, there is an induced action of $\Cee^*$ on
$\mathcal{E}\times_{\Aff^2}\mathcal{A}_\mathcal{E}$, which we denote by
$\lambda\times \lambda'$.

Finally, we must show that the action $\lambda\times \lambda'$ lifts to an
action on $\widehat\Xi_{\hat P}$. To do so, note that by definition $(\Id
\times \lambda')\widehat\Xi_{\hat P} = p^*\widehat\Xi_{\hat P}$, and hence
$(\lambda\times \lambda')^*\widehat\Xi_{\hat P}$ is identified with the
pullback of $\widehat\Xi_{\hat P}$ under the composition of the morphisms
$$\Cee^*\times (\mathcal{E}\times _{\Aff^2}\mathcal{A}_\mathcal{E}) \to
\mathcal{E}\times _{\Aff^2}(\Cee^*\times \mathcal{A}_\mathcal{E}, \tau) \to
\mathcal{E}\times _{\Aff^2}\mathcal{A}_\mathcal{E}$$
where the first morphism is $(\mu,(e,a)) \mapsto (\mu\cdot e,(\mu, a))$ and
the second is $(e,(\mu, a)) \mapsto (\mu^{-1}\cdot e, a)$. Thus the
composition is just projection onto the last two factors. Thus we have found
an isomorphism from $(\lambda\times \lambda')^*\widehat\Xi_{\hat P}$ to the
bundle over $\Cee^*\times (\mathcal{E}\times
_{\Aff^2}\mathcal{A}_\mathcal{E})$ obtained by pulling back
$\widehat\Xi_{\hat P}$ by projection onto the last two factors. Again, it is
straightforward to check that this gives an action of $\Cee^*$ on the
principal bundle $\widehat\Xi_{\hat P}$.
\end{proof}

\begin{remark} The proof actually shows that, in the above notation, the
linearization on the
$\hat G$-bundle $\widehat \Xi$ is induced from one on $\widehat\Xi_{\hat P}$
which fixes the trivialization $\Phi\colon \widehat\Xi_{\hat P}/U\to
\pi_1^*\hat\eta_0$.
\end{remark}

Thus there are two commuting actions of $\Cee^*$ on $\mathcal{A}_\mathcal{E}
= \Spec
\Cee[g_2, g_3, z_0, \dots, z_r]$. Here   the first  is
the action of the identity component of the center of
$L$, and this action is equivariant with respect to the trivial action of
$\Cee^*$ on $\Spec
\Cee[g_2, g_3]$. The other action is the lifting of the $\Cee^*$-action
$\lambda$ on $\mathcal{E}$ given by Theorem~\ref{451} and is equivariant with
respect to the  action on
$\Spec \Cee[g_2, g_3]$ given by
$\mu\cdot (g_2, g_3) = (\mu^4g_2, \mu^6g_3)$. The two commuting actions of
$\Cee^*$ define an action of $\Cee^*\times \Cee^*$ on
$\mathcal{A}_\mathcal{E}$, preserving the zero-section $e$ of
$\mathcal{A}_\mathcal{E}\to \mathcal{E}$. The dual action of $\Cee^*\times
\Cee^*$ on the coordinate ring
$\Cee[g_2, g_3, z_0, \dots, z_r]$ is equivalent to a bigrading on $\Cee[g_2,
g_3, z_0, \dots, z_r]$. 

\begin{proposition}\label{linearization} Let
$\mathcal{T}=e^*T_{\mathcal{A}_\mathcal{E}/\Aff^2}$ be the restriction of the
relative tangent bundle to the zero section $e$ . Since $e$ is
preserved by 
$\Cee^*\times \Cee^*$, there is an induced action of $\Cee^*\times
\Cee^*$  on $\mathcal{T}$. There is a
$\Cee^*\times \Cee^*$-equivariant isomorphism of spaces over $\Aff^2$ from
$\mathcal{A}_\mathcal{E}$ to $\mathcal{T}$.
\end{proposition}
\begin{proof} The argument is very similar to that for \cite[Lemma
4.2.3]{FMII}. First, $\mathcal{A}_\mathcal{E} =\Spec \Cee[g_2, g_3, z_0,
\dots, z_r]$, and the construction of the Appendix to \cite{FMII} shows that
we can assume that the zero-section
$e$ is defined by $z_0 =\cdots = z_r=0$. Let $M$ be the projective $\Cee[g_2,
g_3]$-module corresponding to $\mathcal{T}^*$. The $(\Cee^*\times
\Cee^*)$-action on $\mathcal{T}$ defines a bigrading $M= \bigoplus
_{p,q}M_{p,q}$, and the homomorphism $\Cee[g_2, g_3, z_0,
\dots, z_r]\to M$ defined by sending a function to its differential along $e$
is a surjective homomorphism of bigraded $\Cee[g_2, g_3]$-modules. If
$f\in \Cee[g_2, g_3]$ is weighted homogeneous of degree $n$, then $f\cdot
M_{p,q} \subseteq M_{p, q+n}$. Let $M(p) = \bigoplus _qM_{p,q}$. Then $M(p)$
is a $\Cee[g_2,
g_3]$-module and  $M =\bigoplus_pM(p)$, where only finitely many summands are
nonzero. Each
$M(p)$ is a projective graded $\Cee[g_2, g_3]$-module, compatibly with the
grading on
$\Cee[g_2, g_3]$. Thus, by Lemma~\ref{graded}, each $M(p)$ is a free
$\Cee[g_2, g_3]$-module, and moreover there is a homogeneous $\Cee[g_2,
g_3]$-basis for $M(p)$. It follows that $M$ is free and that there exists a
bihomogeneous basis for
$M$. Lifting these to bihomogeneous elements of $\Cee[g_2, g_3, z_0, \dots,
z_r]$ defines a homomorphism of bigraded $\Cee[g_2, g_3]$-algebras from
$\Sym_{\Cee[g_2, g_3]}^*M$ (the symmetric algebra on the free $\Cee[g_2,
g_3]$-module $M$) to
$\Cee[g_2, g_3, z_0, \dots, z_r]$, corresponding to a $\Cee^*\times
\Cee^*$-equivariant morphism from $\mathcal{A}_\mathcal{E}$ to
$\mathbb{V}(\mathcal{T})$. The arguments of \cite[Lemma
4.2.3]{FMII} show that this morphism is an isomorphism over every point of
$\Aff^2$, and hence it is an isomorphism.
\end{proof}

We have thus linearized the action of $\Cee^*\times
\Cee^*$ on $\mathcal{A}_\mathcal{E}$. The first $\Cee^*$ acts as a group of
vector bundle isomorphisms, and hence with various weights, which we know to
be
$g_0, \dots, g_r$. The corresponding eigenspaces are subbundles on which the
second
$\Cee^*$ acts via a linearization covering the action on the base, and by
Lemma~\ref{graded} each such may be written as a direct sum of linearized line
bundles.  Since the quotient of $\mathcal{A}_\mathcal{E}$ minus the zero
section by the first $\Cee^*$-action is the weighted projective bundle
$\mathcal{W}\mathcal{P}(G)$, we have:

\begin{corollary} Let $Z\to B$ be a Weierstrass fibration, and let
$\mathcal{W}\mathcal{P}(G)\to B$ be the bundle of weighted projective
spaces corresponding to $G$ and to the special root $\alpha$. Then
$\mathcal{W}\mathcal{P}(G)$ is isomorphic to the weighted projective space
bundle associated to
$$\mathcal{L}^{-d_0} \oplus \cdots \oplus \mathcal{L}^{-d_r}$$
for some integers $d_i$, where $\Cee^*$ acts on
$\mathcal{L}^{-d_i}$ with weight $g_i$.
\qed
\end{corollary}

To determine the integers $d_i$, or equivalently the characters in the
linearizations for the family $\mathcal{A}_\mathcal{E} \to \mathcal{E}$, it
is enough to look at the induced action on the fiber of
$\mathcal{A}_\mathcal{E}$ over the unique fixed point of the
$\Cee^*$-action. We find a conceptual way
to do this in the next subsection. Note however that, as in \cite[Lemma
4.1.2]{FMII}, there is a canonical isomorphism (of $(\Cee^*\times
\Cee^*)$-linearized vector bundles)
$$\mathcal{T} \cong R^1\pi_*\frak u(\hat \eta_0),$$
and using this description the linearization can be described directly via a
case-by-case analysis.

\subsection{The line bundles and the Casimir weights}\label{4.6}

We continue to assume that $G$ is not of type $E_8$. This section has two
goals. First, we want to determine the integers
$d_i$ above and relate them to the Casimir weights, which we shall define
below. Second, we shall give a direct proof that the parabolic section
defines a Kostant-type section in the cuspidal case (without using
Looijenga's theorem for a smooth curve).

We begin with a general remark. Suppose that $B$ is a scheme and that
$\bigoplus _i\mathcal{L}_i$ is a direct sum of line bundles over $B$. Let
$\Cee^*$ act on the total space $\mathbb{V}(\bigoplus _i\mathcal{L}_i)$ by
acting on each line factor
$\mathbb{V}(\mathcal{L}_i)$ with weight $w_i$. We can then form the bundle of
weighted projective spaces $\WP(\bigoplus _i\mathcal{L}_i)$. Let $M$ be a
line bundle on $B$. Let $\Cee^*$ act on $\mathbb{V}(\bigoplus
_i(\mathcal{L}_i\otimes M^{w_i}))$, again by acting on
$\mathbb{V}(\mathcal{L}_i\otimes M^{w_i})$ with weight
$w_i$. Clearly this does not change $\WP(\bigoplus _i\mathcal{L}_i)$. More
generally, if all of the $w_i$ are divisible by some integer $n$, we can make
the same construction with
$\bigoplus _i(\mathcal{L}_i\otimes M^{w_i/n})$ and we will not change the
corresponding weighted projective space bundle. A similar construction holds
if $B$ is a scheme with a $\Cee^*$-action, the $\mathcal{L}_i$ are linearized
bundles which are trivial as line bundles, and similarly for $M$. We will
refer to the bundle $\bigoplus _i(\mathcal{L}_i\otimes M^{w_i/n})$ as the
\textsl{weighted tensor product} of the bundle $\bigoplus _i\mathcal{L}_i$ by
the line bundle $M^{w_i/n}$.

We now specialize to $B = \Aff^2$, together with the $\Cee^*$-action defined
by $\mu\cdot(g_2, g_3) =(\mu^4g_2, \mu^6g_3)$, and to the case of linearized
bundles.  We can further restrict the discussion
above to the fiber over the origin, corresponding to a cuspidal curve $C$. Let
$\mathcal{L}_0$ denote the representation
$(\Cee,
\chi_1)$. Using the linearization of Proposition~\ref{linearization},
$\Aff^{r+1}= (\mathcal{A}_{\mathcal{E}})_0$, and by Lemma~\ref{linwts} this
can be  written as a direct sum of lines $\bigoplus
_{i=0}^r\mathcal{L}_0^{a_i}$. Here
$(\nu,\mu) \in \Cee^*\times \Cee^*$ acts on $\mathcal{L}_0^{a_i}$ as
$\nu^{n_\alpha\cdot g_i}\mu^{a_i}$. We denote the action of $(1,\mu)$ on $p\in
\Aff^{r+1}$ by
$\mu\cdot p$. Each point of
$\Aff^{r+1}$ corresponds to a pair $(\xi_P, \varphi)$, where $\xi_P$ is a
$G$-bundle. Let $\xi =\xi_P\times_PG$. Purely in terms of $G$-bundles, the
$\mu$-action sends
$\xi$ to
$\mu^*\xi$, where $\mu \colon C \to C$ is the automorphism of $C$
given by the action of $\mu\in \Cee^*$ given in Definition~\ref{action}.  The
action of $\nu$ is solely on $\varphi$, and so leaves the isomorphism
class of $\xi$ unchanged.

Let $\Aff^1=H^1(\widetilde C; U(\iota^*\hat \eta_0))$. We can also linearize 
the ($\Cee^*\times
\Cee^*$)-action on $\Aff^1$. The
morphism $\pi\colon \Aff^{r+1} \to
\Aff^1$ defined by the normalization map is then  a ($\Cee^*\times
\Cee^*$)-equivariant morphism.  Note that the first factor acts with weight
$n_\alpha$, so that we may assume that $(\nu,\mu)$ acts on $\Aff^1$ as
$\nu^{n_\alpha}\mu^{a_0}$. Thus  $\Aff^1\cong \mathcal{L}_0^{a_0}$.

Fix $x\in \Aff^1-\{0\}$ and let $\Aff^r$ be  the affine subspace 
$\pi^{-1}(x)$. Multiplication by $\mu$ does not preserve this subspace.
However, if we consider the weighted tensor product $\Aff^{r+1}\otimes
\mathcal{L}_0^{-a_0}$, then $\pi$ defines an equivariant morphism from this
space to the one-dimensional vector space with trivial $\Cee^*$-action, and
$\Aff^r$ is invariant under this $\Cee^*$-action.  Another way to view this
weighted tensor product is as follows: given $\mu$, choose $\nu\in \Cee^*$
such that $\nu^{n_\alpha} = \mu^{-a_0}$. Then $\mu$ acts on $\Aff^{r+1}$ via
the action of
$(\nu, \mu)\in \Cee^*\times \Cee^*$.

In
Definition~\ref{parasectdef}, we have defined a morphism
$\sigma$ from the affine space
$\Aff^r$ to
$\frak g^{\mathrm{reg}}$, the parabolic section. It is well-defined up to
the adjoint action of an element of $G$. The adjoint quotient of $\frak g$ is
$\frak h/W$, and hence there is a well-defined morphism $\Aff^r\to \frak h/W$.
The linear action of $\Cee^*$ on $\frak h$ and on $\frak g$ induces an
action of
$\Cee^*$ on $\frak h/W$. The weights of this action are given as follows.  By
Chevalley's theorem, $(\Sym^*\frak h^*)^W$ is a polynomial algebra with
homogeneous generators of degrees $d_1, \dots, d_r$, where the $d_i$ are
the \textsl{Casimir weights}. In other words, there is a linear structure
on $\frak h/W$ for which $\Cee^*$ acts   with the
weights $d_i$.
With this said, we have the following:

\begin{theorem}\label{equivisom} The morphism $\sigma$ induces an isomorphism
$\ov \sigma \colon \Aff^r\to \frak h/W$ with the property that $\ov \sigma
(\mu\cdot x) =
\mu^{-1}\ov\sigma (x)$. Thus  the
$\Cee^*$-action on
$\Aff^r$ is linearizable, with negative weights $-d_1, \dots, -d_r$.
\end{theorem}
\proof We first show the $\Cee^*$-equivariance. Let $x\in \Aff^r$ and let
$\xi$ be the corresponding
$G$-bundle. Choose a trivialization $\tau$ of $\iota^*\xi$ and let $X\in
\frak g$ be the element associated to the pair $(\xi, \tau)$ by
Theorem~\ref{cuspbundles}. Then $\ov\sigma(x) =[X]$, the class of $X$ in the
adjoint quotient of $\frak g$. To calculate
$\ov\sigma (\mu\cdot x)$, by the above remarks, choose $\nu\in \Cee^*$
such that $\nu^{n_\alpha} = \mu^{-a_0}$ and consider the corresponding bundle
$(\nu,\mu)\cdot \xi$. As an abstract $G$-bundle, this is isomorphic to
$\mu^*\xi$, and hence corresponds by Corollary~\ref{muinverse} under the
pullback of the trivialization $\tau$ to $\mu^{-1}X$. Thus $\ov \sigma
(\mu\cdot x) =[\mu^{-1}X]= \mu^{-1}[X]= \mu^{-1}\ov\sigma (x)$.

It thus suffices to show that the morphism $\ov\sigma\colon \Aff^r\to \frak
h/W$ is an isomorphism. First we claim:

\begin{lemma} The differential of $\ov\sigma$ is an isomorphism at every
point.
\end{lemma}  
\begin{proof} By Corollary~\ref{KSisom}, the Kodaira-Spencer
homomorphism induced from the universal bundle over $C\times \Aff^r$ is an
isomorphism. Given a $G$-bundle $\xi$ corresponding to a point $x$ of
$\Aff^r$, let $X\in \frak g^{\mathrm{reg}}$ be the corresponding regular
element of $\frak g$. It is easy to check (cf.\ the proof of
Theorem~\ref{cuspbundles}) that
$H^1(C; \ad \xi)\cong\Ker(\ad X)$ and that the differential of $\sigma\colon
\Aff^r \to \frak g$ at $x$ identifies the tangent space to $\Aff^r$ at $x$
with $\Ker(\ad X )\subseteq \frak g$. On the other hand, by a result of
Kostant
\cite{Kostant}, since $X$ is regular, the differential of the adjoint
quotient morphism $\frak g \to
\frak h/W$ is surjective at $X$ and thus induces an isomorphism from  $\Ker
\ad X$ to the tangent space to $\frak h/W$ at the image of $X$. Combining
these two results, we see that the differential of $\ov\sigma$ is an
isomorphism.
\end{proof}

The theorem then follows from the more general lemma:

\begin{lemma} Suppose that $f\colon \Aff^r \to \Cee^r$ is a
$\Cee^*$-equivariant morphism, where $\Cee^*$ acts linearly on $\Cee^r$ with
strictly negative weights. If the differential of $f$ is everywhere an
isomorphism, then $f$ is an isomorphism.
\end{lemma}
\begin{proof} Directly via a topological argument, or by a theorem of
Bia\l ynicki-Birula \cite{BB}, there is a fixed point $x_0$ for the action of
$\Cee^*$ on
$\Aff^r$, and $f(x_0) =0$. Since the fibers of $f$ are finite, there are only
finitely many fixed points, say $x_0, \dots, x_N$. There exists a classical
neighborhood
$U$ of $0\in \Cee^r$ such that, for all $u\in U$ and $t\in \Cee^*, |t|\geq 1$,
$tu\in U$ and moreover such that $f^{-1}(U) =\coprod _iU_i$, a disjoint union
of open sets $U_i$ in $\Aff^r$ such that $x_i\in U_i$. Clearly, if $v\in U_i$,
then
$tv\in U_i$ for all
$t\in \Cee^*, |t|\geq 1$, and $\lim_{t\to \infty}tv = x_i$. It follows that 
$$\Cee^*\cdot U_i = \{v\in \Aff^r: \lim_{t\to \infty}tv = x_i\}$$
and that $\Aff^r$ is the disjoint union of the open sets $\Cee^*\cdot U_i$.
By connectedness, there is a unique fixed point $x_0$ and every
$\Cee^*$-orbit contains it in its closure. A standard argument (cf.\ the
proof of \cite[Lemma 4.2.3]{FMII}) then shows that
$f$ is a bijection, and hence an isomorphism.
\end{proof}

\begin{remark} One can show quite explicitly that, for appropriate choices,
the inverse of the morphism $\Aff^r\to \frak h/W$, composed with
$\sigma\colon \Aff^r \to \frak g^{\rm reg}$, can be identified with the
Kostant section constructed in
\cite{Kostant}.
\end{remark} 

As a corollary of Theorem~\ref{equivisom}, we have the corresponding result
for the moduli stack:

\begin{corollary} Suppose that $G$ is not of type $E_8$. Let $Z\to B$ be
a Weierstrass fibration, and let
$\mathcal{W}\mathcal{P}(G)\to B$ be the bundle of weighted projective
spaces corresponding to $G$ and to a special root $\alpha$. Then
$\mathcal{W}\mathcal{P}(G)$ is isomorphic to the weighted projective space
bundle associated to
$$\scrO_B\oplus \mathcal{L}^{-d_1} \oplus \cdots \oplus \mathcal{L}^{-d_r},$$
where the $d_i$ are the Casimir weights of $G$. Moreover, $\Cee^*$ acts on
$\mathcal{L}^{-d_i}$ with weight $g_i$ for some ordering of the  simple roots,
and $\Cee^*$ acts on $\scrO_B$ with weight one.
\qed
\end{corollary}

Of course, the corollary does not  determine which $g_i$ correspond
to which $d_i$. Direct computation as in \S\ref{E8revisited} shows the
following: if the
$d_i$ are arranged in non-decreasing order, then so are the $g_i$, except for
the case of $D_n$, where the ordering of the $g_i$ as $1,1,1,1,2, \dots,2$
corresponds to the ordering of the $d_i$ as $0,2,4,n, 6, \dots,2n-2$.

\section{Comparison with other constructions}

Our goal in this section is to compare the family of weighted projective 
bundles given by the parabolic construction with other compactifications of
the moduli space of semistable $G$-bundles in families. For example,
Wirthm\"uller
\cite{Wirt} has constructed a compactification in case all fibers are smooth 
or nodal by toric methods. We shall show that, under very general
circumstances, all such constructions lead to the same family of
compactifications.

 \subsection{A lemma on bundles of weighted projective spaces}

We begin with some general remarks on weighted projective spaces. Let $w_i\in
\Zee^+$ for $i=0, \dots, r$ and let $\WP$ be a  weighted projective space of type
$(w_0, \dots, w_r)$, in other words the quotient of $\Cee^{r+1}-0$ by $\Cee^*$
acting via $\lambda\cdot(z_0, \dots, z_r) = (\lambda^{w_0}z_0, \dots,
\lambda^{w_r}z_r)$. We assume that the gcd of $w_0, \dots, w_r$ is one. For every
$n\in
\Zee$, there is the torsion free, reflexive coherent sheaf
$\scrO_{\WP}(n)$, whose global sections are weighted homogeneous polynomials on
$\Cee^{r+1}$ of degree $n$. Moreover, every effective Weil divisor on
$\WP$ is the zero locus of a weighted homogeneous polynomial. However, the sheaves
$\scrO_{\WP}(n)$ can behave in a rather complicated way. For example, if $r=1$ so
that $\WP \cong
\Pee^1$, and
$w_0$ and $w_1$ are relatively prime, $\scrO_{\WP}(w_0)\cong
\scrO_{\WP}(w_1)\cong\scrO_{\Pee^1}$ whereas
$\scrO_{\WP}(w_0w_1)\cong\scrO_{\Pee^1}(1)$. Hence the map $n\in \Zee \mapsto
\scrO_{\WP}(n)$ from $\Zee$ to the class group of Weil divisors on $\WP$
need not be injective, nor a homomorphism. As a consequence, not every
automorphism of
$\WP$ need be induced from a graded automorphism of the graded ring $\Cee[z_0,
\dots, z_r]$, as is clear from the above example.

On the other hand, suppose that the singular locus of $\WP$, in the sense of
orbifolds, has codimension at least two, in other words that there is no prime
which divides $r$ of the positive integers $w_i$. Let $\WP_{\mathrm{reg}}$ be the
regular locus of $\WP$. Then since $\WP$ is normal and $\WP - \WP_{\mathrm{reg}}$
has codimension at least two in $\WP$, a Weil divisor is determined by its
restriction to $\WP_{\mathrm{reg}}$. It follows that  the map
$n\in \Zee \mapsto
\scrO_{\WP}(n)$ from $\Zee$ to the class group of Weil divisors on $\WP$ is an
isomorphism of groups, and that moreover the class group is also isomorphic to
$\Pic (\WP_{\mathrm{reg}})$ via the restriction map. In this case, if $f$ is
an automorphism of $\WP$, then $f^*\scrO_{\WP}(k) = \scrO_{\WP}(k)$ for every
$k\in \Zee$, and hence $f$ is induced from a graded automorphism of
$\Cee[z_0,
\dots, z_r]$. 

For $B$ a connected complex space or
scheme, we define a \textsl{fibration of weighted projective spaces over $B$ of
type
$(w_0, \dots, w_r)$}  to consist of a complex space $P$ and a morphism $p\colon
P\to B$ such that there exists an open cover $\{\Omega_i\}$ of $B$ and
isomorphisms (of spaces over
$\Omega_i$) $p^{-1}(\Omega_i) \to \WP\times \Omega_i$, where $\WP$ is a weighted
projective space of type
$(w_0, \dots, w_r)$.

\begin{theorem}\label{wtdprojthm} Let $B$ be a smooth scheme and let
$p_i\colon P_i\to B$ be two fibrations of weighted projective spaces over $B$
of type
$(w_0, \dots, w_r)$. Assume that no prime divides $r$ of the positive
integers $w_0, \dots, w_r$. Suppose that there exists a Zariski open and  dense
set
$U\subseteq P_1$ and a holomorphic map $f\colon U \to P_2$ of spaces over
$B$ such that 
\begin{enumerate}
\item[\rm (i)] $U$ contains an open set of the form $p_1^{-1}(V)$, where $V$ is a
Zariski dense subset of $B$;
\item[\rm (ii)] $U$ has nonempty intersection with every fiber;
\item[\rm (iii)] $f$ restricts to an isomorphism from $p_1^{-1}(V)$ to
$p_2^{-1}(V)$.
\item[\rm (iv)] For all $b\in B$, $f(U\cap p_1^{-1}(b))$ is not contained in a
proper subvariety of $p_2^{-1}(b)$.
\end{enumerate}
Then $f$ extends to a holomorphic isomorphism of fibrations of weighted projective
spaces over $B$.
\end{theorem}
\begin{proof} The hypotheses imply that $P_1$ and $P_2$ are normal, that
$\operatorname{codim}(P_1-U) \geq 2$, and that the orbifold singular locus of
$\WP$ has codimension at least two. Moreover, the question is
local around
$b\in B$, and thus we may assume that $B=\Spec R$ is affine and that both $P_1$
and $P_2$ are products $\WP\times B$. Thus $P_1=P_2=\mathbf{Proj}\,R[s_0, \dots,
s_r]$, where the weight of $s_i$ is $w_i$, and we may assume that the $w_i$ are
arranged in nondecreasing order.

For each $i$, $f^*s_i$ is a section of $\scrO_{P_1}(w_i)\otimes
\pi_2^*M_i|U$, where $M_i$ is a line bundle over $B$. After further shrinking
$B$, we may assume that $M_i$ is trivial. Since
$\scrO_{P_1}(w_i)$ is reflexive, $f^*s_i$ extends to a global holomorphic section
of
$\scrO_{P_1}(w_i)$. Clearly, to prove the theorem, it suffices to prove that, for
all $b\in B$, the sections $(f^*s_i)_b$ generate the ring $\Cee[s_0, \dots,
s_r]$. Order the $s_i$ so that $s_0, \dots, s_{a_0}$ have weight $w_0$,
$s_{a_0+1}, \dots, s_{a_1}$ have the next highest weight $w_{i_1}$, and so on. It
will suffice to prove that, for all $j$, $f^*s_{a_{j-1}+1}, \dots, f^*s_{a_j}$
span the space of all homogeneous polynomials of weight $w_{i_j}$ modulo those
which are polynomials in variables of smaller weight. It is enough to prove that 
$f^*s_{a_{j-1}+1}, \dots, f^*s_{a_j}$
are linearly independent modulo  polynomials in variables of smaller weight. The
proof is by induction on the weight $w_{i_j}$, and we can begin the induction at
weight $0$, i.e.\ for the constants, where the result is clear. Suppose that there
is a relation
$$\sum _k\lambda_kf^*s_{a_{j-1}+k} + P(s_0, \dots, s_{a_{j-1}}) =0,$$ 
where the $\lambda_i$ are not all $0$ and $P$ is homogeneous of weight $w_{i_j}$. 
Then by the inductive hypothesis we can write $P(s_0, \dots, s_{a_{j-1}}) =
f^*Q(s_0, \dots, s_{a_{j-1}})$ for some homogeneous polynomial $Q$. It then
follows that $f$, which is defined on $U$  maps the open dense subset
$U_b=p_1^{-1}(b)\cap U$ onto the proper subvariety of the weighted
projective space $p_2^{-1}(b)$ defined by $\sum _k\lambda_ks_{a_{j-1}+k} +
Q(s_0, \dots, s_{a_{j-1}}) =0$. This contradicts Condition (iv).
\end{proof}

 \subsection{Wirthm\"uller's construction}

We briefly review the construction of Wirthm\"uller \cite{Wirt}. Let $\Gamma$
be a torsion free subgroup of finite index in $SL_2(\Zee)$. Let $\frak H$
denote the upper half plane and let $B_0 = \frak H/\Gamma$ be the
corresponding affine curve, with $\mathcal{E}_\Gamma^0\to B_0$ the universal
elliptic curve over $B_0$. If
$B$ is the projective completion of
$B_0$, then there is a fibration of Weierstrass cubics
$\pi\colon \mathcal{E}_\Gamma
\to B$, whose singular fibers are nodal. Let
$(\mathcal{E}_\Gamma)_{\textrm{reg}}
\to B$ be the restriction of this fibration to the smooth locus of $\pi$,
so that $(\mathcal{E}_\Gamma)_{\textrm{reg}}$ is a group scheme over $B$. The
construction of
\cite{Wirt} is a toroidal compactification of the fibration $\frak{A}^0  =
\mathcal{E}_\Gamma^0\otimes \Lambda \to B_0$ to a fibration $\frak{A}\to
B$. Unwinding the definitions of \cite{Wirt}, one checks that $\frak{A}$
contains $(\mathcal{E}_\Gamma)_{\textrm{reg}} \otimes \Lambda$ as a dense
open set, and in fact the fibers of $\frak{A}$ over the cusps are toric
compactifications of the corresponding fibers of
$(\mathcal{E}_\Gamma)_{\textrm{reg}} \otimes \Lambda$, which are algebraic
tori $\Cee^*\otimes \Lambda$. The action of
$W$ on
$\frak{A}^0  =
\mathcal{E}_\Gamma^0\otimes \Lambda$ defined by the action of $W$ on
$\Lambda$ extends to an action of $W$ on $\frak{A}$.

Wirthm\"uller then proves:

\begin{theorem} Suppose that $G$ is not of type $E_8$. Then the  fibration
$\frak{A}/W \to B$ is a fibration of weighted projective spaces of type
$(g_0, \dots, g_r)$.
\qed
\end{theorem}

In light of Theorem~\ref{wtdprojthm}, we then have:

\begin{theorem} Suppose that $G$ is not of type $E_8$, and let
$\mathcal{W}\mathcal{P}(G)\to B$ be the fibration of weighted projective
spaces associated to the fibration
$\mathcal{E}_\Gamma \to B$ and the group $G$. Then the inclusions of
$((\mathcal{E}_\Gamma)_{\rm reg } \otimes \Lambda)/W$ in
$\mathcal{W}\mathcal{P}(G)$ and in
$\frak{A}/W$ extend to define an isomorphism of fibrations of weighted
projective spaces. \qed
\end{theorem}

Even in the case of $E_8$, the parabolic construction gives an extension of
$\frak{A}^0/W$ to a fibration of weighted projective spaces over $B$.
Thus it is likely that $\frak{A}/W$ is also a fibration of weighted
projective spaces over $B$. If so, then Theorem~\ref{wtdprojthm} would
imply that it is the same as the fibration given by the parabolic
construction. However, 
$\frak{A}/W$ is not the fibration of weighted projective spaces associated to
$\bigoplus _i\mathcal{L}^{a_i}$ as is required by the method of proof of
\cite{Wirt}.

\section{Inclusions of subgroups}

In this section, we analyze the relation between the parabolic construction for a
group $G$ and for certain subgroups $G'$. The prototype for this question
is the inclusion of $SL_{n-1}(\Cee)$ in $SL_n(\Cee)$, and the problem,
which was studied in \cite{FMW}, is that, if $V$ is a regular semistable
vector bundle of rank 
$n-1$ and trivial determinant over $C$, then $V\oplus \scrO_C$ is not always
regular. We generalize the study of this problem to a wide variety of
inclusions $G'\subseteq G$. One application, to be given in the next
section, is to the characterization of $G$-bundles $\xi$ arising from the
parabolic construction and representations $\rho\colon G \to GL_N(\Cee)$ such
that the vector bundle $\xi\times _G\Cee^N$ is unstable.  Very related methods
give a rigidification of the construction for $G$ of type
$E_8$ and $C$ a cuspidal curve, allowing us to construct a $9$-dimensional
weighted projective space together with an action of a one-dimensional
unipotent group.  Inside the  $9$-dimensional
weighted projective space, there is a  $9$-dimensional
affine space corresponding to bundles which become trivial on $\widetilde C$,
and the group action on this affine space has an $8$-dimensional affine slice
which plays the role of the Kostant section.

\subsection{Statement of the main theorem}

For the moment, we consider the case where $X$ is allowed to be an arbitrary
scheme or analytic space.

\begin{theorem}\label{Thm1} Suppose that the algebraic group $P$ is a 
semidirect product $L\ltimes U$, where
$U$ is a normal unipotent subgroup of $P$. Let  $U'$ be a subgroup of $U$ and 
let
$Q$ be a subgroup of $L$. Suppose that $U'$ is invariant under the
conjugation  action of $Q$, and that this action factors through a surjection
$r\colon Q\to \ov Q$ with kernel $N$. Let $\eta$ be an
$L$-bundle over
$X$, and suppose that there is a fixed reduction of the structure group of
$\eta$ to $Q$, i.e.\ an isomorphism
$$\eta_Q\times _QL\cong \eta ,$$ where $\eta_Q$ is a $Q$-bundle. Let 
$\eta_{\ov Q}= r_*\eta_Q$. Then:
\begin{enumerate}
\item[\rm (i)] As sheaves of unipotent groups over $X$, $U'(\eta_Q) \cong
U'(\eta_{\ov Q})$, and thus
$H^1(X;U'(\eta_Q)) \cong H^1(X;U'(\eta_{\ov Q}))$, where this isomorphism is
equivariant with respect to the  action of
$\Aut_Q\eta_Q$; 
\item[\rm (ii)] There is a natural inclusion of sheaves of unipotent groups
$U'(\eta_Q) \subseteq U(\eta)$, and hence there is an 
$\Aut_Q\eta_Q$-equivariant function 
$H^1(X;U'(\eta_Q)) \to H^1(X;U(\eta))$. Moreover, if $\xi'$ is the
$\ov QU'$-bundle corresponding  to a class in $H^1(X;U'(\eta_{\ov Q}))$ and
$\xi$ is the $QU'$-bundle corresponding to its image in $H^1(X;U'(\eta_Q))$,
then $\xi/N \cong \xi'$; 
\item[\rm (iii)] Suppose that $\eta$ satisfies the hypotheses of {\rm
\cite[Theorem A.2.2]{FMII}}  with respect to the filtration
$\{U_i\}$ and that either $\dim X =1$ or  the vector group $(U'\cap U_i)/(U'\cap
U_{i+1})$ is a direct summand of $U_i/U_{i+1}$ for every $i$. Then
$\eta_Q$ and
$\eta_{\ov Q}$ satisfy these hypotheses with respect to $\{U'\cap U_i\}$, so that
the functors corresponding to $H^1(X;U'(\eta_Q))$ and $ H^1(X;U'(\eta_{\ov Q}))$
are isomorphic and representable by affine spaces. Moreover the functions
$H^1(X;U'(\eta_Q)) \cong H^1(X;U'(\eta_{\ov Q}))$ and 
$H^1(X;U'(\eta_Q))
\to H^1(X;U(\eta))$ correspond to morphisms of functors which are
represented by morphisms of affine spaces, equivariant with respect to the action
of $\Aut_Q\eta_Q$.
\end{enumerate}
\end{theorem}
\begin{proof} Part (i) follows from the isomorphism $\eta_Q\times _QU' \cong
\eta_{\ov Q}\times _{\ov Q}U'$ (cf.\ 
\cite[Lemma A.1.4]{FMII}). The inclusion of Part (ii) is induced from the
inclusion
$U'(\eta_Q)  =
\eta_Q\times _QU'\subseteq \eta_Q\times_QU$ and the isomorphism
$$U(\eta_Q) = \eta_Q\times_QU \cong (\eta_Q\times_QL)\times _LU = U(\eta),$$
and the second statement in (ii) is clear.
The first statement in  Part (iii) follows since  it is easy to check that the
representability conditions of
\cite[Theorem A.2.2]{FMII} are satisfied. It remains to see that the given
functions are represented by morphisms of affine spaces. It suffices to exhibit
morphisms of the corresponding functors. To do so, for example for the case of
the function
$H^1(X;U'(\eta_Q))
\to H^1(X;U(\eta))$, suppose that $S$ is an affine scheme and that $(\xi',
\varphi')$ is a pair consisting of a $QU'$-bundle $\xi'$ and an isomorphism
$\varphi'\colon \xi'/U' \to \pi_1^*\eta_Q$. Define $\xi =\xi'\times _{QU'}LU$ and
let $\varphi\colon \xi/U \to \eta$ be the isomorphism given by the composition of
the canonical isomorphisms
$$\xi/U = \xi'\times _{QU'}L \cong (\xi'/U')\times _QL$$  with the isomorphism
$\varphi'$ and the fixed identification of $\eta_Q\times _QLU$ with $\eta$.
The map $(\xi', \varphi')\mapsto (\xi, \varphi)$ is then the required
morphism of functors.
\end{proof} 

We note  the following, whose proof is
standard and left to the reader:

\begin{proposition}\label{geomquot} Under the hypotheses  of
Theorem~\ref{Thm1}, suppose that
$U'$ is a normal subgroup of $U$ and that $H^1(X;
(U/U')(\eta_Q))=0$. Then the  group $H^0(X;
(U/U')(\eta_Q))$ acts on $H^1(X;U'(\eta_Q))$ and the quotient of this action
is $H^1(X;U(\eta_Q))= H^1(X;U(\eta))$. If $\eta_Q$ and $U'$ satisfy the
conditions of {\rm \cite[Theorem A.2.2]{FMII}}, then 
$H^0(X; (U/U')(\eta_Q))$ acts as a unipotent algebraic group on the affine
space representing the functor corresponding to $H^1(X;U'(\eta_Q))$. Finally,
if
$\eta$ also satisfies the conditions of {\rm \cite[Theorem A.2.2]{FMII}},
then the affine space corresponding to $H^1(X;U(\eta))$ is a geometric
quotient for the action of $H^0(X; (U/U')(\eta_Q))$ on 
$H^1(X;U'(\eta_Q))$. \qed
\end{proposition}

 In many  cases of interest to us, the homomorphism
$r\colon Q\to \ov Q$ will be a split surjection, i.e.\ there is a homomorphism
$i\colon \ov Q\to Q$ such that
$r\circ i=\Id$. We identify $\ov Q$ with its image in $Q$. Via the isomorphism
$H^1(X;U'(\eta_Q)) \cong H^1(X;U'(\eta_{\ov Q}))$, elements of
$H^1(X;U'(\eta_Q))$, which define $QU'$-bundles $\xi_{QU'}$, also define $\ov
QU'$-bundles $\xi_{\ov QU'}$. Their images in $H^1(X;U(\eta))$ define
$LU$-bundles $\xi_{QU'}\times_{QU'}LU$. However, the bundle
$\xi_{QU'}\times_{QU'}LU$ is not in general isomorphic to the bundle $\xi_{\ov
QU'}\times _{\ov QU'}LU$. Instead, we have the following criterion for when the
bundles
$\xi_{QU'}\times_{QU'}LU$ and
$\xi_{\ov QU'}\times _{\ov QU'}LU$ are S-equivalent in  a generalized sense.

\begin{proposition}\label{Thm2}  Suppose that there is a homomorphism $i\colon
\ov Q\to Q$ such that
$r\circ i=\Id$, and that there exists a subgroup of  $Q$ isomorphic to $\Cee^*$
with the following property: $\Cee^*$ acts trivially on
$U'$ and the conjugation action $\psi \colon Q\times
\Cee^*
\to Q$ extends to a morphism $\ov \psi \colon Q\times \Cee \to Q$ such that
$\ov \psi (q,0) = i\circ r(q)$. Let $\xi_{QU'}$ be a $QU'$-bundle corresponding to
an element of $H^1(X;U'(\eta_Q))$ and let $\xi_{\ov QU'} = \xi_{QU'}\times
_{QU'}\ov QU'$. Then there is an
$QU'$-bundle
$\Xi \to X\times \Cee$ such that, for $t\neq 0$, $\Xi_t=\Xi|X\times \{t\} \cong
\xi_{QU'}$ and such that $\Xi_0 \cong \xi_{\ov QU'}\times _{\ov QU'}QU'$.
\end{proposition}
\begin{proof} Suppose that $\{\Omega_i\}$ is an open cover of $X$ and that
$\xi_{QU'}$ is defined by the cocycle
$q_{ij}u_{ij}$, where $q_{ij}\colon \Omega_i\cap \Omega_j\to Q$ and  $u_{ij}\colon
\Omega_i\cap \Omega_j
\to U'$ are morphisms. Define $\Xi$ by the cocycles $\ov \psi(q_{ij},t)u_{ij}$.
Since $t\in \Cee^*$ acts trivially on $u_{ij}$, $\ov \psi(q_{ij},t)u_{ij} =
t(q_{ij}u_{ij})t^{-1}$ for $t\neq 0$. It follows that the
$\ov \psi(q_{ij},t)u_{ij}$ are indeed cocycles, that $\Xi_t=\Xi|X\times \{t\}
\cong
\xi_{QU'}$ for $t\neq 0$,  and  that $\Xi_0 \cong \xi_{\ov QU'}\times _{\ov
QU'}QU'$.
\end{proof}

One circumstance where we can verify the hypotheses of Proposition~\ref{Thm2}
is the following:

\begin{proposition}\label{Thm3} With the above notation, let $N=\Ker r$, and
let
$\Lie(Q) =
\frak q$,
$\Lie(U') = \frak u'$, $\Lie(N) = \frak n$.   Suppose that $N$ is unipotent
and that, as in Proposition~\ref{Thm2}, there is a homomorphism 
$i\colon \ov Q\to Q$ such that
$r\circ i=\Id$, with $\Lie(\ov Q) =\ov{\frak q}$, which we identify via $i$
with a subalgebra of $\frak q$,  and that there exists a subgroup of 
$Q$ isomorphic to  $\Cee^*$ which acts trivially on $\ov{\frak q}$ and on
$\frak u'$ and with positive weights on $\frak n$. Then:
\begin{enumerate}
\item[\rm (i)] The conjugation action $\psi \colon Q\times
\Cee^*
\to Q$ extends to a morphism $\ov \psi \colon Q\times \Cee \to Q$ such that
$\ov \psi (q,0) = i\circ r(q)$.
\item[\rm (ii)] If $\rho\colon QU'\to GL(V)$ is a finite-dimensional
representation, then there is a filtration on the vector bundle
$\xi_{QU'}\times _{QU'}V$ whose associated graded is isomorphic to
the vector bundle $\xi_{\ov QU'}\times _{\ov QU'}V$.
\item[\rm (iii)] Suppose that $QU'$ is contained in an algebraic group $G$,
that $G'$ is an algebraic subgroup of $G$ containing $\ov QU'$, and that the
inclusion $i\colon \ov Q\to Q$ is induced from the inclusion of $G'$ in $G$.
Suppose further that the subgroup $\Cee^*$ of $Q$ centralizes $G'$, and that 
$\rho\colon QU'\to GL(V)$ is a finite-dimensional
representation which is the restriction of a representation of $G$. Then the
filtration defined in {\rm (ii)} has the property that its associated graded
is isomorphic to the vector bundle $(\xi_{\ov QU'}\times_{\ov QU'}G')\times
_{G'}V$.
\end{enumerate}
\end{proposition}
\begin{proof} To see (i), note that as varieties with a $\Cee^*$-action, $Q
\cong \ov Q \times N\cong \ov Q \times \frak n$ via the exponential map, where
the
$\Cee^*$-action is trivial on the first factor. Since 
$\Cee^*$   acts  with positive weights on $\frak n$, the conjugation
function
$\psi\colon Q\times
\Cee^* \to Q$ defined by $\psi(q,t) = tqt^{-1}$ extends to a morphism  $\ov
\psi
\colon Q\times \Cee \to Q$ such that $\ov \psi (q,0) = i\circ r(q)$. To see
(ii), let $V_i$ be the subspace of $V$ where $\Cee^*$ acts via $\rho$ with
weight $i$. Thus, 
there is a direct sum decomposition $V =\bigoplus _iV_i$, and this direct
sum decomposition is preserved by $\ov QU'$ since  $\Cee^*$ acts trivially on
$\ov{\frak q}$ and on $\frak u'$. On the other hand, since $\Cee^*$ acts with
positive weights on $\frak n$, $NV_i \subseteq \bigoplus _{j>i}V_j$. The
filtration of $V$ defined by $F^i = \bigoplus _{j\geq i}V_j$ is preserved by
$QU' = \ov QNU'$, and the action on the associated graded factors through the
action of $\ov QU'$ on $V$. Part (ii) then follows. Finally, under the
hypotheses of (iii),  the $V_i$ are
$G'$-invariant subspaces, and thus the associated gradeds come from
representations of $G'$.
\end{proof}

\subsection{Parabolic induction}

In this section,  we shall use the above
results to relate
$G'$-bundles to $G$-bundles, where $G'$ is a semisimple subgroup of $G$
whose Dynkin diagram is obtained by deleting one end of the Dynkin diagram
of $G$. A key ingredient in the construction is  the fact that the bundle
$\eta_0$ reduces to a Borel subgroup of $L$. 

We begin by introducing some notation. If $\gamma \in \Delta$ and $\beta \in R$,
let $c_\gamma(\beta)$ be the coefficient of $\gamma$ in the expression of $\beta$
as a sum of simple roots, i.e.\ $c_\gamma(\beta) =\beta(\varpi_\gamma\spcheck)$.
As usual, let $\alpha$ be the special root. Now choose a simple root 
$\alpha_1\neq \alpha$ such that
$\langle \alpha,
\alpha_1\rangle
\neq 0$, i.e.\ the vertex of the Dynkin diagram corresponding to $\alpha_1$ is
adjacent to the vertex corresponding to
$\alpha$. Let $n(\alpha, \alpha_1) = 2\langle \alpha, \alpha_1\rangle/\langle
\alpha_1, \alpha_1\rangle$ be the Cartan integer, and let $m=-n(\alpha,
\alpha_1) >0$. Thus, if we normalize the inner product so that $\langle
\alpha, \alpha\rangle =2$, then  $\langle
\alpha, \alpha_1\rangle =-1$,  $\langle
\alpha_1, \alpha_1\rangle =2/m$, and $\nu= \alpha + m\alpha_1$ is a long root.
Suppose that  $\alpha_1, \dots, \alpha_k$ are the simple
roots such that
$\{\alpha_1, \dots, \alpha_k\}$ is the set of vertices of a connected component
of the Dynkin  diagram of $\Delta -\{\alpha\}$ and such that $\langle
\alpha_i,\alpha_{i+1}\rangle \neq 0$.

 Let $R'$ be the set
of all roots $\beta$ such that $mc_\alpha(\beta) = c_{\alpha_1}(\beta)$. Clearly,
$R'$ is a root system in a hyperplane in the real span of $R$, and $\Delta' =
(\Delta -\{\alpha, \alpha_1\})\cup\{\nu=\alpha+m\alpha_1\}$ is a set of simple
roots for $R'$. Let $\frak h'\subseteq \frak h$ be the complex span of the
corresponding coroots. Then 
$$\frak g' =\frak h'\oplus \bigoplus_{\beta\in R'}\frak g^\beta$$ is a semisimple
subalgebra of $\frak g$. Let $G'$ be the corresponding subgroup of $G$, and
denote by $i\colon G' \hookrightarrow G$ the inclusion. It is easy to see that
$G'$ is again simple and simply connected, and that the Dynkin diagram for
$G'$ is obtained from that of $G$ by replacing $\alpha, \alpha_1$ and the edge
connecting them by a single long root $\nu=\alpha+m\alpha_1$, with
$\nu\spcheck = \alpha \spcheck + \alpha_1\spcheck$. Note that
$n(\gamma, \nu) = n(\gamma, \alpha)$ for $\gamma\in \Delta -\{\alpha,
\alpha_1,
\alpha_2\}$ and $n(\alpha_2, \nu) = -m$. Hence $\nu$ is a special root for
$\Delta'$. The Dynkin diagram for $\Delta'$ is the same as that for
$\Delta -\{\alpha_k\}$. In fact:

\begin{lemma} There exists a $w\in W$ such that $w\Delta'=\Delta
-\{\alpha_k\}\subseteq \Delta$.  Thus, if we let
 $G''$ be the standard simple subgroup of $G$ corresponding to $\Delta
-\{\alpha_k\}$, there is a $g\in G$ normalizing $H$  which conjugates $G'$
onto $G''$.
\end{lemma}
\begin{proof} Let $r_{\alpha_i}\in W$ be the reflection corresponding to the
simple root $\alpha_i$, and let $w = r_{\alpha_k}\circ r_{\alpha_{k-1}}\circ
\cdots \circ r_{\alpha_1}$. Then it is straightforward to check that
$w(\alpha_k) = \alpha_{k-1}, \dots, w(\alpha_2) = \alpha_1, w(\nu) =\alpha$,
and $w(\beta) = \beta$ for $\beta \in \Delta -\{\alpha, \alpha_1, \dots,
\alpha_k\}$. Thus $w$ is as required.
\end{proof}

In the simply laced case, the next lemma says that, for every positive root
$\beta$, the coefficients of $\beta$ expressed as a linear combination of the
simple roots decrease along each arm of the Dynkin diagram.

\begin{lemma}\label{decrease} In the above notation, suppose that $\beta$ is a
positive root and that $c_\alpha(\beta) > 0$. Then $mc_\alpha(\beta) \geq
c_{\alpha_1}(\beta)$. Moreover, for all $i, 1\leq i\leq k-1$, 
$c_{\alpha_i}(\beta)\geq c_{\alpha_{i+1}}(\beta)$.
\end{lemma}
\begin{proof} Let $\Delta _1= \Delta-\{\alpha, \alpha_1, \dots, \alpha_k\}$. Write
$\beta = t\alpha + \sum _it_i\alpha_i + \sum _{\gamma
\in
\Delta_1}s_\gamma\gamma$. The statement of
the lemma is that $mt\geq t_1 \geq t_2 \geq \cdots \geq t_k$. The proof is by
induction on the length
$t+\sum _it_i+
\sum _{\gamma} s_\gamma$ of
$\beta$. If the length is one, then $\beta =\alpha$ and the lemma is clear. Thus
we may assume that $\beta = \beta'+\gamma$ for some $\gamma\in \Delta$ and
$\beta'\in R$. If
$c_\alpha(\beta') =0$, then $\gamma =\alpha$. In this case, since
$\Delta-\{\alpha\}$ is a union of diagrams of type $A$,
$\beta'$ is a sum of simple roots and there exists an $\ell, 1\leq \ell <k$, such
that
$t_i=1$ for $i<\ell$ and $t_i-0$ for $i>\ell$. The lemma is clearly true in this
case. Thus we may assume that
$c_\alpha(\beta')>0$, so that by the inductive hypothesis $\beta' = t'\alpha +
\sum _it_i'\alpha_i + \sum _{\gamma \in
\Delta_1}s_\gamma'\gamma$ with $mt'\geq t_1'\geq t_2' \cdots$.
If
$\gamma
\neq \alpha_i$, then $t_i'=t_i$ and $t\geq t'$ and the inductive step holds. If
$\gamma =\alpha_1$, then $t_1 = t_1'-1$, $t_i=t_i'$ for $i> 1$, and $t=t'$. Since
$\beta'+\alpha_1$ is a root, $n(\beta', \alpha_1) \leq 0$. If $n(\beta',
\alpha_1) < 0$, then
$$-1 \geq \frac{2t'\langle \alpha, \alpha_1\rangle}{\langle \alpha_1,
\alpha_1\rangle} +2t_1'-t_2' = (t_1'-mt') + (t_1'-t_2') = t_1-1-mt +
(t_1'-t_2').$$ Since $t_1'-t_2'\geq 0$, we must have $t_1\leq mt$, which is the
statement of the lemma in this case. If $n(\beta', \alpha_1) = 0$, then
$\beta-\alpha_1$ is also a root, and so by the inductive hypothesis
$(t_1'-1)-t_2'\geq 0$, and a similar argument shows that $t_1\leq mt$. A similar
but simpler argument handles the case where
$\gamma =\alpha_i$ for some $i>1$.
\end{proof}

Let $P=LU$ be the maximal parabolic subgroup associated
to $\alpha$, with Levi factor $L$ and unipotent radical $U$. Recall that
$\frak u =\Lie (U) =
\bigoplus_{c_\alpha(\beta)>0}\frak g^\beta$. Let
$Q$ be the Borel subgroup of $L$ with
$$\frak q =\Lie(Q) = \frak h\oplus \bigoplus _{c_\alpha(\beta) =0, \beta >
0}\frak g^\beta.$$ 

  Define
the subgroup
$U'$ of
$U$ as follows:
$U'$ is the connected unipotent subgroup of $U$ such that
$$\frak u' =\Lie(U') = \bigoplus _{mc_\alpha(\beta) = c_{\alpha_1}(\beta)>
0}\frak g^\beta.$$ 
Let $N$ be the unipotent  subgroup
of $Q$ such that
$$ \frak n = \Lie (N) = \bigoplus_{c_\alpha(\beta) =0,
c_{\alpha_1}(\beta)>0}\frak g^\beta.$$

\begin{lemma} We have $[\frak q, \frak u']\subseteq \frak u'$. Moreover, if
$c_\alpha(\beta) =0$ and
$c_{\alpha_1}(\beta)>0$, then  $[\frak g^\beta, \frak u'] =0$. Thus $[ \frak n,
\frak u'] =0$, and $ \frak n$ is normal in $\frak q$.
\end{lemma}
\begin{proof} The proof is immediate from Lemma~\ref{decrease}.
\end{proof}

\begin{corollary} Conjugation by $Q$ leaves the subgroup $U'$ invariant.
Moreover, the normal subgroup $N$ of $Q$  acts trivially on $U'$. Thus $N$
is a normal subgroup of $QU'$. \qed
\end{corollary}

Let
$P'=(P')^\nu$ be the standard maximal parabolic subgroup of $G'$ corresponding to
the simple root
$\nu$. Then the group
$U'$ is the unipotent radical of $P'$. Let $L^{\{\alpha, \alpha_1\}} = 
L_{\Delta-\{\alpha,\alpha_1\}}$ be the closed connected subgroup of $G$ whose Lie
algebra is
$$\Lie(L^{\{\alpha, \alpha_1\}})= \frak h \oplus \bigoplus_{c_\alpha(\beta) =
c_{\alpha_1}(\beta) =0}\frak g^\beta.$$ There is a Levi factor $L'$ of $P'$ whose
Lie algebra is given by
$$\Lie(L')=\frak h' \oplus \bigoplus_{\beta\in R', c'_\nu(\beta)=0}\frak
g^\beta,$$ where the coefficient $c'_\nu(\beta)$ refers to writing $\beta$ as a
sum of elements of $\Delta'$.  Hence there is a finite group
$F$ such that
$L^{\{\alpha,
\alpha_1\}} \cong L'\times _F\Cee^*$. Here $\Lie(\Cee^*)\subseteq \frak h$ is
the kernel of all simple roots except for $\alpha$ and $\alpha_1$, as well as
the root $\nu$, and a nonzero element of this Lie algebra is given by
$-m\varpi_\alpha\spcheck + \varpi_{\alpha_1}\spcheck$. The group
$\Cee^*$ thus acts trivially on
$U'$ and with positive weights on $ \frak n$, and centralizes $G'$.

Quite concretely, $L$ is the subgroup of
 $\prod_{i=1}^tGL_{n_i}(\Cee)$, where the first factor $GL_{n_1}(\Cee)$
corresponds to the simple roots $\alpha_1,
\dots, \alpha_k\in \Delta$, consisting of all tuples $(A_1, \dots, A_t)$ such
that $\det A_1 = \cdots = \det A_t$. Then
$L^{\{\alpha, \alpha_1\}}$ is the subgroup of  $L$ of all such tuples where
$$A_1 = \begin{pmatrix}
\lambda &0\\ 0& A_1'
\end{pmatrix}.$$ In particular, the corresponding roots for the first factor are
$\alpha_2, \dots,
\alpha_k$. On the other hand, $L'$ is the subgroup of $GL_{n_1-1}(\Cee) \times
\prod_{i=2}^tGL_{n_i}(\Cee)$ consisting of all tuples $(A_1', \dots, A_t)$ such
that $\det A_1' = \cdots = \det A_t$. Given an upper triangular $n_1\times n_1$
matrix, we can map it to the lower right hand block of size
$(n_1-1)\times (n_1-1)$, and this is a homomorphism onto the group of upper
triangular matrices in $GL_{n_1-1}(\Cee)$. 

Returning to the general theory, we have
$$\frak q/ \frak n \cong \frak h \oplus \bigoplus_{c_\alpha(\beta) =
c_{\alpha_1}(\beta) = 0, \beta >0}\frak g^\beta.$$
Thus

\begin{lemma} There is a surjection $r$ from $Q$ to a Borel subgroup $\ov Q$ of
$L^{\{\alpha,\alpha_1\}} \cong L'\times _F\Cee^*$, which identifies $\ov Q$
with $Q/N$. Moreover, if $i\colon \ov Q\to Q$  is induced by the inclusion
$i\colon G' \to G$, then 
$r\circ i=\Id$. Conjugation by the above $\Cee^*$ in the center of
$L^{\{\alpha,\alpha_1\}}$ acts trivially on $G'$, normalizes $N$, and acts
with positive weights on $\frak n$. \qed
\end{lemma}

An $L$-bundle $\eta$ consists of $t$ vector bundles $V_i$ of rank $n_i$, such 
that
$\det V_1 = \cdots = \det V_t$. An $L^{\{\alpha,
\alpha_1\}}$-bundle consists of a line bundle $\lambda$, a vector bundle 
$V_1'$ of rank $n_1-1$, and bundles $V_2, \dots, V_t$ of rank $n_i$ such that
$\lambda\otimes
\det V_1' = \det V_2  = \cdots = \det V_t$. An $L'$-bundle consists of a 
vector  bundle $V_1'$ of rank $n_1-1$, and bundles $V_2, \dots, V_t$ of rank
$n_i$ such that $\det V_1' = \det V_2  = \cdots = \det V_t$.

Let $\eta_0'$ be the unique stable
$L'$-bundle of determinant $\scrO_C(-p_0)$. Because of the inductive 
construction of $W_k\spcheck$ via the exact sequence
$$0\to W_k\spcheck \to W_{k+1}\spcheck \to \scrO_C \to 0,$$ we see:

\begin{lemma}\label{eta0reduces} There is a holomorphic reduction of structure
$\eta_Q$ of
$\eta_0$ to a $Q$-bundle so that $r_*\eta_Q$ lifts to $\eta_0'$. \qed
\end{lemma}

By Theorem~\ref{Thm1}, we have:

\begin{corollary} Let  notation be as above. Then there is a 
$\Cee^*$-equivariant function  from $H^1(C; U'(\eta_Q))$ to
$H^1(C; U(\eta_0))$, and a $\Cee^*$-equivariant isomorphism from $H^1(C;
U'(\eta_Q))$ to $H^1(C; U'(\eta_0'))$. Except in the case where
$C$ is cuspidal and
$G=E_8$,  the cohomology sets are represented by affine spaces and the functions
between them are  represented by $\Cee^*$-equivariant morphisms. \qed
\end{corollary}

In the statement of the above corollary, one has to be a little careful about
what is meant by the $\Cee^*$-actions, because the centers of $L$ and $L'$ do
not correspond to the same subgroup of $Q$. Instead, there is a homomorphism
from
$Z_1\times Z_2$ into the center of $L^{\{\alpha, \alpha_1\}}$, where
$Z_1\cong \Cee^*$ is the identity component of the center of $L$ and the
$Z_2\cong \Cee^*$ is the identity component of the center of $L'$. The first
$\Cee^*$ factor acts on $H^1(C; (\frak u')^k(\eta_Q))$  and $H^1(C; \frak
u^k(\eta_0))$ via $kn_\alpha$ and the second acts on  $H^1(C; (\frak 
u')^k(\eta_Q))$  and 
$H^1(C; (\frak u')^k(\eta_0'))$ via $kn_\alpha'$, where $n_\alpha'=n_\nu$,
and these weights cannot agree unless $n_\alpha = n_\nu$. The point is that
there is a third $\Cee^*$ in $Z(L^{\{\alpha, \alpha_1\}})$ such that
$L^{\{\alpha, \alpha_1\}} \cong L'\times _F\Cee^*$, and this $\Cee^*$ acts
trivially on $U'$. Thus, there is a positive integer $n$ such that, for
all $t \in Z_1$, $t^n =s \mu$ where $s\in Z_2$ and $\mu $ acts trivially on
$U'$, and in this way we can line up the two different $\Cee^*$-actions up to
a positive multiple.

Applying  Proposition~\ref{Thm2}, and Proposition~\ref{Thm3} in this
situation gives:

\begin{lemma}\label{limit} Let $x\in  H^1(C; U'(\eta_Q))$, and let $\xi$ be
the
$G$-bundle corresponding to the image of $x$ in $H^1(C; U(\eta_0))$. Let
$\xi'$ be the $G'$-bundle corresponding to the image of $x$ in $H^1(C;
U'(\eta_0'))$. Then there is a $G$-bundle $\Xi \to C\times \Cee$ such that,
for $t\neq 0$, 
$\Xi_t=\Xi|C\times \{t\} \cong
\xi$ and such that $\Xi_0 \cong \xi'\times _{G'}G$. Moreover, if
$\rho\colon G \to \Cee^N$ is a representation of $G$, then there is a
filtration on the vector bundle $\xi\times _G\Cee^N$ whose associated
graded is isomorphic to $\xi'\times _{G'}\Cee^N$, where $G'$ acts on
$\Cee^N$ via $\rho|G'$. \qed
\end{lemma}

If $C$ is smooth, the above lemma says that the semistable and regular 
$G$-bundle
$\xi$ is S-equivalent to the semistable, but not necessarily regular,
$G$-bundle $\xi'\times _{G'}G$. In other words, the morphism of affine spaces
$H^1(C; U'(\eta_0')) \to H^1(C; U(\eta_0))$ defined by the previous corollary
induces a corresponding map on weighted projective spaces $\WP(G') \to 
\WP(G)$, which via 
Looijenga's theorem is the same as the natural morphism $\mathcal{M}(G') \to 
\mathcal{M}(G)$ of moduli spaces. Of course, this morphism is also induced
from the obvious morphism
$(C\otimes \Lambda')/W' \to (C\otimes \Lambda)/W$, and in particular it is
finite.  However, this finite morphism need not be injective. A similar
picture holds in the singular case on the open set $\WP(G) -D_\infty
\cong (C_{\rm reg}\otimes \Lambda)/W$.

   We now show that, even when $C$ is singular, for most pairs of groups, the
function
$H^1(C; U'(\eta_0'))
\to H^1(C; U(\eta_0))$ induces a morphism on the level of weighted projective
spaces, if this is meaningful:

\begin{proposition}\label{isnotzero} Suppose that $G$ is not of type of type
$E_8$, and  that the pair $(G',G)$ is not of type  
$(D_6,E_7)$. If
$x\in  H^1(C; U'(\eta_Q)) -\{0\}$, then the image of
$x$ in $H^1(C; U(\eta_0))$ is nonzero.
\end{proposition}
\begin{proof} We will deal with the case where $(G',G)$ is of type  
$(E_6,E_7)$ in
\cite{FMIV} via the link with del Pezzo surfaces; the argument will also
handle the case where $(G',G)$ is of type $(D_5,E_6)$ or $(A_4,D_5)$.  Here
we will consider the case where
$(G',G) = (Spin(2n-2), Spin(2n))$ and $n\geq 4$. The other classical cases or
$(G',G)$ of type $(A_1,G_2)$ are similar and for the most part
easier. The case where $G$ is of type $E_6$ or $F_4$ involves a detailed
analysis of the cubic form on the corresponding $27$ or $26$-dimensional
representation but is otherwise similar, and the case of $(A_6,E_7)$ can also
be handled by these methods.

Consider the standard representation of $Spin(2n)$. The vector bundle
$V$ associated to a $Spin(2n)$-bundle $\xi$ arising from the parabolic
construction has a quadratic form and a filtration $F_{-1}\subseteq F_0 \subseteq
F_1=V$ such that $F_{-1}\cong W_{n-2}\spcheck$, $F_{-1}^\perp = F_0$, and
$F_0/F_{-1}\cong Hom(W_2, W_2)$ with its standard quadratic form. The bundle
$\xi\cong
\eta_0\times_LG$ if and only if this filtration is split compatibly with the
quadratic form, i.e.\ $V\cong W_{n-2}\spcheck \oplus Hom(W_2, W_2)\oplus
W_{n-2}$, where $W_{n-2}\spcheck $ and $W_{n-2}$ are isotropic and both are
orthogonal to $Hom(W_2, W_2)$.

If $\xi$ arises via parabolic induction from a $Spin(2n-2)$-bundle $\xi'$, then
there is a second filtration $G_{-1}\subseteq G_0 \subseteq
G_1=V$ with $G_{-1}\cong \scrO_C$, $G_{-1}^\perp = G_0$, and
$G_0/G_{-1}=V'$, where $V'$ is an $SO(2n-2)$-bundle associated to $\xi'$. 
Furthermore, $G_{-1}\subseteq V$ projects to $F_1/F_0$ to give the
natural inclusion of $\scrO_C\to W_{n-2}$. Thus, the filtration
$$F_{-1}\cap G_0\subseteq F_0\cap G_0\subseteq G_0/G_{-1}$$
is the parabolic filtration of $\xi'$.
Now suppose that $\xi$ is given by $0\in H^1(C;U(\eta_0))$ so that
$V=W_{n-2}\spcheck\oplus Hom(W_2,W_2)\oplus W_{n-2}$ compatibly with the
quadratic form. The image of $G_{-1}$ in the
$W_{n-2}\spcheck$ factor is trivial by stability. Hence
$G_{-1}$ is contained in the sum of the last two factors and   the
restriction of the quadratic form of $V$ to $G_{-1}$ is the same as the
restriction of the quadratic form to the projection of $G_{-1}$ into
$Hom(W_2,W_2)$. But there is up to scale a unique map of
$\scrO_C$ to $Hom(W_2,W_2)$, and the quadratic form on the image
of any nonzero map is nontrivial. Since   the restriction
of the quadratic form to $G_{-1}$ is trivial, it follows that the
projection of $G_{-1}$ to $Hom(W_2,W_2)$ is zero, and hence  
$G_{-1}$ is contained in the last factor $W_{n-2}$ of the direct sum
decomposition of $V$. It then follows that the filtration of $V'$
induced as above is also split and hence that the parabolic bundle $\xi'$ is
given by $0\in H^1(C;U'(\eta_Q))$.
\end{proof}

\begin{remark} (1) It is reasonable to expect that more complicated arguments
will extend  Proposition~\ref{isnotzero} to the remaining  cases of parabolic
induction, and that there is in fact a unified proof. The cases listed above
will cover all of the applications in the next section.

(2) Similar arguments work in the case of a family $Z\to B$ or
the linearized case  of $\mathcal{E} \to \mathbb{A}^2$, except that we must
work in various conformal forms  of the groups $G$ and $G'$, using the bundle
$\hat \eta_0$, and a somewhat involved argument shows that there is a morphism
between the corresponding weighted projective bundles. However, we shall not
give the details here. 

(3) If we iterate the above construction, we obtain an embedding of
$SL_2(\Cee)$ in $G$ corresponding to the coroot $\lambda_1(\alpha)\spcheck$.
It is easy to see that the corresponding morphism of $\Pee^1$ to the
weighted projective space corresponding to $G$ is an embedding onto a
subprojective space in weight one. From this, one can show the following,
assuming that
$G$ is not of type $E_8$: For every identification of the weighted
projective bundle over the moduli stack with the weighted projective bundle
associated to
$\scrO_B\oplus \mathcal{L}^{-d_1} \oplus \cdots \oplus \mathcal{L}^{-d_r}$,
the section
$\Pee(\scrO_B) \subseteq \WP(\scrO_B\oplus \mathcal{L}^{-d_1} \oplus \cdots
\oplus \mathcal{L}^{-d_r})$, which is independent of the choice of the
identification, is the section corresponding to the S-equivalence class of the
trivial bundle, i.e.\ the image of the origin in
$(C_{\textrm{reg}}\otimes \Lambda)/W$ for every fiber $C$. Of course, over the
cusp fibers, this is consistent with the discussion of \S\ref{4.6}.
\end{remark}

\begin{defn} If $G'\subseteq G$ is an inclusion of subgroups as in this
section and the conclusion of Proposition~\ref{isnotzero} holds,   let
$D_{G'} \subseteq \WP(G)$ be the image of $\WP(G')$ via the above
construction. Clearly, $D_{G'}$ is an irreducible hypersurface in $\WP(G)$.
We say that the
$G$-bundles corresponding to points of $D_{G'}$ are obtained by
\textsl{parabolic induction}.
\end{defn}

\subsection{The case of $E_8$}

We turn now to $E_8$. We will again use Lemma~\ref{eta0reduces}, that there is
a holomorphic reduction of structure
$\eta_Q$ of
$\eta_0$ to a $Q$-bundle, where $Q$ is an appropriate Borel subgroup of $L$.
In this case, $U'$ will be a normal subgroup of $U$. In terms of Lie
algebras, we have $\frak u = \bigoplus _{k\geq 1}\frak u^k$. Moreover,
$\bigoplus _{k> 1}\frak u^k$ is an ideal of $\frak u$, and the quotient is
abelian. Thus, any vector subspace $\frak u'$ of $\frak u$ containing
$\bigoplus _{k> 1}\frak u^k$ is in fact  an ideal of $\frak u$, and the
corresponding subgroup $U'$ of $U$ is a normal unipotent subgroup with
abelian quotient. Of course, we also need $Q$ to act on $U'$. In our case,
$\frak u^1$ is an
$L$-module. Viewing $L$ as a $\Cee^*$-extension of $SL_2\times SL_3\times
SL_5$, we can identify $\frak u^1$ with $V_2\otimes V_3^*\otimes V_5^*$,
where $V_2$ is the standard representation of $SL_2$ and similarly for $V_3$
and $V_5$. The Borel subgroup $Q$ fixes a filtration of $V_5$,  and hence it
fixes a four-dimensional quotient $V_4$ of $V_5$. Thus there is an induced
subspace
$V_4^*\subseteq V_5^*$. We then have the subspace $V_2\otimes V_3^*\otimes
V_4^*\subseteq V_2\otimes V_3^*\otimes V_5^*$.  Let
$\frak u' = \bigoplus _{k> 1}\frak u^k
\oplus  (V_2\otimes V_3^*\otimes V_4^*)$. Thus $\frak u'$ is an ideal of
$\frak u$ and the quotient $\frak u/\frak u'$ is abelian and isomorphic to
$V_2\otimes V_3^*$. If $U'$ is the corresponding subgroup of $U$, then $U'$
is normal in $U$ and $U/U'\cong V_2\otimes V_3^*$. 

Let $\eta_Q$ be the reduction of $\eta_0$ to $Q$. There is thus an exact
sequence
$$1\to U'(\eta_Q) \to U(\eta_0) \to W_2\otimes W_3\spcheck \to 0.$$
Let $\frak u_i = \bigoplus _{k\geq i}\frak u^i$ and let $U_i$ be the
corresponding subgroup of $U$.  The filtration $\{U_i\}$ induces a decreasing
filtration $\{U_i'\}$ of $U'$. Note that $U_i = U_i'$ for
$i\geq 2$.

\begin{proposition} In the above notation, let $\pi\colon Z\to B$ be a
Weierstrass fibration, let
$\eta_0$ be the
$L$-bundle  over $Z$  corresponding to the triple 
$((\mathcal{W}_2\otimes\pi^*\mathcal{L}^7)\spcheck,
(\mathcal{W}_3\otimes\pi^*\mathcal{L}^4)\spcheck, 
(\mathcal{W}_5\otimes\pi^*\mathcal{L})\spcheck)$ and let $\eta_Q$ be the
natural reduction of $\eta_0$ to a $Q$-bundle.
\begin{enumerate}
\item[\rm (i)] The group $U'$, the bundle $\eta_Q$, and the filtration
$\{U_i'\}$ satisfy the hypotheses of {\rm \cite[Theorem A.2.2]{FMII}}.
\item[\rm (ii)] The functor  corresponding to $H^1(Z; U'(\eta_Q))$ is
representable by a bundle of affine spaces $\mathcal{A}$ over $B$ with a
$\Cee^*$-action induced by the center of $L$ which is
$\Cee^*$-equivariantly isomorphic to
$$(\mathcal{L}^4\oplus \mathcal{L}^2) \oplus (\mathcal{L}^2\oplus
\scrO_B)\oplus (\scrO_B\oplus \mathcal{L}^{-2})\oplus (\mathcal{L}^{-2}\oplus
\mathcal{L}^{-4}) \oplus \mathcal{L}^{-4} \oplus \mathcal{L}^{-6},$$
where the summands have $\Cee^*$-weight $1,2,3,4,5,6$ respectively.
\item[\rm (iii)] The unipotent group scheme $R^0\pi_*(U/U')(\eta_Q)$ over
$B$  is isomorphic to the total space of the line bundle
$\mathcal{L}^{-2}$. The $\Cee^*$-action induced by the center of $L$ is the
standard linear action on $\mathcal{L}^{-2}$ of weight one.
\item[\rm (iv)] The   action of the unipotent group
scheme of {\rm (iii)} on 
$\mathcal{A}$ given in Proposition~\ref{geomquot} is $\Cee^*$-equivariant. The
differential of this action along the zero-section is given by
$(c_1G_3,c_2G_2, 0,0,0,0,0,0,0)$, where the
$c_i$ are nonzero. Away from the cuspidal fibers, the bundle
of  $9$-dimensional affine spaces constructed in Corollary~\ref{affinebundle}
is a geometric quotient of this action.
\end{enumerate}
\end{proposition}
\begin{proof} The vector bundles $(U_i'/U_{i+1}')(\eta_Q)$ are isomorphic to
$\frak u^i(\eta_0)$ for $i> 1$ and to $W_2\otimes W_3\spcheck \otimes
W_4\spcheck$ for $i=1$. Thus, the vanishing of $H^0(C;
(U_i'/U_{i+1}')(\eta_Q))$ follows from the computations used to prove
Theorem~\ref{casebycase}. This proves Part (i). The representability
statement in Part (ii) follows from \cite[Theorem A.2.6]{FMII}, and the
arguments of \cite[Lemma 4.2.3]{FMII} show that this affine space is
$\Cee^*$-equivariantly isomorphic to $R^1\pi_*\frak u'(\eta_Q)$. The rest of
Part (ii) follows by direct computation.

To see Part (iii), we have 
$$U/U'(\eta_Q) \cong \mathcal{W}_2 \otimes \pi^*\mathcal{L}^7\otimes
(\mathcal{W}_3 \otimes \pi^*\mathcal{L}^4)\spcheck \otimes
\pi^*\mathcal{L}^{-5} \cong \mathcal{W}_2\otimes
\mathcal{W}_3\spcheck \otimes \pi^*\mathcal{L}^{-2}.$$
The calculation of
$R^0\pi_*(\mathcal{W}_2\otimes
\mathcal{W}_3\spcheck)$ is given in Lemma~\ref{somemorecomps} and the
identification of the differential of the action follows from the calculation
of the coboundary map in Lemma~\ref{E8}.  Given that 
$R^1\pi_*(\mathcal{W}_2\otimes
\mathcal{W}_3\spcheck)=0$ by Lemma~\ref{somemorecomps}, the final
statement in (iv) is a consequence of 
Proposition~\ref{geomquot}.
\end{proof}

Given a fibration  $\pi\colon Z \to
B$, it is natural to take the weighted tensor product with  
$\mathcal{L}^{-4}$ so as to get the family of weighted projective bundles
associated to
$$(\scrO_B \oplus \mathcal{L}^{-2})\oplus (\mathcal{L}^{-6}
\oplus \mathcal{L}^{-8})
\oplus (\mathcal{L}^{-12}\oplus \mathcal{L}^{-14})\oplus
(\mathcal{L}^{-18}\oplus \mathcal{L}^{-20})\oplus
\mathcal{L}^{-24}\oplus \mathcal{L}^{-30},$$
where the summands have weights $1,2,3,4,5,6$ respectively. The
unipotent group scheme then becomes the total space of the line bundle
$\mathcal{L}^{-2}\otimes \mathcal{L}^{-4}\cong \mathcal{L}^{-6}$.  Except for
the term
$\mathcal{L}^{-6}$ in the above direct sum, we see the Casimir weights of
$E_8$. We now give an explanation for this, and will show that the action of
the unipotent group scheme accounts for the extra term
$\mathcal{L}^{-6}$.

For $C$  cuspidal, let $\Aff^{10}$ be the $10$-dimensional affine space
and let $\WP^9$ be the 9-dimensional weighted projective space  constructed
above from
$H^1(C; U'(\eta_Q))$. As in Theorem~\ref{451} and the discussion following it,
there is a
$\Cee^*$-equivariant morphism from $\Aff^{10}$ to $\Aff^1$ and compatible 
actions on the universal $E_8$-bundles $\Xi\to C\times \Aff^{10}$ and
$(\iota\times \Id)^* \Xi \to  \widetilde C\times \Aff^1 $,
covering the  product of this action with the 
$\lambda$-action on $C$ given in Definition~\ref{action}. The
method of Proposition~\ref{reponnorm} constructs a Zariski open  subset
$\WP^9_0$ isomorphic to $\Aff^9$, consisting of bundles whose pullbacks to
the normalization $\widetilde C$ are trivial. Let us denote  the
one-dimensional unipotent group scheme  by 
$\mathbb{G}_a$. The map from $H^1(C; U'(\eta_Q))$ to $H^1(C; U(\eta_Q))$
is constant on $\mathbb{G}_a$-orbits. Hence, given
$x\in
\WP^9$, if
$\xi_x$ is the associated $G$-bundle, then the action of$\mathbb{G}_a$ does
not change the isomorphism class of $\xi_x$. Thus
$\mathbb{G}_a$ acts on
$\WP^9_0$, and the induced morphism from $\WP^9_0$ to $\frak h/W$ is 
$\mathbb{G}_a$-invariant. To obtain the parabolic construction of an analogue
of the Kostant section for
$G$ of type
$E_8$, we use the following:

\begin{theorem}\label{longcomp}  Fix a basis
$e_0,
\dots, e_{9}$ of
$\Aff^{10}$ corresponding to the ordering of the line bundles above. Then the
subset
$\mathbb{S}$ of
$\Aff^{10}$ where the coefficient of
$e_0$ is $1$ and the coefficient of $e_2$ is $0$ is a slice for the action of
$\mathbb{G}_a$ on 
$\WP^9_0$.
\end{theorem}
\begin{proof} The proof is for the most part a very involved calculation, and
we shall just sketch the general outlines. Note that $\mathbb{G}_a = H^0(C;
(U/U')(\eta_Q))$ and the natural morphism $H^1(C; U'(\eta_Q)) \to H^1(C;
(U'/U_k)(\eta_Q))$ is equivariant with respect to the actions of
$\mathbb{G}_a$ and both actions of $\Cee^*$. The action of $\mathbb{G}_a$ on
$H^1(C; (U'/U_2)(\eta_Q))$ is given by the coboundary and is trivial, by
Lemma~\ref{E8}. The
fibers  of the morphism $H^1(C; (U'/U_3)(\eta_Q))\to  H^1(C;
(U'/U_2)(\eta_Q))$ are principal homogeneous spaces over 
$H^1(C;(U_2/U_3)(\eta_Q))$, which is a vector space.
Given   $\zeta \in H^1(C; (U'/U_2)(\eta_Q))$ and a
lift of $\zeta$ to an element $\widetilde \zeta \in H^1(C;
(U'/U_3)(\eta_Q))$, and given $\gamma \in \mathbb{G}_a$, $\gamma \cdot
\widetilde \zeta$ lies in the same fiber as $\widetilde \zeta$. Thus the
difference $\gamma \cdot
\widetilde \zeta - \widetilde \zeta$ defines an element 
$p(\gamma, \widetilde \zeta) \in  H^1(C;(U_2/U_3)(\eta_Q))$. Now 
$\mathbb{G}_a \cong H^0(C; (U/U')(\eta_Q))$ and the  homomorphism
$H^0(C; (U/U_2)(\eta_Q)) \to H^0(C; (U/U')(\eta_Q))$ is an isomorphism.
Furthermore, $H^0(C; (U/U_2)(\eta_Q)) \cong H^0(C; W_2\otimes W_3\spcheck
\otimes W_5\spcheck)$. Thus we can identify $\gamma$ with an element of 
$H^0(C; W_2\otimes W_3\spcheck
\otimes W_5\spcheck)$, which we shall continue to denote by $\gamma$.
Moreover, $\zeta \in H^1(C; (U'/U_2)(\eta_Q))\cong H^1(C; W_2\otimes
W_3\spcheck \otimes W_4\spcheck)$. Of course, there is the surjection
$j\colon H^1(C; W_2\otimes W_3\spcheck \otimes W_4\spcheck) \to H^1(C;
W_2\otimes W_3\spcheck \otimes W_5\spcheck)$. Finally,
$$H^1(C;
(U_2/U_3)(\eta_Q)) \cong H^1(C; \bigwedge^2W_2\otimes \bigwedge^2W_3\spcheck
\otimes \bigwedge^2W_5\spcheck).$$

With this said, a cocycle computation gives the following:

\begin{lemma} Under the above identifications, $p(\gamma, \widetilde \zeta) =
\gamma \smile j(\zeta)$, where the cup product 
$$H^0(C; W_2\otimes W_3\spcheck \otimes W_5\spcheck) \otimes H^1(C;
W_2\otimes W_3\spcheck \otimes W_5\spcheck) \to H^1(C;\bigwedge^2W_2\otimes 
\bigwedge^2W_3\spcheck
\otimes \bigwedge^2W_5\spcheck)$$
 is the product corresponding to the homomorphism
$$\left(W_2\otimes W_3\spcheck \otimes W_5\spcheck\right) \otimes
\left(W_2\otimes W_3\spcheck \otimes W_5\spcheck\right) \to 
\bigwedge^2W_2\otimes  \bigwedge^2W_3\spcheck \otimes
\bigwedge^2W_5\spcheck$$
which is induced from the natural homomorphisms $W_i\otimes W_i \to \bigwedge
^2W_i$. In particular, $p(\gamma, \widetilde \zeta)$ only depends on the
image $\zeta$ of $\widetilde \zeta$ in $H^0(C; (U'/U_2)(\eta_Q))$ and it
is linear in $\gamma$. \qed 
\end{lemma}

Next, detailed computations as in the proof of Lemma~\ref{E8comps}, using the
description of $W_n$ given in Lemma~\ref{cuspcomps}, the recipe for computing
$H^0(C;V)$ given in Corollary~\ref{linalg} and an analogous formula for
$H^1(C;V)$, where $V$ is a vector bundle on $C$ described by a bundle on
$\widetilde C$ together with an endomorphism of the fiber over $0$, give:

\begin{lemma}\label{coboun} Suppose that $\tau \in H^1(C;
W_2\otimes W_3\spcheck \otimes W_5\spcheck)$ maps to a nonzero element of 
$H^1(\widetilde C; \iota^*(W_2\otimes W_3\spcheck \otimes W_5\spcheck))$ and
that $\gamma \in H^0(C; W_2\otimes W_3\spcheck \otimes W_5\spcheck)$ is not
zero. Then  $\gamma \smile  \tau\neq 0$. \qed
\end{lemma}

To complete the proof of Theorem~\ref{longcomp}, choose a linearization
$$H^1(C; (U'/U)(\eta_Q)) \cong \mathcal{L}_0^{0} \oplus \mathcal{L}_0^{-2}
\oplus \mathcal{L}_0^{-6} \oplus \mathcal{L}_0^{-8}.$$
(Here, the actual weights are not important.) The above lemmas imply that the
$\mathbb{G}_a$-action is of the form
$$\gamma\cdot (x_0, x_1, x_2, x_3) = (x_0, x_1, x_2 + c_2x_0\gamma,
x_3+c_3x_1\gamma)$$ for some constants $c_2, c_3$ such that $c_2\neq 0$.
Fixing
$x_0 =1$, then, the subset where $x_2=0$ gives a slice to the
$\mathbb{G}_a$-action, as required.
\end{proof}

Let $\mathbb{S}$ be the $8$-dimensional affine slice of
Theorem~\ref{longcomp}. After taking a suitable weighted tensor product,
there is a $\Cee^*$-action on $\mathbb{S}$ covering the $\Cee^*$-action
on $C$. Since $\mathbb{S} \subseteq
\Aff^{10}$, there is a universal $E_8$-bundle over $C\times \mathbb{S}$ which
is trivial over
$\widetilde C \times \mathbb{S}$. By Theorem~\ref{cuspfamilies}, choosing a
trivialization gives  an induced morphism
$\sigma
\colon \mathbb{S} \to \frak g$, where $\frak g$ is the Lie algebra of $E_8$.
Let $\ov\sigma\colon \mathbb{S}\to \frak h/W$ be the induced morphism.

\begin{proposition} In the above notation, the image of $\sigma$ is contained
in $\frak g^{\mathrm{reg}}$. The differential of $\ov \sigma$ is
everywhere an isomorphism. Finally, $\ov \sigma$  is 
$\Cee^*$-equivariant with respect to the action on $\mathbb{S}$ defined above
and the inverse of the usual action on $\frak h/W$.
\end{proposition}
\begin{proof} Given $x\in \mathbb{S}$, let $\xi_{x,P}$ be the corresponding
$P$-bundle. First, we claim that
$H^0(C;
\frak u(\xi_{x,P})) =0$. Indeed, there is an exact sequence
$$0 \to (\frak u_2/\frak u_3)(\eta_Q) \to (\frak u/\frak u_3)(\xi_{x,P}) \to
(\frak u/\frak u_2)(\eta_Q)\to 0,$$
and a straightforward calculation shows the coboundary 
$$H^0(C; (\frak u/\frak
u_2)(\eta_Q))
\to H^1(C; (\frak u_2/\frak u_3)(\eta_Q))$$ is given by taking cup product
with the class of the image of $x$ in $H^1(C; (\frak u/\frak u_2)(\eta_Q))$.
By Lemma~\ref{coboun}, the coboundary map is  an isomorphism for $x\in
\mathbb{S}$. Thus $H^0(C; (\frak u/\frak u_3)(\xi_{x,P}))=0$, and it follows
easily that $H^0(C; \frak u(\xi_{x,P})) =0$. Given this, the arguments of
\cite[4.5.1 and 4.5.2]{FMII} imply that the homomorphism $H^1(C; \frak
u(\xi_{x,P})) \to H^1(C; \ad_G(\xi_{x,P}\times _PG))$ is surjective. But as
$\dim H^1(C;
\frak u(\xi_{x,P})) = r+1$, $\dim H^1(C; \ad_G(\xi_{x,P}\times_PG)) \geq r$,
and the homomorphism is constant on the image of
$H^0(C; \ad_L\eta_0)$, it follows that $\dim H^1(C; \ad_G(\xi_{x,P}\times_PG))
= r$. Hence
$\xi_{x,P}\times_PG$ is regular, and the Kodaira-Spencer homomorphism defines
an isomorphism from the tangent space of $\mathbb{S}$ at $x$ to $H^1(C;
\ad_G(\xi_{x,P}\times_PG))$. The rest of the argument is identical to the
proof of Theorem~\ref{equivisom}.
\end{proof}

In particular, the $\Cee^*$-weights on $\mathbb{S}$ are given by the Casimir
weights of $E_8$. Thus, given a fibration $Z\to B$, when we write the
$10$-dimensional bundle of affine spaces as a direct sum of line bundles,
normalized to begin with the factor
$\scrO_B$, the remaining line bundle summands are of the form
$\mathcal{L}^{-d_i}$, where the $d_i$ are the Casimir weights, together
with $\mathcal{L}^{-6}$, where $-6$ is the suitably normalized weight of the
unipotent group scheme, thought of as a line bundle over the moduli stack.

Clearly,
similar constructions with appropriate unipotent subgroups $U'\subseteq U$ can
be made with groups other than
$E_8$. For example, for $E_6$, if we again look at the weight one summand in
$\frak u$ and use the quotient corresponding to the subbundle
$W_3\spcheck
\otimes W_3\spcheck \subseteq W_2\otimes W_3\spcheck
\otimes W_3\spcheck$, the result is a 10-dimensional affine space with the
action of a 3-dimensional unipotent group. These actions seem closely
connected to the equations defining the corresponding del Pezzo surfaces, and
will be discussed in more detail in \cite{FMIV}. For example, in the case of
$E_8$, the 9-dimensional weighted projective space seems very related to the
equations of rational elliptic surfaces with a given fiber, and in the case
of $E_6$, the affine space mentioned above seems connected with the space of
cubic equations with a given hyperplane section.

\section{Unstable bundles and minuscule representations}

Throughout this section, we assume that $C$ is a singular Weierstrass cubic and
that $\iota\colon \widetilde C \to C$ is the normalization.  As we mentioned in
the introduction, there is no general definition for semistability of a $G$-bundle
over a singular curve. Of course, for any reasonable definition of semistability,
it is natural to require that a
$G$-bundle  whose pullback  to the normalization is semistable should itself be
semistable. But the set of all such bundles will never give a compact moduli
space. For $SL_n$, we can define semistable vector bundles by the usual slope
definition applied to torsion free subsheaves, and, in the case of a
Weierstrass cubic, these will form a compact moduli space. On the other hand,
tensor operations, for example exterior powers, applied to semistable
vector bundles do not always yield semistable vector bundles. In this section,
we answer the following question: given a
$G$-bundle $\xi$ arising from the parabolic construction and an irreducible
representation
$\rho\colon G \to GL_N(\Cee)$, when is the vector bundle $\xi\times_G\Cee^N$
semistable? We may as well restrict attention to the divisor $D_\infty$
defined in Corollary~\ref{Dinfty}, since for $\xi$ corresponding
to a  point of
$\WP(G) -D_\infty$, the vector bundle
$\xi\times_G\Cee^N$ is semistable for every $\rho$. As we shall see, over
$D_\infty$, the  $\rho$ for which 
$\xi\times_G\Cee^N$ is semistable are in a very short list of ``small"
representations of $G$. Even for such $\rho$, the question is rather delicate.

\subsection{A preliminary reduction}

For the definitions of minuscule and quasi-minuscule, we refer to
\cite{Bour}. The following result shows that, if
$\xi$ is a
$G$-bundle on
$C$ coming from the parabolic construction such that $\iota^*\xi$ is not
trivial, then for almost every irreducible representation $\rho\colon G \to
GL_N(\Cee)$, the vector bundle
$\xi \times _G \Cee^N$ is unstable. 

\begin{lemma}\label{prelimred} Let $\rho\colon G \to GL_N(\Cee)$ be an
irreducible representation of
$G$ with highest weight $\varpi$. Let $\xi$ be a $G$-bundle arising from the
parabolic construction, and suppose that $\xi$ corresponds to a point in the
divisor
$D_\infty$, i.e.\ the pullback $\iota^*\xi$ is not the trivial $G$-bundle. If
$\xi \times _G \Cee^N$ is a semistable vector bundle, then, for all short coroots
$\beta\spcheck$, $\varpi(\beta\spcheck) \leq 1$. In particular, $\varpi$ is a
fundamental weight $\varpi_\delta$ for some $\delta \in \Delta$, and the
coefficient
$g_\delta$ of
$\delta
\spcheck$ in the coroot $\widetilde \alpha\spcheck$ dual to the highest root is
one. Finally, if $\varpi_\delta$ is minuscule, then $g_\delta=1$, and if $G$
is simply laced and $g_\delta=1$, then  conversely  $\varpi_\delta$ is
minuscule.
\end{lemma}
\begin{proof} It follows from
Lemma~\ref{whenunstable} that, if $V$ is a semistable vector bundle of
degree zero on
$C$ and
$\iota^*V \cong \bigoplus _i\scrO_{\Pee^1}(a_i)$, then $a_i\leq 1$ for every 
$i$.

By Corollary~\ref{onedim} and Lemma~\ref{lambda1}, if $\iota^*\xi$ is not the
trivial
$G$-bundle, then 
$\iota^*\xi \cong \iota^*\eta_0\times_LG =
\gamma_{-\lambda_1(\alpha)\spcheck}\times _HG$, where 
$\lambda_1(\alpha)\spcheck=\sum _{\beta \in \Delta}\beta\spcheck$ is the
highest coroot such that the coefficient of
$\alpha\spcheck$ in $\lambda_1(\alpha)\spcheck$ is one. Hence, the pullback of
$\xi \times _G \Cee^N$ to $\widetilde C\cong \Pee^1$ contains as direct summands
bundles of the form $\scrO_{\Pee^1}(w\varpi(-\lambda_1(\alpha)\spcheck))$, for
$w\in W$. Thus, there are summands of the form $\scrO_{\Pee^1}(n)$ for every
integer $n$ of the form $\varpi(w^{-1}(-\lambda_1(\alpha)\spcheck)$. On the
other hand, $-\lambda_1(\alpha)\spcheck$ is a short coroot, and so its orbit
under $W$ is the set of all short coroots. Hence $\varpi(\beta\spcheck) \leq 1$
for every short coroot $\beta\spcheck$. In particular, $\varpi(\widetilde
\alpha\spcheck)
\leq 1$.

Write $\varpi = \sum _{\delta\in \Delta}n_\delta\varpi_\delta$, where
$\varpi_\delta$ is the fundamental weight corresponding to $\delta$. Then 
$$1\geq \varpi(\widetilde\alpha\spcheck) =\sum _{\delta \in \Delta}
n_\delta g_\delta.$$
It follows that $n_\delta $ is nonzero for exactly one $\delta \in \Delta$, and
for this $\delta$ we have $n_\delta= g_\delta =1$. The final statement is then
clear.
\end{proof} 

For example, in case $G$ is of type $E_8$, $C$ is nodal,  and $\xi$ is a
$G$-bundle arising from the parabolic construction which does not pull back to
the trivial bundle on
$\widetilde C$, then $\xi\times _G\Cee^N$ is unstable for every irreducible
representation $\rho\colon G \to GL_N(\Cee)$. In fact, using the results of
Ramanathan \cite{Ra2} (or directly), the same is true for every $G$-bundle
$\xi$ which does not pull back to the trivial bundle on
$\widetilde C$.

The representations satisfying the conclusions of Lemma~\ref{prelimred} are
easy to tabulate:

\begin{lemma} Suppose that $G$ is not  simply laced. Then the 
non-minuscule representations with highest weight
$\varpi_\delta$,
$\delta\in
\Delta$ such that $g_\delta =1$  are as follows: they are either the
quasi-minuscule representations, in case
$G$ is of type $B_n$, $F_4$, or $G_2$, or any of the fundamental
representations in case
$G$ is of type $C_n$. \qed
\end{lemma}

Note that, in every case but type $C_n$, the highest weight
$\varpi$  of a quasi-minuscule representation is the highest short root, and that
$\varpi$ is the fundamental weight dual to a simple coroot corresponding to an
endpoint of the Dynkin diagram. The same is true in the minuscule case if $G$
is   not of type $A_n$, or if $G=SL_n(\Cee)$ and $\rho$ is either the
standard representation or its dual. Hence, if we are not in one of these
exceptional cases and if
$G'\subseteq G$ is the simple subgroup corresponding to the root system
defined by the kernel of
$\varpi_\delta$, or equivalently whose simple roots are $\Delta -\{\delta\}$,
then $G'\subseteq G$ is one of the inclusions of groups arising via parabolic
induction. Recall that, in this case, there is the hypersurface
$D_{G'}\subseteq
\WP(G)$ arising from parabolic induction.

\subsection{The minuscule case}

In this subsection, we assume that $\rho$ is minuscule. Our goal is to show
the following inductive result:

\begin{theorem}\label{minunst}  Suppose that $G$ is not of type $A$. Let
$\rho\colon G
\to GL_N(\Cee)$ be a minuscule representation of $G$, with highest weight
$\varpi_\delta$, and let
$G'\subseteq G$ be the simple subgroup whose simple roots are $\Delta
-\{\delta\}$. Suppose that $x\in
\WP(G)$ corresponds to the $G$-bundle $\xi$, where $\xi\times_G\Cee^N$ is
unstable. Then $x\in D_{G'}$ is the image of $x'\in \WP(G')$. Moreover, if
$\xi'$ is the 
$G'$-bundle corresponding to $x'$, then there exists an irreducible
subrepresentation 
$\rho'\colon G'\to GL_{N'}(\Cee)$ of
$\rho|G'$, necessarily minuscule as well, such that $\xi'\times _{G'}\Cee^{N'}$
is unstable. Conversely, if $\xi$ is obtained via parabolic induction from such a
$\xi'$, then $\xi\times_G\Cee^N$ is
unstable.
\end{theorem}
\begin{proof} If $\xi\times_G\Cee^N$ is
unstable, then $H^0(C; \xi\times_G\Cee^N) \neq 0$. The following result, which is
of independent interest and also holds in case $C$ is smooth,  then tells us that
$x\in D_{G'}$:

\begin{theorem}\label{Deerho} Suppose that $\rho\colon G
\to GL_N(\Cee)$ is a minuscule representation of $G$, with highest weight
$\varpi_\delta$ and let
$G'\subseteq G$ be the simple subgroup whose simple roots are $\Delta
-\{\delta\}$.
Then the divisor
$D_{G'}$ is exactly the set of $x\in \WP(G)$ such that, if $\xi_x$ is the
corresponding
$G$-bundle, then
$H^0(C; \xi_x\times_G\Cee^N) \neq 0$.
\end{theorem}

Assuming Theorem~\ref{Deerho}, let us finish the proof of the first two
conclusions of Theorem~\ref{minunst}. It suffices to show that, if 
$\xi\times_G\Cee^N$ is unstable and if $\xi'$ is the
corresponding
$G'$-bundle, then there exists an irreducible subrepresentation  $\rho'\colon
G'\to GL_{N'}(\Cee)$ of
$\rho|G'$  such that $\xi'\times _{G'}\Cee^{N'}$
is unstable. Clearly, it suffices to prove that $(i_*\xi')\times_G\Cee^N$ is
unstable. By  Lemma~\ref{limit}, there is a family of
$G$-bundles $\Xi \to C\times
\Cee$ such that, for $t\neq 0$, $\Xi_t\cong \xi$, and such that $\Xi_0 \cong
i_*\xi'$. It follows  that the vector bundle $(i_*\xi')\times_G\Cee^N$ is
the limit of the unstable bundles $\xi\times_G\Cee^N$. Hence it is unstable.

Finally, the last statement of Theorem~\ref{minunst} follows from the fact that,
by Proposition~\ref{Thm3}, there is a filtration of $\xi\times_G\Cee^N$ whose
associated graded is the bundle
$\xi'\times _{G'}\Cee^N$. Moreover, by Part (iii) of Proposition~\ref{Thm3},
the successive   summands of the associated graded are associated to
representations of $G'$ and thus have degree zero, since $G'$ is simple.
Since  at least one summand in the associated graded is unstable and of
degree zero, it follows that $\xi\times_G\Cee^N$ is unstable.
\end{proof}

\noindent \textbf{Proof of Theorem~\ref{Deerho}.} Let $D_\rho$ be the set of $x\in
\WP(G)$ such that
$H^0(C;
\xi_x\times_G\Cee^N) \neq 0$.  The argument will be in the following steps: First
we show that $D_\rho$ is a hypersurface containing $D_{G'}$. Then we show that
$D_\rho$ is irreducible, by showing that $D_\rho \cap (\WP(G) -D_\infty)$ is
irreducible and that $D_\rho$ does not contain the irreducible hypersurface
$D_\infty$. It then follows that $D_\rho$ is irreducible, and hence equal to
$D_{G'}$.

We first claim  that $D_\rho$ is a either a hypersurface, is all of
$\WP(G)$,  or is empty. It suffices to work in
$\Aff^{r+1}$ with the corresponding $\Cee^*$-invariant subvariety. There is a
universal bundle over
$C\times \Aff^{r+1}$, and the inverse image of $D_\rho$ is clearly the 
divisor associated to the canonical section of the inverse of the determinant line
bundle of this family. Thus, either
$D_\rho$ is empty, is all of
$\WP(G)$,  or  is a hypersurface. 

Next, we claim that, under the identification of $\WP(G) -D_\infty$ with 
$(C_{\rm reg}\otimes\Lambda)/W$, $D_\rho \cap (\WP(G)
-D_\infty)$ is the image of $C_{\rm reg}\otimes\Lambda'$, where $\Lambda'$ is
the coroot lattice of $G'$. In particular, this implies that $D_\rho \cap
(\WP(G) -D_\infty)$ is irreducible.  If
$t\in C_{\rm reg}\otimes
\Lambda$ is a lift of a point of $(C_{\rm reg}\otimes
\Lambda)/W$ corresponding to $\xi$, the corresponding vector bundle $\xi\times
_G\Cee^N$  is semistable, with support the divisor in $C_{\rm reg}$ given by $\sum
_{\lambda}\lambda(t)$, where the sum is over the weights of $\rho$. In
particular, the origin $p_0$ is in the support of $\xi\times
_G\Cee^N$ if and only if $H^0(C; \xi\times
_G\Cee^N) \neq 0$ if and only if $\lambda(t) = 0$ for some weight $\lambda$ of
$\rho$. Since $\rho$ is minuscule, all weights are conjugate under $W$, so the
hypersurface $D_\rho \cap (\WP(G) -D_\infty)$ is is the image in  $(C_{\rm
reg}\otimes\Lambda)/W$ of the kernel of $\lambda \colon C_{\rm reg}\otimes
\Lambda \to C_{\rm reg}$, where $\lambda$ is any fixed weight of $\rho$, for
example $\lambda =\varpi_\delta$.  So it suffices to show that the kernel of 
$\varpi_\delta
\colon C_{\rm reg}\otimes \Lambda \to C_{\rm reg}$ is $C_{\rm
reg}\otimes\Lambda'$. But $\varpi_\delta$ is a fundamental weight and hence
primitive, and its kernel  is $\Lambda'$, the coroot lattice of $G'$.
Moreover the sequence
$$0 \to \Lambda' \to \Lambda \xrightarrow{\varpi_\delta} \Zee \to 0$$ is exact.
Hence the kernel of $\varpi_\delta \colon C_{\rm reg}\otimes
\Lambda \to C_{\rm reg}$ is  $C_{\rm reg}\otimes \Lambda'$, and so $D_\rho
\cap  (\WP(G) -D_\infty)$ is irreducible. It also follows that $D_\rho \cap 
(\WP(G) -D_\infty) = D_{G'}\cap (\WP(G) -D_\infty)$. Since $D_{G'}$ is
irreducible,   $D_{G'}\subseteq D_\rho$.

Finally, we claim that $D_\rho$ is itself irreducible and hence is equal to
$D_{G'}$. It suffices to show that
$D_\rho$ does not contain the irreducible hypersurface
$D_\infty$. The proof is by induction on the rank of $G$. In order to carry out
the induction, we shall have to consider groups of type $A$, since it is
possible that the group $G'$ in the above notation is of type $A$. If the
rank of
$G$ is one, then $G\cong SL_2(\Cee)$, there is a unique minuscule representation,
the standard one, $D_\infty$ is a point, and the corresponding vector bundle,
which is 
$V_{1,1}$ if $C$ is nodal and $V_2$ if $C$ is cuspidal, is semistable and
supported at the singular point of
$C$. Hence
$H^0(C; V)=0$. Now suppose that the rank of $G$ is at least $2$. If $G$ is of
type $A_{n-1}$, then it is easy to see that there are vector bundles of the
form 
$$V=(V_0\otimes \lambda)\oplus \bigoplus _{i=1}^{n-2}\lambda_i,$$
where $V_0$ is either $V_{1,1}$ or $V_2$ and  $\lambda$ and $\lambda_i$
are general line bundles of degree zero such that
$\lambda^2\otimes
\lambda_1\otimes \cdots \otimes \lambda_{n-2}$ is trivial, satisfying:
$h^0(C;\bigwedge ^kV) = 0$ for all $k$ such that $1\leq k \leq n-1$. In
particular, such bundles define points in $D_\infty - D_\rho$ for every
representation $\rho$ corresponding to a minuscule (i.e.\ fundamental) weight
of $SL_n(\Cee)$.

Thus we can
apply the inductive hypothesis to the simple subgroup
$G'$ of $G$ defined by $\rho$. We view $\rho$ as a representation of $G'$ by 
restriction. For every $G'$-bundle $\xi'$, the bundle $\xi'\times _{G'}\Cee^N$ is
a direct sum of bundles of the form $\xi'\times _{G'}\Cee^{N_i}$ associated to
minuscule representations of
$G'$ and trivial factors, corresponding to weights of $\rho$ which are trivial on
the span of the simple coroots in $G'$. Choose a $G'$-bundle $\xi'$, obtained by
the parabolic construction for $G'$, such that $\iota^*\xi'$ is not the trivial
bundle and such that, for every minuscule representation $\rho_i \colon G' \to
GL_{N_i}(\Cee)$, $H^0(C; \xi'\times _{G'}\Cee^{N_i}) = 0$. Let $\xi$ be the
$G$-bundle obtained from $\xi'$ by parabolic induction. We first claim that
$\iota^*\xi$ is not trivial, so that $\xi$ corresponds to a point of $D_\infty$:

\begin{lemma}\label{Thislemma} Let $\xi'$ be a $G'$-bundle, obtained by the
parabolic construction for $G'$, such that $\iota^*\xi'$ is not the trivial
$G'$-bundle, and let  $\xi$ be the
$G$-bundle obtained from $\xi'$ by parabolic induction. Then $\iota^*\xi$ is not
the trivial $G$-bundle.
\end{lemma}
\begin{proof} The hypothesis that $\iota^*\xi'$ is not trivial easily implies
that, for every nontrivial finite-dimensional representation  of
$G'$ on a vector space $V$, $\iota^*\xi'\times _{G'}V$ is an unstable bundle of
degree zero over
$\Pee^1$. Suppose that $\iota^*\xi$ is trivial. In particular, if we choose some
nontrivial finite-dimensional representation $\tau\colon G \to GL(V)$, then
$\iota^*\xi\times_GV$ is the trivial bundle.  By Proposition~\ref{Thm3}, there
is a filtration of $\iota^*\xi\times_GV$ whose associated graded is the bundle
$\iota^*\xi'\times _{G'}V$. Moreover, by (iii) of Proposition~\ref{Thm3} the
successive   summands of the associated graded are associated to
representations of $G'$ and thus have degree zero, since $G'$ is simple.
Since $V$ is nontrivial as a representation of
$G'$, at least one summand in the associated graded is unstable and of degree
zero, which is impossible. Hence $\iota^*\xi$ is nontrivial.
\end{proof}

Because of the trivial factors in the representation $\rho|G'$, it is always the
case that $H^0(C; \xi'\times _{G'}\Cee^N) \neq 0$. However, we can modify the
construction as follows. Let $\lambda$ be a line bundle of degree zero on $C$,
which we identify with a holomorphic principal $\Cee^*$-bundle, and denote by
$\xi'\otimes \lambda$ the $(G'\times_F\Cee^*)$-bundle induced by the pair $(\xi',
\lambda)$. Clearly, $\iota^*(\xi'\otimes \lambda) \cong \iota^*\xi'$, so that
$\iota^*(\xi'\otimes \lambda)$ is not the trivial $(G'\times_F\Cee^*)$-bundle. 
Note that
$\xi'\otimes
\lambda$ has a holomorphic reduction of structure to $\ov QU'$, covering the
bundle
$\eta_0'\otimes\lambda$. 

Since $\Cee^*$ is in the center of $G'\times _FG$, it acts on each irreducible
summand $\Cee^{N_i}$ of $\Cee^N$ via  a character $\chi_i$. Moreover, if
$\Cee^{N_i}$ is a trivial representation of $G'$, then $\chi_i$ is not the trivial
character of $\Cee^*$, because no nonzero element of $\Cee^N$ is fixed by $H$.
Thus, 
$$(\xi'\otimes \lambda)\times _{(G'\times_F\Cee^* )}\Cee^N\cong \bigoplus
_i(V_i\otimes \lambda^{n_i}),$$
where the $V_i$ are vector bundles over $C$ such that either $H^0(C; V_i) =0$ or
$V_i\cong \scrO_C$ and $n_i\neq 0$. It follows that, for generic $\lambda$,
$H^0(C; (\xi'\otimes \lambda)\times _{(G'\times_F\Cee^*) }\Cee^N) =0$.

Now let us show that we can carry out a modified version of   parabolic induction
for $\xi'\otimes \lambda$. First note that $\xi'\otimes \lambda$ is a lift to the
parabolic subgroup $P'\times _F\Cee^*$ of $G'\times _F\Cee^*$ of the bundle
$\eta_0'\otimes
\lambda$. Let $\xi_\lambda$ be a $QU'$-bundle lifting $\xi'\otimes \lambda$. The
argument of Lemma~\ref{Thislemma} shows that $\iota^*\xi_\lambda$ is not the
trivial $G$-bundle for any such lift $\xi_\lambda$. It now suffices to prove that
there is a lift
$\xi_\lambda$  given by the parabolic construction for $G$, at least for
$\lambda$ close to the trivial bundle, for then $\xi_\lambda$ satisfies:
$\xi_\lambda \in D_\infty -D_\rho$, as desired. By Corollary~\ref{KSisom}, it is
enough to show that, at least for
$\lambda$ close to the trivial bundle, there is a lift $\xi_\lambda$ which is a
small deformation of $\xi$. First, it is easy to check directly or via
deformation theory that there is a family $\eta_\lambda$ of $QU'$-bundles, such
that $\eta_0$ agrees with the usual definition of $\eta_0$, and such that
$\eta_\lambda/N \cong \eta_0'\otimes \lambda$. Since $H^0(C; \frak u(\eta_0)) =
0$,  $H^0(C; \frak u(\eta_\lambda)) =
0$ for generic values of $\lambda$. Hence, for such $\lambda$, there is a bundle
of affine spaces parametrized by $\lambda$ which restricts for each $\lambda$ to
give the affine space representing $H^1(C; U'(\eta_\lambda))$.  On the other
hand, $H^1(C; U'(\eta_\lambda)) \cong H^1(C; U'(\eta_0'\otimes \lambda))$.
The map $\lambda \mapsto \xi'\otimes \lambda$ defines a section of the affine
space bundle corresponding to $H^1(C; U'(\eta_0'\otimes \lambda))$. Let
$\xi_\lambda$ be the corresponding section of the affine bundle corresponding
to $H^1(C; U'(\eta_\lambda)) $. For $\lambda$ close to the trivial line
bundle, the bundle
$\xi_\lambda$ is a small deformation of
$\xi$.  The arguments of (ii) of Theorem~\ref{Thm1} show that
$\xi_\lambda/N\cong \xi'\otimes
\lambda$, and so $\xi_\lambda$ is the desired lift.
\qed

\subsection{The quasiminuscule case}

In this subsection, in case $G$ is of type $E_8$, we  assume that
$C$ is not cuspidal, and shall let $\rho\colon G \to GL_N(\Cee)$  be a
quasiminuscule, non-minuscule representation.  Thus
$\rho$ is the adjoint representation if
$G$ is simply laced, and in general the highest weight of $\rho$ is the highest
short root. Let $s$ be the multiplicity of the trivial weight in $\rho$, so that
$s=r$ if $G$ is simply laced, and in general $s$ is the number of short simple
roots.

The first lemma shows that, by contrast to the minuscule case, the ranks of the
associated vector bundles in the quasiminuscule case are constant:

\begin{lemma} For every $G$-bundle $\xi$  coming from the parabolic
construction, $$\dim H^0(C; \xi \times _G\Cee^N) = s.$$
\end{lemma}
\begin{proof} In case $\rho$ is the adjoint representation, this follows from
Corollary~\ref{KSisom}. Thus we may assume that $G$ is not simply laced. On
the open dense subset of
$\WP(G)$ corresponding to regular semisimple elements of $(C_{\rm reg}\otimes
\Lambda)/W$, it is clear that $\dim H^0(C; \xi \times _G\Cee^N) = s$. Thus by
semicontinuity
$\dim H^0(C; \xi \times _G\Cee^N) \geq s$ for all $\xi$ arising from the
parabolic construction. On the other hand, viewing $\rho$ as a representation of
$L$ by restriction, a case-by-case argument  shows that $h^0(C;
\eta_0\times_L\Cee^N) = s$. Thus, for all $\xi$ corresponding to points of the
affine space
$H^1(C; U(\eta_0))$ close to the origin, $h^0(C;
\eta_0\times_L\Cee^N) \leq s$. Using the $\Cee^*$-action, the same statement
holds for all $\xi$. Combining the two inequalities, we see that $h^0(C;
\eta_0\times_L\Cee^N) = s$ for all $\xi$.
\end{proof}

For any vector bundle $V$, we have the natural map $H^0(C;V) \otimes
\scrO_C \to V$ and hence an induced map $\psi_V\colon H^1(C; H^0(C;V) \otimes
\scrO_C) \to H^1(C; V)$. For example, if $V$ is semistable and of degree zero,
then
$\psi_V$ is an isomorphism if and only if
$V$ has no summands of the form $I_n$ for $n>1$, where $I_n$ is the unique
indecomposable bundle which has a filtration whose associated graded is
$\scrO_C^n$. 

In particular, there is a natural map
$\psi_\xi\colon H^1(C;
\scrO_C^s)
\to H^1(C;\xi
\times _G\Cee^N)$.

\begin{proposition}\label{isadivisor} In the above notation, 
\begin{enumerate}
\item[\rm (i)] Suppose that $\chi(C;V) =0$. Then
$\psi_V$ is an isomorphism if and only if $V$ is isomorphic to $\scrO_C^s\oplus
V'$, where
$V'$ is a vector  bundle such that $H^0(C; V') =0$.
\item[\rm (ii)] Let $\mathcal{V}\to C\times T$ be a family of vector bundles of
index zero over $C$, and suppose that $h^0(C; V_t)=s$ is independent of $t\in
T$. Then
$\mathcal{D}=\{t\in T:
\psi_{V_t}
\text{ is not an isomorphism}\}$ is either a divisor in $T$, is all of $T$, or is
empty.
\end{enumerate}
\end{proposition} 
\begin{proof} Clearly, if $V\cong \scrO_C^s\oplus
V'$, where  $H^0(C; V') =0$, then $\psi_V$ is an isomorphism. Conversely, suppose
that $\psi_V$ is an isomorphism and define $V'$ by the exact sequence
$$0 \to \scrO_C^s\to V \to V'\to 0,$$
where $s=h^0(C;V)$.
Since by construction $H^0(C; \scrO_C^s) \to H^0(C;V)$ is an isomorphism, we see
that $\psi_V$ is an isomorphism if and only if $H^0(C;V') = H^1(C; V') = 0$. In
particular, by Serre duality $\Ext^1(V', \scrO_C)$ is dual to $H^0(C;V')$ and
hence is zero. It follows that $V\cong \scrO_C^s\oplus V'$, and in particular
$V'$ is locally free. This proves (i).

To prove (ii), consider the parametrized analogue of the above exact sequence:
$$\scrO_{C\times T}^s \to \mathcal{V} \to \mathcal{V}' \to 0.$$
By standard base change results, restricting to every fiber $C\times \{t\}$  gives
the exact sequence $0 \to \scrO_C^s\to V_t \to V'_t\to 0$. The proof of the first
part shows that $\psi_{V_t}$ fails to be an isomorphism if and only if $h^0(C;
V_t')
\neq 0$. By the local criterion of flatness, the map
$\scrO_{C\times T}^s
\to
\mathcal{V}$ is injective and
$\mathcal{V}'$ is flat over $T$. Since each $V'_t$ is of index zero, the set of
$t$ such that $h^0(C; V_t')
\neq 0$ is then a divisor in $T$, is all of $T$, or is
empty, as required.
\end{proof}

\begin{defn} Let $\mathcal{D}_\rho\subseteq \WP(G)$ correspond to the set
of bundles 
$\xi$ such that the corresponding homomorphism $\psi_\xi \colon H^1(C; \scrO_C^s)
\to H^1(C; \xi
\times _G\Cee^N)$ is not an isomorphism.
\end{defn} 

A different characterization of $\mathcal{D}_\rho$ is as follows:

\begin{proposition} The bundle $\xi$ corresponds to a point of $\WP(G)
-\mathcal{D}_\rho$ if and only if the vector bundle $\xi \times _G\Cee^N$ is
isomorphic to $\scrO_C^s\oplus V'$, where $V'$ is a vector bundle such that
$H^0(C; V') =0$. In particular, if $\xi \times _G\Cee^N$ is unstable, then $\xi$
corresponds to a point of $\mathcal{D}_\rho$.
\end{proposition} 
\begin{proof} The first statement follows from Proposition~\ref{isadivisor}. To
see the second, if $\xi\notin\mathcal{D}_\rho$, then $\xi \times _G\Cee^N
\cong \scrO_C\oplus V'$, where $h^0(C;V') = 0$. In particular $V'$ is
semistable. Thus $\xi
\times _G\Cee^N$ is semistable as well.
\end{proof}

With this said, there is the following analogue of Theorem~\ref{Deerho}:

\begin{theorem}\label{calDeerho} Let $\rho$ be a quasiminuscule representation
of $G$, and let $G'$ be the stabilizer of a nonzero weight of $\rho$. Then
$\mathcal{D}_\rho$ is a hypersurface in $\WP(G)$. If in addition $G$ is
not simply laced or of type $C_n$, so that one can do parabolic induction from
$G'$ to $G$, then
$\mathcal{D}_\rho= D_{G'}$, and in particular it is irreducible.
\end{theorem}
\begin{proof}  In this case, $\mathcal{D}_\rho\cap (\WP(G) -D_\infty)$ is
exactly the locus corresponding to the points in $(C_{\rm reg}\otimes\Lambda)/W$
where a weight (i.e.\ a short root) vanishes, and hence $\mathcal{D}_\rho\cap
(\WP(G) -D_\infty) = D_{G'}\cap (\WP(G) -D_\infty)$. Given
this, minor modifications of the proof of Theorem~\ref{Deerho} handle this case
as well.
\end{proof} 

\begin{remark} In case $G$ is of type $C_n$, it is easy to check that the divisor
$\mathcal{D}_\rho$ is still irreducible. On the other hand, if $G$ is simply
laced, then $\mathcal{D}_\rho$ always contains $D_\infty$ and is thus reducible,
with two components, except for $G=SL_2(\Cee)$ and $C$ nodal, where it consists
of three points, one of which has multiplicity two.
\end{remark} 

Arguments as in the proof of Theorem~\ref{minunst} then show:

\begin{corollary} Suppose that $G$ is not simply laced or of type $C_n$, and that
$\rho\colon G \to GL_N(\Cee)$ is a quasiminuscule representation of $G$. Let $G'$
be the corresponding subgroup of $G$. Let $x\in \WP(G)$ be a point
corresponding to the bundle $\xi$, and suppose that $\xi\times _G\Cee^N$ is
unstable. Then $x\in D_{G'}$, i.e.\ $\xi$ is obtained via parabolic induction from
a
$G'$-bundle
$\xi'$. Moreover, there exists an irreducible subrepresentation
$\rho'\colon G' \to GL_{N'}(\Cee)$, necessarily
quasiminuscule, of $\rho|G'$ such that the bundle $\xi'\times _{G'}\Cee^{N'}$ is
unstable. Conversely, for every $G'$-bundle $\xi'$, obtained via the parabolic
construction for $G'$ and with this property, if $\xi$ is the $G$-bundle obtained
via parabolic induction from such a $\xi'$, then $\xi\times_G\Cee^N$ is
unstable.\qed
\end{corollary}

\subsection{Examples}

\noindent\textbf{The case of $\boldsymbol{SL_n(\Cee)}$}

In this case, for the standard representation, the set of vector bundles
associated to the parabolic construction is exactly the set of regular semistable
bundles \cite[Theorem 3.2]{FMW}. In particular, every such vector bundle is
semistable. For the standard inclusion of
$SL_{n-1}(\Cee)
\subseteq SL_n(\Cee)$, the image of parabolic induction corresponds to the set of
rank $n$ regular semistable vector bundles $V$ of degree zero such that $h^0(C;
V) \geq 1$. For such bundles
$V$, there is a (two step) filtration on
$V$, $\scrO_C\subseteq V$, with associated graded $\scrO_C \oplus V'$, where $V'$
is the corresponding regular
semistable vector bundle of rank $n-1$. To deal with the remaining minuscule
representations, we have the following, which is an immediate consequence
of Corollary~\ref{wedgeunstable}  and  Corollary~\ref{wedgesemistable}.

\begin{proposition}\label{sln} Suppose that $C$ is nodal. Let $V$ be a regular
semistable vector bundle on $C$ with $\det V \cong\scrO_C$, and write $V\cong
V_{a,b;\lambda} \oplus V'$, where $\lambda\in \Cee^*$ and $V'$ has no support at
the singular point. If either
$a$ or $b$ is $1$, then $\bigwedge^kV$ is semistable for all $k$. If $a,b\geq 2$,
then $\bigwedge^kV$ is unstable for all $k$, $2\leq k\leq n-2$. Similarly,
suppose that
$C$  is cuspidal and that $V\cong
V_{a; \lambda} \oplus V'$, where $\lambda\in \Cee$ and $V'$ has no support at
the singular point.  If $a\leq 3$ then  $\bigwedge^kV$ is semistable for all
$k$, and if
$a\geq 4$, then $\bigwedge^kV$ is unstable for all $k$, $2\leq k\leq n-2$.\qed
\end{proposition}

\medskip
\noindent\textbf{The case of $\boldsymbol{Sp(2n)}$}

The parabolic construction for $Sp(2n)$-bundles is a special subset of 
extensions given by the parabolic construction for $SL_{2n}(\Cee)$-bundles,
coming from the inclusion 
$$H^1(C; \Sym^2W_n\spcheck)\subseteq H^1(C;
W_n\spcheck\otimes W_n\spcheck).$$
In particular, every associated vector
bundle in the standard representation is semistable, and parabolic induction
from $Sp(2n-2)$ corresponds to the set of regular semistable symplectic
bundles $V$ of rank $2n$ such that
$H^0(C;V) \neq 0$. In this case, there is a self-dual complex 
$$0 \to \scrO_C \xrightarrow{f} V \xrightarrow{g}  \scrO_C \to 0,$$
such that $\Ker g/\operatorname{Im} f\cong V'$ is a  regular semistable
symplectic bundle of rank $2n-2$. Finally, those bundles in $D_\infty$ are
necessarily of the form $V_{a,a} \oplus V'$, where $V'$ is a regular
semistable symplectic bundle of rank $2n-2a$ with no support at the singular
point, if $C$ is nodal, or $V_{2a} \oplus V'$, where $V'$ is a regular
semistable symplectic bundle of rank $2n-2a$ with no support at the singular
point, if $C$ is cuspidal.

Given a representation $\rho$ satisfying the conditions of
Lemma~\ref{prelimred} and a $Sp(2n)$-bundle $\xi$ coming from the parabolic
construction, we must decide when the associated vector bundle is unstable.
The representations $\rho$ in question are exactly those whose highest weight
is a   fundamental weight. Thus, they are given as follows: Let $\omega \in
\bigwedge ^2\Cee^{2n}$ be the symplectic form. Then, for all $k$ with $1\leq
k\leq n$, $Sp(2n)$ acts on $\bigwedge ^k\Cee^{2n}$, preserving the subspace
$\omega
\wedge \bigwedge ^{k-2}\Cee^{2n}$,  the quotient $\bigwedge ^k\Cee^{2n}/\omega
\wedge
\bigwedge ^{k-2}\Cee^{2n}$ is an irreducible $Sp(2n)$-module with highest weight
a fundamental weight, and these are exactly a set of fundamental representations.

\begin{proposition} Let $\xi$ be a principal $Sp(2n)$-bundle arising from the
parabolic construction such that the associated rank $2n$-bundle 
$V$ is of the form $V_{a,a}
\oplus V'$, where
$V'$ is a  regular semistable symplectic bundle of rank $2n-2a$ with no 
support at the singular point, if $C$ is nodal, or $V_{2a} \oplus V'$, where
$V'$ is a regular semistable symplectic bundle of rank $2n-2a$ with no
support at the singular point, if $C$ is cuspidal. Let $\rho$ be a
fundamental representation of
$Sp(2n)$. Then the vector bundle associated to $\xi$ via $\rho$ is semistable 
if
$a=1$, and is unstable if $a> 1$ and $\rho$ is not the standard
representation.
\end{proposition}
\begin{proof} The proof is clear if $a=1$. To prove the result in case $a>1$, 
it suffices to consider the nodal case, since semistability is an open
condition in families of curves. Clearly, it is enough to show the following:
for $a>1$ and  for all
$k$ with $2\leq k\leq a$, the quotient 
$$\left(\bigwedge ^kV_{a,a}\right)\Bigg/\omega
\wedge\left( \bigwedge
^{k-2}V_{a,a}\right)$$   is unstable. This can be proved using the methods of
Lemma~\ref{ssetale} and Proposition~\ref{wedgess}. We shall omit the somewhat
tedious details.
\end{proof}

\medskip
\noindent\textbf{The case of $\boldsymbol{Spin(2n+1)}$}

In this case, there is the spin representation, which is minuscule, and the
standard ($SO$) representation, which is quasi-minuscule. We consider these
separately.

In the case of the spin representation $\rho$, the group $G'$ which stabilizes a
weight is the subgroup $SL_n(\Cee)$. The restriction of $\rho$ to $SL_n(\Cee)$
defines the representation $\bigwedge^*\Cee^n$, the direct sum of the exterior
powers of the standard representation of $SL_n(\Cee)$. Moreover, the
restriction of the
$SO$-representation to $G'$ is $\Cee^n\oplus (\Cee^n)^*\oplus \Cee$, the
direct sum of the standard representation of $SL_n(\Cee)$, its dual, and a
trivial summand. Thus we see:

\begin{proposition}\label{oddspin} Given  $x\in
\WP(Spin(2n+1))$, let  $\xi$ be the $Spin(2n+1)$-bundle corresponding to $x$
and let
$S$ be the vector bundle associated to $\xi$ by the spin representation. Then
$x\in D_{SL_n(\Cee)}$ if and only if $H^0(C; S) \neq 0$. In this case,
the vector bundle associated to $\xi$ via the $SO$-representation is
semistable. Suppose that $x\in D_{SL_n(\Cee)}$, and let $V$ be
the regular semistable vector bundle of rank $n$ associated to any point of
$\WP(SL_n(\Cee))$ mapping to $x$. Then $S$ is unstable if and only if
$\bigwedge ^kV$ is unstable for some $k$ with $2\leq k\leq n-2$. 
\qed
\end{proposition}

Of course, the regular semistable vector bundles $V$ such that $\bigwedge ^kV$
is unstable for some $k$ have been classified in Proposition~\ref{sln}. Thus,
for example, if
$n\leq 3$ then every bundle
$S$ associated to
$\xi$ from the spin representation is semistable.

To describe   $D_{SL_n(\Cee)}\subseteq \WP(Spin(2n+1))$ in terms of the
$SO$-representation, recall that we have the $\Cee^*$ commuting
with
$SL_n(\Cee)$, which then defines a filtration on any associated vector
bundle. In this case, using the standard description of the roots of
$Spin(2n+1)$ as $\alpha_1 = e_1-e_2, \dots, \alpha_{n-1} = e_{n-1}-e_n,
\alpha_n =e_n$, so that the highest weight of the spin representation is
$\frac12\sum_ie_i =\varpi_{\alpha_n}$,  the Lie algebra of the
$\Cee^*$ in question is
$\Cee\cdot\varpi_{\alpha_n}\spcheck=\Cee\cdot (\sum_ie_i)$. Since the weights of
the standard representation are $\pm e_i$ and the trivial weight, we see
that, if $\xi$ is a $Spin(2n+1)$-bundle arising from parabolic induction
from $SL_n(\Cee)$, then  the associated orthogonal vector bundle
$V=\xi\times_{Spin(2n+1)}\Cee^{2n+1} $  has a filtration
$F_{-1}\subseteq F_0\subseteq F_1 =V$, where $(F_{-1})^{\perp} = F_0$, such
that the successive quotients are
$V',\scrO_C, (V')^*$, where
$V'$ is a regular semistable vector bundle of rank $n$ and trivial 
determinant. Conversely, if $V$ is an  orthogonal bundle with such a
filtration, then it is not difficult to show that there is a lifting of $V$
to a $Spin(2n+1)$-bundle $\xi$ (uniquely, in general), such that the
associated spin bundle has a section. In particular, if $\xi$ arises from the
parabolic construction, then in fact it arises from parabolic induction
from $SL_n(\Cee)$.

We turn now to the quasi-minuscule representation of $Spin(2n+1)$, i.e.\ the
standard representation. Arguments similar to those above show:

\begin{proposition}\label{oddspin2} Suppose that $n\geq 2$. Given  $x\in
\WP(Spin(2n+1))$, let  $\xi$ be the $Spin(2n+1)$-bundle corresponding to $x$
and let
$V$ be the vector bundle associated to $\xi$ by the $SO$-representation. Then
$x\in D_{Spin(2n-1)}$ if and only if there is a self-dual complex
$$0 \to \scrO_C \to V \to \scrO_C \to 0.$$ In this case,
let
$V'$ be the orthogonal bundle of rank $2n-1$ associated to any point of
$\WP(Spin(2n-1))$ mapping to $x$ under parabolic induction, so that $V'$ is
the cohomology of the above complex. Then $V$ is unstable if and only if $V'$
is unstable.\qed
\end{proposition}

\begin{corollary}\label{oddspin3} For all $n\geq 1$, there is a unique point
$x\in
\WP(Spin(2n+1))$ such that, if $\xi$ is the corresponding $Spin(2n+1)$-bundle
and $V$ is the vector bundle of rank $2n+1$  associated to $\xi$ via the
standard representation, then $V$ is unstable. The point $x$ comes via
repeated parabolic induction from $Spin(3) = SL_2(\Cee)$ and  the
bundle corresponding to $V_{1,1}$ in case $C$ is nodal or $V_2$ in case $C$ is
cuspidal. Finally, the vector bundle
$S$ associated to 
$\xi$  via the spin  representation is semistable.
\end{corollary}
\begin{proof} By repeated applications of Proposition~\ref{oddspin2},
$x$ comes via iterated parabolic
induction from $Spin(3) = SL_2(\Cee)$. There is a unique point of
$\WP(SL_2(\Cee))$ which lies in $D_\infty$, and the corresponding
representation of $SL_2(\Cee)$ is the adjoint representation. Hence the
associated rank $3$ vector bundle is unstable. To see the last statement, note
that it follows from Proposition~\ref{oddspin} that, if $S$ were unstable, then
$V$ would be semistable. Thus conversely if $V$ is unstable then $S$ is
semistable.
\end{proof}

\medskip
\noindent\textbf{The case of $\boldsymbol{Spin(2n)}$}

In this case, there are three minuscule representations: the two half-spin
representations $\rho_{\pm}$ and the standard representation $\rho_0$. If $n$ 
is odd, then $\rho_+$ is the dual of $\rho_-$, whereas if $n$ is even then
$\rho_+$ and $\rho_-$ are both isomorphic to their own duals. We shall
consider first the half-spin representations and then the standard ($SO$)
representation. Note that there is always an outer automorphism of $Spin(2n)$
which exchanges $\rho_+$ and
$\rho_-$, so that it suffices to study one of them. However, if $n$ is odd,
$\xi$ is a  $Spin(2n)$-bundle, and $S^+$ and $S^-$ are the two vector bundles
associated to $\xi$ by $\rho_+$ and $\rho_-$ respectively, then $S^+$ has a
section if and only if $S^-$ has a section, whereas if $n$ is even this is no
longer necessarily the case.

Consider then $\rho_+$. The stabilizer of its highest weight is isomorphic to
$SL_n(\Cee)$, and $\rho_+|SL_n(\Cee)$ is the representation
$\bigwedge^{2*}\Cee^n$, where $\Cee^n$ is the standard representation. Likewise,
$\rho_-|SL_n(\Cee)$ is the representation
$\bigwedge^{2*-1}\Cee^n$. Finally the
restriction of the
$SO$-representation to $G'$ is $\Cee^n\oplus (\Cee^n)^*$. Thus:

\begin{proposition}\label{evenspin} Given  $x\in
\WP(Spin(2n))$, let  $\xi$ be the $Spin(2n)$-bundle corresponding to $\xi$
and let
$S^+$ be the vector bundle associated to $\xi$ by $\rho_+$. Then
$x\in D_{SL_n(\Cee)}$ if and only if $H^0(C; S^+) \neq 0$. In this case,
the vector bundle associated to $\xi$ via the standard representation is
semistable. Suppose that $x\in D_{SL_n(\Cee)}$, and let $V$
be the regular semistable vector bundle of rank $n$ associated to any point of
$\WP(SL_n(\Cee))$ mapping to $x$. Then $S^+$ is unstable if and only if
$\bigwedge ^{2k}V$ is unstable for some $k$. 
\qed
\end{proposition}

Note that, in the above notation, if $S^+$ is unstable, then so is $S^-$, 
except in the case where $n=4$ and $V$ is of the form $V_{2,2}$ or $V_4$. In
this case,
$S^-$ is semistable. Of course, the analogous construction for $\rho_-$ yields a
$Spin(8)$-bundle  for which $S^-$ is unstable and $S^+$ is semistable.

As in the case of $Spin(2n+1)$, if $V$ is an orthogonal bundle of rank $2n$
arising from a $Spin(2n)$-bundle which is in the image of $\WP(SL_n(\Cee))$,
then there is a self-annihilating subspace
$V'$ of $V$, where $V'$ is a regular semistable vector bundle of rank $n$ and
trivial determinant. Thus there is an exact sequence
$$0 \to V' \to V \to (V')^* \to 0,$$
with $(V')^{\perp} = V'$. Conversely, if $V$ is an orthogonal bundle with
such a filtration, the subspace $V'$ determines exactly one of the maximal
parabolic subgroups corresponding to the highest weights of $\rho_+$ and
$\rho_-$. If say it determines $\rho_+$, then there is a lifting of $V$ to a 
$Spin(2n)$-bundle $\xi$ (uniquely, in general) for which $S^+$ has a section,
and if in addition $\xi$ arises from the parabolic construction then it arises
via parabolic induction from $SL_n(\Cee)$.

Turning now to the case of the standard representation, the image of parabolic
induction from $Spin(2n-2)$ corresponds to those $\xi$ such that, if $V$ is the
vector bundle associated to the standard representation, then there is a
self-dual complex 
$$0 \to \scrO_C \to V \to \scrO_C \to 0,$$
and the corresponding $SO(2n-2)$-bundle is $(\scrO_C)^{\perp}/\scrO_C$. In this
case,  we can use repeated parabolic induction to reduce to the case of
$SL_4(\Cee)$, where there is a unique unstable bundle of the form
$\bigwedge^2V$, for $V$ rank $4$, regular semistable and with trivial
determinant. Thus:

\begin{proposition}\label{evenspin2} For $n\geq 3$, there is a unique point
$x\in
\WP(Spin(2n))$ such that, if $\xi$ is the corresponding $Spin(2n)$-bundle
and $V$ is the vector bundle of rank $2n$  associated to $\xi$ via the
$SO$-representation, then $V$ is unstable. The point $x$ comes via
repeated parabolic induction from $Spin(6) = SL_4(\Cee)$ and the bundle
described above. Finally, the vector bundle $S$ associated to 
$\xi$  via the spin  representation is semistable. \qed
\end{proposition}

In fact, Proposition~\ref{evenspin2} continues to hold for the group
$Spin(4)\cong SL_2(\Cee) \times SL_2(\Cee)$, where $\WP(Spin(4))\cong
\Pee^1\times \Pee^1$ and the unique unstable bundle in the
$SO$-representation is $V_{1,1}\otimes V_{1,1}$ if $C$ is nodal and
$V_2\otimes V_2$ if $C$ is cuspidal.

In particular, for $Spin(8)$, triality exchanges the three representations
$\rho_+, \rho_-, \rho_0$, and, for each of the three representations, there is a
unique point in $\WP(Spin(8))$ for which the vector bundle associated to
the given representation is unstable and the remaining two vector bundles are
semistable.

\medskip
\noindent\textbf{The case where $\boldsymbol{G}$ is of type $\boldsymbol{E_6}$}

There are two minuscule representations of $G$, each of dimension $27$. Since
they are exchanged under duality, it suffices to consider just one of them, say
$\rho$. The stabilizer of a weight of $\rho$ is $Spin(10)$. Moreover,
$\rho|Spin(10)$ decomposes into one trivial summand, one copy of the standard
$10$-dimensional representation of $Spin(10)$, and one of the spin
representations, of dimension $16$. If
$x\in
\WP(G)$ corresponds to the bundle $\xi$, and $V$ is the associated
$27$-dimensional vector bundle, then $x\in D_{Spin(10)}$ if and only if $h^0(C;
V) \neq 0$. If $V$ is unstable, then $x\in  D_{Spin(10)}$, and, if $\xi'$ is
the corresponding
$Spin(10)$-bundle, then $\xi'\times _G\Cee^{27} \cong \scrO_C \oplus V'\oplus
S$, where $V'$ is the rank $10$ vector bundle associated to $\xi'$ via the
standard representation of $Spin(10)$, and $S$ is the rank $16$-vector bundle
associated to $\xi'$ via the appropriate spin representation. Hence $V$ is
unstable if and only if either $S$ is unstable or $V'$ is unstable. Thus the
locus of points where $V$ is unstable is the image of a rational curve
$\Gamma$ in
$\WP(Spin(10))$, corresponding to the set of regular $SL_5(\Cee)$-bundles
of the form
$V_{2,2;\lambda}\oplus \lambda'$, plus at most one more point $x$, the image
under iterated parabolic induction from the point of $\WP(SL_4(\Cee))$
corresponding either to $V_{2,2}$ in the nodal case or $V_4$ in the cuspidal
case. However, one can check that the morphism from
$\WP(Spin(10))$ to
$\WP(G)$ identifies the   point $x$ with a point on   $\Gamma$.  The point
$x$ comes from parabolic induction from $A_3$ to $D_4$ and thence to $D_5$.
But we can also do parabolic induction from $A_3$ to $A_4$ to $D_5$. The two
embeddings of $A_3$ in $D_5$ are not conjugate, but one can check that they
become conjugate in $E_6$. Hence there is an irreducible rational curve in
$\WP(G)$ corresponding to  unstable bundles $\xi$.   The  
$G$-bundles $\xi$ are unstable in both minuscule representations of $G$ and
hence in all nontrivial representations of $G$.

\medskip
\noindent\textbf{The case where $\boldsymbol{G}$ is of type $\boldsymbol{E_7}$} 

If $\rho$ is the unique minuscule representation, then the stabilizer of a weight
of $\rho$ is a group $G_0$ of type $E_6$. The restriction $\rho|G_0$ is a
direct sum of the two minuscule representations of $G_0$ plus two copies of
the trivial representation. Thus, the
$G$-bundles $\xi$ arising from the parabolic construction for which the
corresponding vector bundle
$V$ of rank $56$ is unstable are the image via parabolic induction of the
corresponding locus for $G_0$, and hence are the image of the rational curve
$\Gamma$. (The morphism $\WP(G_0) \to \WP(G)$ is $2$-$1$ and the curve
$\Gamma \subseteq \WP(G_0)$ is invariant under the associated involution.) 
The corresponding
$G$-bundles are unstable in the minuscule representation of
$G$ and hence in all nontrivial representations of $G$.

\medskip
\noindent\textbf{The case where $\boldsymbol{G}$ is of type $\boldsymbol{F_4}$}

As a
$Spin(7)$-module, the $26$-dimensional quasi-minuscule representation $\rho$
of $G$ splits into   three copies of the trivial representation, one copy
of the standard $7$-dimensional representation of $Spin(7)$, and two copies
of the
$8$-dimensional spin representation. Now suppose that $\xi$ is a $G$-bundle,
arising from the parabolic construction, such that the  vector bundle $V$
associated to $\xi$ via $\rho$ is unstable. Then $\xi$ is obtained via
parabolic induction from a $Spin(7)$-bundle $\xi'$. The vector bundle
$V'$ associated to $\xi'$ via $\rho$ then decomposes in the same way as
$\rho|Spin(7)$. By Proposition~\ref{oddspin}, the two $8$-dimensional pieces are
semistable. Thus, the bundle $\xi'\times_{Spin(7)}\Cee^7$ obtained from the
standard representation of $Spin(7)$ is unstable. It follows from
Corollary~\ref{oddspin3} that there is a unique such $Spin(7)$-bundle. Hence
there is a unique $G$-bundle up to isomorphism arising from the parabolic
construction such that the associated $26$-dimensional vector bundle is
unstable. Thus it is unstable in all nontrivial representations of $G$.

\medskip
\noindent\textbf{The case where $\boldsymbol{G}$ is of type $\boldsymbol{G_2}$}

If $\xi$ is a $G_2$-bundle  arising from the parabolic construction, for which
the associated $7$-dimensional vector bundle $V$ is unstable, then $\xi$ comes
via parabolic induction from an $SL_2(\Cee)$-bundle. The quasi-minuscule
representation of
$G_2$ decomposes as an $SL_2(\Cee)$-module into three trivial summands and two
copies of the standard representation. It follows that, for every
$G_2$-bundle $\xi$ arising from the  parabolic construction, the 
vector bundle associated to $\xi$ via the $7$-dimensional
irreducible representation is semistable.

\bigskip
\noindent
Department of Mathematics \\
Columbia University \\
New York, NY 10027 \\
USA

\bigskip
\noindent
{\tt rf@math.columbia.edu, jm@math.columbia.edu}

\end{document}